\documentclass{article}

\usepackage{amsthm}

\usepackage{amsfonts}
\usepackage{graphicx}
\usepackage{amsmath}
\usepackage{amssymb}
\DeclareGraphicsExtensions{.pdf,.png,.jpg} 
\usepackage{epstopdf}
\usepackage{sectsty}

\subsubsectionfont{\large}

\newtheorem{teorema}{Theorem}[section]
\newtheorem{lema}[teorema]{Lemma}

\newtheorem{corolario}[teorema]{Corollary}
\newtheorem{proposicion}[teorema]{Proposition}
\newtheorem{pregunta}[teorema]{Question}

\newtheorem{asuncion}[teorema]{Assumption}

\newtheorem*{lema.sn}{Lemma}
\newtheorem*{def.sn}{Definition}
\newtheorem*{teoa}{Theorem A}
\newtheorem*{addb}{Addendum B}
\newtheorem*{teoc}{Theorem C}
\newtheorem*{claim.sn}{Claim}

\theoremstyle{definition}
\newtheorem{remark}[teorema]{Remark}

\theoremstyle{definition}
\newtheorem{definicion}[teorema]{Definition}

\theoremstyle{definition}
\newtheorem{claim}[teorema]{Claim}

\newenvironment{prueba}[1][Proof.]{\begin{trivlist}
   \item[\hskip \labelsep {\itshape #1}]}{\rule{0.5em}{0.5em}\end{trivlist}}

\newcommand{\ol}{\overline}
\newcommand{\wt}{\widetilde}
\newcommand{\wh}{\widehat}
\newcommand{\ra}{\rightarrow}
\newcommand{\pr}{\partial}
\newcommand{\txt}{\textnormal}
\newcommand{\tl}{\tilde}
\newcommand{\empt}{\emptyset}
\newcommand{\e}{\epsilon}
\newcommand{\minus}{\setminus}

\newcommand{\R}{\mathbf R}

\newcommand{\N}{\mathbf N}
\newcommand{\T}{\mathbf T}
\newcommand{\Z}{\mathbf Z}
\newcommand{\Q}{\mathbf Q}
\newcommand{\cl}{\mathcal}

\title{On torus homeomorphisms whose rotation set is an interval}

\author{Pablo D\'avalos}

\date{}

\begin{document}

\maketitle


\begin{abstract}
We prove that for a homeomorphism $\tl{f}:\T^2\ra\T^2$ in the homotopy class of the identity and with a lift $f:\R^2\ra\R^2$ whose rotation set $\rho(f)$ is an interval, either every rational point in $\rho(f)$ is realized by a periodic orbit, or the dynamics of $\tl{f}$ is \textit{annular}, in the sense that there exists a periodic, essential, annular set for $\tl{f}$. In the latter case we also give a qualitative description of the dynamics. 
\end{abstract}

\tableofcontents

\section{Introduction.}

In \cite{po}, Poincar\'e defined the rotation number for circle homeomorphisms, and he showed it to be a topological invariant carrying dynamical information. For a circle homeomorphism $\tl{f}$ with a lift $f:\R\ra\R$, the rotation number of $f$, denoted $\rho(f)$, is rational if an only if $f$ has periodic points, and is irrational if and only if there exists a model for the dynamics of $\tl{f}$, in the sense that $\tl{f}$ is semi-conjugate to the irrational rotation $x\mapsto x+\rho(f) \ \mod 1$.

In \cite{mz} Misiurewicz and Ziemian generalize the concept of the rotation number for homeomorphisms of $\T^n$, for any $n\in\N$. For a torus homeomorphism $\tl{f}:\T^n\ra\T^n$, the \textit{rotation set} of some lift $f:\R^n\ra\R^n$, denoted $\rho(f)$, is defined as the set of accumulation points of sequences of the form 
$$\left\{ \frac{f^{m_i}(x_i)-x_i}{m_i} \right\}_{i\in\N}$$
where $m_i\ra\infty$ and $x_i\in\R^n$. The set $\rho(f)$ is a compact subset of $\R^n$, and in the case $n=2$ it is also convex. In \cite{mz}, and in many other subsequent articles it is studied the relation between the rotation set and the dynamics of $\tl{f}$ (see for example \cite{mz2}, \cite{lm}, \cite{f1},\cite{f2},\cite{kk}, \cite{j1}, etc.). A lot more is known in the case that $n=2$ thanks to the theory of surface homeomorphisms, like Brouwer theory, Thurston's classification theory, etc. For this reason, we will restrict ourselves to the case $n=2$.

A basic question, making an analogy with the theory of the circle, is whether there are periodic points associated to points with rational coordinates in $\rho(f)$. This problem has been extensively studied. For a point $v\in\rho(f)\cap\Q^2$, expressed in the form $v=(p_1/q,p_2/q)$ with gcd$(p_1,p_2,q)=1$, we say that $v$ is realized by a periodic orbit of $\tl{f}$ if there exists $x\in\R^2$ such that
$$f^q(x)= x + (p_1,p_2).$$
In \cite{f1}, Franks proved that rational extremal points in the rotation set are realized by periodic orbits, and in \cite{f2} he proved that rational points in the interior\footnote{with respect to the topology of $\R^2$.} of the rotation set are also realized by periodic orbits. In the case that $\rho(f)$ is an interval, it is not true that rational points are always realized by periodic orbits. However, in \cite{f3} and \cite{kk}, are given sufficient conditions under which, if $\rho(f)$ is an interval, every rational point in $\rho(f)$ is realized by a periodic orbit.

A second basic question is whether there are dynamical models associated to certain rotation sets. Unfortunately, this is not always the case. An example of this are the \textit{pseudo-rotations}, that is, homeomorphisms whose rotation set is reduced to a single point, called \textit{rotation vector}. There are examples of pseudo-rotations with the same rotation vector, but with very different behavior. The simplest example of a pseudo-rotation is a translation $T_v:x\mapsto x+v \ \mod \Z^2$. Unlike the case of the circle, one may have that the deviations
$$D(x,n)=\left| f^n(x)-x- nv \right|$$
are unbounded, and this allows to create examples of pseudorotations with exotic properties, like Lebesgue weak-mixing \cite{fay1}, topological expansive-type properties \cite{kk2}, etc. Therefore, in the case that the rotation set is reduced to a point, there seems not to be models for the dynamics associated to the rotation set. 

In this work, we study the case that the rotation set is an interval. Suppose that $\tl{f}$ is a torus homeomorphism with a lift $f:\R^2\ra\R^2$ whose rotation set is an interval of the form $\{0\}\times I$, with $0\in\txt{int}\, I$ (the simplest example of such a homeomorphism is the twist $(x,y)\mapsto (x,y+\sin(2\pi x)) \mod \Z^2$). As above, we could wonder if the horizontal deviations
$$D_1(x,n)=\left| f^n(x)_1-x_1 \right|$$
can be unbounded, and in this way construct examples with qualitatively different dynamics. We will show that if $\rho(f)=\{0\}\times I$, then, either every rational point in $\rho(f)$ is realized by a periodic orbit, or deviations $D_1(x,n)$ are uniformly bounded in $x$ and $n$. With this, we will obtain the following, for the case that $\rho(f)$ is any interval: 
\begin{center}
\textit{If $\rho(f)$ is an interval, then either every rational point in $\rho(f)$ is realized by a periodic orbit, or there exists a `model' for the dynamics.}
\end{center}

\begin{figure}[h]        
\begin{center} 
\includegraphics{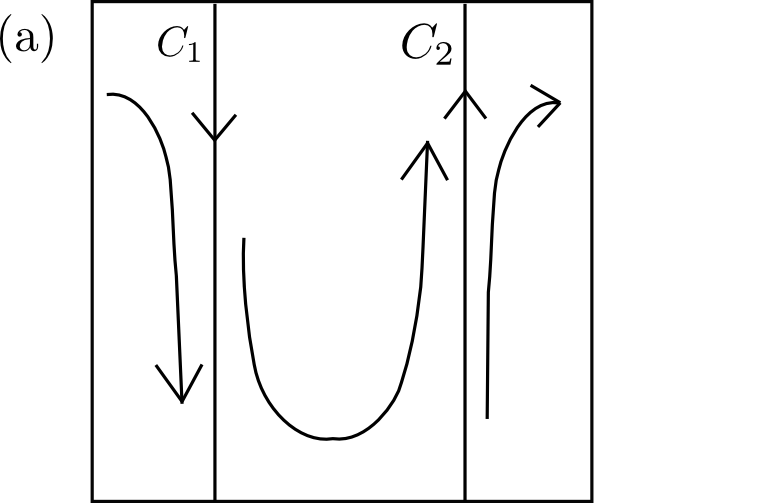}
\includegraphics{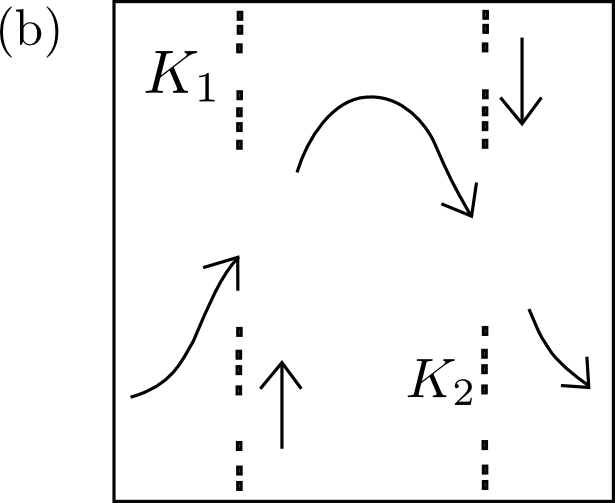}
\caption{$(a)$ $\tl{f}(x,y)=(\varphi(x),y+\pi\sin(2\pi x))$, with $\varphi:S^1\ra S^1$ Morse-Smale. $(b)$ a tentative to obtain unbounded horizontal displacements $D_1(x,n)$.}
\label{fig.ex1}
\end{center}  
\end{figure}

A precise meaning of a `model' for the dynamics is given in Theorems A and B. The simplest example of a homeomorphism with rational points in the rotation set not realized by a periodic orbit is a skew product $\tl{f}$ of a Morse-Smale circle homeomorphism (with fixed points), and a twist torus homeomorphism, as illustrated in Fig. \ref{fig.ex1}a. This example has a lift $f:\R^2\ra\R^2$ with $\rho(f)=[-\pi,\pi]$, and any rational point contained in $\rho(f)$ is not realized by a periodic orbit, since $\tl{f}$ has no periodic points.

One could try to modify this example, replacing the two invariant circles $C_1$ and $C_2$ by two non-connected invariant sets $K_1$ and $K_2$, so that orbits can pass through in hopes to obtain unbounded horizontal displacements $D_1(x,n)$, \textit{but still having $\rho(f)=\{0\}\times I$} (see Fig. \ref{fig.ex1}b). However, we will see that this is not possible; indeed, our main theorem implies that such a modification always leads to $\max|\txt{pr}_1(\rho(f))|>0$.

In Theorem A we deal with the particular case that $\rho(f)$ is a vertical interval containing the origin in its interior, and the origin is not realized by a periodic orbit, and in Theorem C we deal with the case that $\rho(f)$ is a general interval. In Theorem A we prove that the horizontal displacements $D_1(x,n)$ are uniformly bounded, and to prove this we will show there exists an invariant vertical `wall', that is, an invariant, annular, essential and vertical set for $\tl{f}$. By an \textit{annular} set we mean a nested intersection of compact annuli $A_i\subset\T^2$ such that the inclusion $A_{i+1}\hookrightarrow A_i$ is a homotopy equivalence, and by \textit{essential} and \textit{vertical} we mean that the annuli $A_i$ are homotopic to the annulus $\{x\in\T^2\: : \: 0\leq x_1\leq 1/2\}$. This `wall' will also have the property of being a \textit{semi-attractor}, which we now define:
\begin{definicion}
An annular, essential set $A\subset\T^2$ is a \textit{semi-attractor} for a homeomorphism $h:\T^2\ra\T^2$ if $A$ is $h$-invariant and there exist two simple, closed, essential curves $\gamma_1,\gamma_2\subset\T^2$ disjoint from $A$ and such that:
\begin{itemize}
\item $\omega(x,h)\subset A$ for all $x\in\gamma_1$, and
\item either $\omega(y,h)\subset A$ for all $y\in\gamma_2$ or $\alpha(y,h)\subset A$ for all $y\in\gamma_2$
\end{itemize}
\end{definicion}

We say that a curve $\gamma\subset\T^2$ is \textit{free forever} for $\tl{f}$ if $\tl{f}^n(\gamma)\cap\gamma=\empt$ for all $n\in\Z$. A closed curve $\gamma\subset\T^2$ is \textit{vertical} if it is homotopic to a straight vertical circle. Also, if $\gamma_1,\gamma_2\subset\T^2$ are vertical and disjoint curves, $[\gamma_1,\gamma_2]\subset\T^2$ denotes the closed annulus whose `left' border component is $\gamma_1$ and whose `right' border component is $\gamma_2$ (for precise definitions see Section \ref{sec.notations}). By last, we denote by $\Omega(\tl{f})$ the \textit{non-wandering set} of $\tl{f}$, that is, the set of points $x\in\T^2$ such that for every neighborhood $V$ of $x$, there is $n>0$ such that $\tl{f}^n(V)\cap V\neq\empt$.   

We now state our main theorem.

\begin{figure}[h]        
\begin{center} 
\includegraphics{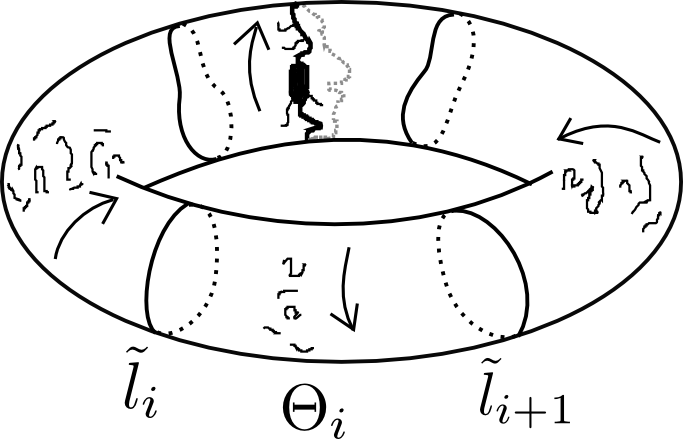}
\caption{The sets $\Theta_i$ and the curves $\tl{l}_i$. At least one of the $\Theta_i$ must be annular and essential.}
\label{fig.teo1}
\end{center}  
\end{figure}

\begin{teoa} \label{teo1}
Let $\tl{f}$ be a homeomorphism of $\T^2$ homotopic to the identity with a lift $f:\R^2\ra\R^2$ such that: 
\begin{itemize} 
\item $\rho(f)=\{0\}\times I$, where $I$ is a non-degenerate interval containing $0$ in its interior, and
\item $(0,0)\in\rho(f)$ is not realized by a periodic orbit.
\end{itemize}
Then, there exists a finite family $\{ \tl{l}_i \}_{i=0}^{r-1}$, $r\geq 2$, of curves in $\T^2$ which are simple, closed, vertical, and pairwise dijoint, and with the following properties. If 
$$\Theta_i:= \bigcap_{n\in\Z} \tl{f}^n \left( [\tl{l}_i,\tl{l}_{i+1}] \right)\ \ \ \txt{for } i\in\Z/r\Z,$$
then,
\begin{enumerate}
\item at least one of the sets $\Theta_i$ is an annular, essential, $\tl{f}$-invariant set which is a semi-attractor, 
\item the curves $\tl{l}_0,\tl{l}_1,\ldots,\tl{l}_{r-1}$ are free forever for $\tl{f}$,
\item there is $\epsilon>0$ such that $\rho(\Theta_i,f)$ is contained either in $\{0\}\times (\epsilon,\infty)$, or in $\{0\}\times (-\infty,-\epsilon)$, and
\item $\Omega(\tl{f})\subset \cup \Theta_i$, (see Fig. \ref{fig.teo1}).
\end{enumerate}
\end{teoa}

\begin{remark}
If the cardinality $r$ of the family $\{\tl{l}_i\}_{i=0}^{r-1}$ is minimal such that it satisfies the conclusion of Theorem A, then the sets $\Theta_i$ rotate alternatively `upwards' and `downwards'. 
\end{remark}

For a definition of the rotation sets $\rho(\Theta_i,f)$ see Section \ref{sec.conj.rot}. In Theorem A, the curves $\tl{l}_i$ decompose the dynamics in a way similar to a filtration. If a point, under iteration by $\tl{f}$, leaves an annulus $[\tl{l}_i,\tl{l}_{i+1}]$, it never enters that annulus again, and if a point enters an annulus $[\tl{l}_{i_0},\tl{l}_{i_0+1}]$ containing an essential set $\Theta_{i_0}$, then it remains in that annulus forever.

Furthermore, we have that the dynamical decomposition given by Theorem A is essentially unique:

\begin{addb}
Let $\tl{f}$ and $f$ be as in Theorem A. Let $\{l_i\}_{i=0}^{r-1}$, and $\{l'_i\}_{i=0}^{r'-1}$ be two families of simple, closed, vertical and pairwise disjoint curves in $\T^2$, and suppose that both families satisfy the conslusion of Theorem A. Suppose also that the cardinalities $r$ and $r'$ are minimal with respect to this property.

Then, if we denote
$$\Theta_i = \bigcap_{n\in\Z} \tl{f}^n \left( [l_i,l_{i+1}] \right) \  \ \txt{for } i\in\Z/r\Z, \txt{ and}$$
$$\Theta'_i = \bigcap_{n\in\Z} \tl{f}^n \left( [l'_i,l'_{i+1}] \right)\ \  \txt{for } i\in\Z/r'\Z,$$
we have that $\{\Theta_i\}_{i=0}^{r'-1}= \{\Theta'_i\}_{i=0}^{r-1}$.
\end{addb}

As a corollary of Theorem A and Addendum C we will obtain the following theorem, dealing with the case that the rotation set is a general interval.

\begin{teoc} 
Let $\tl{f}:\T^2\ra\T^2$ be a homeomorphism homotopic to the identity with a lift $f:\R^2\ra\R^2$ whose rotation set is an interval. 

Then, either every rational point in the rotation set is realized by a periodic orbit, or there is $k\in\N$ such that $\tl{f}^k$ is topologically conjugate to a homeomorphism within the hypotheses of Theorem A. In the latter case:
\begin{itemize}

\item there exists a simple, closed, essential curve in $\T^2$ that is free forever for $\tl{f}$,

\item there exists an annular, essential set $A$ in $\T^2$ that is periodic for $\tl{f}$, and 

\item if $q$ is the minimal period of $A$, the sets $A$, $\tl{f}(A)$, \ldots, $\tl{f}^{q-1}(A)$ are pairwise disjoint.

\end{itemize}
\end{teoc}

\begin{remark} 
We would like to emphasize that the main novelty in Theorem A is the existence of an \textit{annular essential} set, which is invariant for $\tl{f}$. The dynamical decomposition given by Theorem A without the existence of an annular, essential, $\tl{f}$-invariant set is obtained in Section \ref{sec.lema.prel}, and it combines the techniques from a theorem of Le Calvez (Theorem \ref{teo.intro.pat}) with standard techniques from ergodic theory, like Atkinson's Lemma \ref{atk}. The existence of a curve that is free forever for $\tl{f}$ is also a novelty, but by the mentioned techniques it is equivalent to the existence of the annular essential $\tl{f}$-invariant set.  
\end{remark}

\begin{remark}
In Theorem A from \cite{kk} it is proved that, if $\rho(f)$ is an interval and there is a rational point in $\rho(f)$ that is not realized by a periodic orbit, then for all $n\in\N$ there is a simple, closed, essential curve disjoint from its first $n$ iterates. Theorem C in this article generalizes that result, showing that there is actually a (simple, closed, essential) curve that is free forever.
\end{remark}

We do not know if the property of the displacements $D_1(x,n)$ being uniformly bounded is also present in the case that all the rational points in the rotation set are realized by periodic orbits:

\begin{pregunta}
If $\tl{f}:\T^2\ra\T^2$ is a homeomorphism with a lift $f:\R^2\ra\R^2$ such that $\rho(f)$ is an interval of the form $\{0\}\times I$, then, are the deviations $D_1(x,n)$ uniformly bounded?
\end{pregunta}


This work is organized as follows. In section \ref{sec.preliminares} we introduce the preliminary theory used in the proof of Theorem A, which is mainly the following: the rotation set for homeomorphisms of $\T^2$ and some results related to it, the Brouwer theory for planar homeomorphisms developed by Patrice Le Calvez, and Atkinson's Lemma from ergodic theory. In section \ref{sec.teofromteo1} we prove Theorem C assuming Theorem A. The proof of Theorem A is divided in sections \ref{sec.lema.prel} and \ref{sec.teoa2}, and the proof of Addendum B is given in section \ref{sec.addendum}.\\
\\
\textbf{Acknowledgements.} This work is from a Ph.D. thesis under the supervision of A. Koropecki and E. Pujals. I thank A. Koropecki for proposing me to study this problem, for many hours of discussions, for important ideas and simplifying many proofs. I also thank E. Pujals for key comments and ideas, and for his constant support and encouragement. I am grateful to B. Fayad for his hospitality during a visit to the Universit\'e Paris VII and for insightful conversations, and to P. Le Calvez also for important comments. Finally, I would like to thank an anonymous referee for suggesting to incorporate Addendum B, which is used to correct an error in the original proof of Theorem C, for suggesting also a restructuration of the proof of Main Lemma \ref{N6}, which simplifies the original one, and for a careful reading of this paper.

\section{Notations.}   \label{sec.notations}

By $\txt{pr}_1,\txt{pr}_2:\R^2\ra\R$, we will denote the projections to the first and second coordinate, respectively. Also, if $x\in\R^2$, $x_1$ and $x_2$ will denote $\txt{pr}_1(x)$ and $\txt{pr}_2(x)$, respectively. 

For the circle $S^1=\R/\Z$, and the two-torus $\T^2=\R^2/\Z^2$, denote by $\pi,\pi'$ and $\pi''$ the canonical projections 
$$\R^2 \stackrel{\pi}{\ra} \R\times S^1 \stackrel{\pi''}{\ra} \T^2, \ \ \ \txt{and } \ \  \pi'= \pi''\circ \pi.$$

For a set $A\subset\R$, the diameter of $A$ is $\txt{diam}(A)=\sup_{x,y\in A} |x-y|$. For $A\subset\R^2$, the \textbf{horizontal diameter} of $A$ is $\txt{diam}_1(A)= \txt{diam}(\txt{pr}_1(A))$, and the \textbf{vertical diameter} of $A$ is $\txt{diam}_2(A)=\txt{diam}(\txt{pr}_2(A))$. 

For a set $A\subset\R^2$ and $x\in\R^2$, denote $d(x,A)=\inf_{y\in A}|y-x|$. For $x\in\R^2$ and $r>0$, denote $B_r(x)=\{y\in\R^2\, : \, |y-x|<r\}$, and for $A\subset\R^2$, denote $B_r(A)=\{x\in\R^2\, : \, d(x,A)<r\}$.

We will denote also by $d(\cdot,\cdot)$ the metric in $\T^2$ or in $\R\times S^1$ induced by the euclidean metric in $\R^2$. 

Define $T_1,T_2:\R^2\ra\R^2$ to be the translations $T_1:(x_1,x_2)\mapsto (x_1+1, x_2)$, $T_2:(x_1,x_2)\mapsto (x_1,x_2+1)$. Also, $T_1$ will denote the translation in $\R\times S^1$, $T_1:(x_1,x_2)\mapsto (x_1+1, x_2)$. 

By a \textbf{curve} $\gamma:I\ra\R^2$, depending on the context, we mean either $\gamma$ or $\txt{Im}(\gamma)\subset\R^2$. By an \textbf{arc}, we mean a compact injective curve, and if $\alpha$ is an arc, $\dot{\alpha}$ denotes the curve $\alpha$ without its endpoints. A \textbf{line} $\ell$ is a proper embedding $\ell:\R\ra\R^2$. By Shoenflies' Theorem (\cite{cairns}), given a line $\ell$ there exists an orientation preserving homeomorphism $h$ of $\R^2$ such that $h\circ \ell(t)=(0,t)$, for all $t\in\R$. Then, the open half-plane $h^{-1}((0,\infty)\times\R)$ is independent of $h$, and we call it the \textbf{right} of $\ell$, and denote it by $R(\ell)$. Analogously, we define $L(\ell)=h^{-1}((-\infty,0)\times\R)$ the open half-plane to the \textbf{left} of $\ell$.  The sets $\ol{R}(\ell)$ and $\ol{L}(\ell)$ denote the closures of $R(\ell)$ and $L(\ell)$, respectively. 

Let $\gamma:S^1\ra\R^2$ be a closed curve, and $x \in \R^2 \minus \gamma$. By Ind$(\gamma, x)$ we will denote the degree of the map $S^1\ra S^1$ given by 
$t\mapsto \gamma(t)-x / \| \gamma(t)-x \|$, and it will be called the \textbf{index} of $\gamma$ with respect to $x$.

By $\ell\prec\ell'$ we will mean $\ell\subset L(\ell')$. 

A closed curve $\gamma$ in $\T^2$ or in $\R\times S^1$ is \textbf{essential} if it is not homotopic to a point, and we say that $\gamma$ is \textbf{vertical} if $\gamma$ is freely homotopic to a curve of the form $c\beta$, where $c\in\{1,-1\}$ and $\beta(t)=(0,t)$. A curve $\gamma$ in $\T^2$ (or in $\R\times S^1$) is \textbf{free} for $f:\T^2\ra\T^2$ ($f:\R\times S^1\ra\R\times S^1$) if it is simple and closed, and $f(\gamma)\cap\gamma=\emptyset$, and we say it is \textbf{free forever} for $f$ if $\gamma$ is disjoint from all its iterates by $f$. 

If $\ell,\ell'$ are two lines in $\R^2$, we define $(\ell,\ell')=R(\ell)\cap L(\ell')$, and $[\ell,\ell']=\ol{R}(\ell)\cap \ol{L}(\ell)$. Similarly we define $(\ell,\ell'] = R(\ell)\cap \ol{L}(\ell')$ and $[\ell,\ell')=\ol{R}(\ell)\cap L(\ell')$. If $\gamma$ and $\gamma'$ are two disjoint, simple, closed and vertical curves in $\T^2$, we define the topological annuli $(\gamma,\gamma')\subset\T^2$ and $[\gamma,\gamma']\subset\T^2$ in the following way. Let $\tl{\gamma}\subset\R^2$ be any lift of $\gamma$, and let $\tl{\gamma}'$ be the first lift of $\gamma'$ to the right of $\tl{\gamma}$, that is, $\tl{\gamma}'$ is the lift of $\gamma'$ with $\tl{\gamma}\prec \tl{\gamma}'\prec T_1(\tl{\gamma})$. Orient $\tl{\gamma}$ and $\tl{\gamma}'$ as going upwards. Define then $(\gamma,\gamma')=\pi'((\tl{\gamma},\tl{\gamma}'))$ and $[\gamma,\gamma']=\pi'([\tl{\gamma},\tl{\gamma}'])$. In a similar way, if $\gamma$ and $\gamma'$ are disjoint, simple, closed and vertical curves in $\R\times S^1$, we define $(\gamma,\gamma')\subset \R\times S^1$ and $[\gamma,\gamma']\subset\R\times S^1$.     

Let $M$ be either $\T^2$ or in $\R\times S^1$. A topological (open or closed) annulus $B$ contained in $M$ is \textbf{essential} if the inclusion $B\hookrightarrow M$ induces a non trivial map $\pi_1(B)\ra\pi_1(M)$. A set $A\subset M$ is \textbf{annular} if it is a nested intersection of topological compact annuli $A_i$, such that for each $i$, the inclusion $A_{i+1}\hookrightarrow A_i$ is a homotopy equivalence. An annular set is said to be \textbf{essential} if the $A_i$ are essential annuli, and $A$ is called \textbf{vertical} if the $A_i$ are homotopic to the vertical annulus $\{0\leq x\leq 1/2\}\times S^1\subset\T^2$.

For a map $f:X\ra X$, where $X$ is any metric space, we define an $\epsilon$-chain for $f$ as a finite sequence $(x_i)_{i=0}^{i_1}\subset X$, with $i_1 \geq 2$, such that $d(x_{i+1},f(x_i))<\epsilon$ for all $0\leq i<i_1$. An $\epsilon$-chain $\{x_i\}_{i=0}^{i_1}$ is \textbf{periodic} if $x_{0}=x_{i_1}$. A point   $x\in X$ is \textbf{chain recurrent} for $f$ if for all $\epsilon>0$ there exists a periodic $\epsilon$-chain $\{x_i\}_{i=0}^n$ for $f$ with $x_0=x_n=x$. The \textbf{chain recurrent set}, denoted by $CR(f)$, is the set of chain recurrent points for $f$.

\section{Preliminaries.}               \label{sec.preliminares}

\subsection{The rotation set.}    \label{sec.conj.rot}

Denote by $\txt{Homeo}(\T^2)$ the set of homeomorphisms of $\T^2$, and by $\txt{Homeo}_*(\T^2)$ the elements of $\txt{Homeo}(\T^2)$ which are homotopic to the identity. Let $\tl{f}\in\txt{Homeo}_*(\T^2)$ and let $f:\R^2\ra\R^2$ be a lift of $\tl{f}$.

\begin{definicion}[\cite{mz}]     \label{def.rotset}
The \textit{rotation set of $f$} is defined as 
$$\rho(f)=\bigcap_{m=1}^{\infty} \txt{cl} \left( \bigcup_{n=m}^{\infty} \left\{ \frac{f^n(x)-x}{n}\, : \, x\in\R^2 \right\} \right) \subset\R^2. $$
The \textit{rotation set of a point} $x\in\R^2$ is defined by
$$\rho(x,f)= \bigcap _{m=1}^{\infty} \txt{cl} \left\{ \frac{f^n(x)-x}{n}\, : \, n>m \right\}.$$
If the above set consists of a single point $v\in\R^2$, we call $v$ the \textit{rotation vector of} $x$. If $\Lambda\subset\T^2$ is a compact $\tl{f}$-invariant set, we define the \textit{rotation set of} $\Lambda$ as
$$\rho(\Lambda,f)=\bigcap_{m=1}^{\infty} \txt{cl} \left( \bigcup_{n=m}^{\infty} \left\{ \frac{f^n(x)-x}{n}\, : \, x\in\pi'^{-1}(\Lambda) \right\} \right) \subset\R^2. $$
\end{definicion}

\begin{remark}      \label{remark.compact}
One can easily verify that the sets $\rho(f)$ and $\rho(\Lambda,f)$ in Definition \ref{def.rotset} are compact.
\end{remark}

\begin{remark}            \label{r.pot}
It is easy to see that for integers $n,m_1,m_2$,
$$\rho(T_1^{m_1}T_2^{m_2} f^n)= n\rho(f) + (m_1,m_2).$$
Then, the rotation set of any other lift of $\tl{f}$ is an integer translate of $\rho(f)$, and we can talk of the `rotation set of $\tl{f}$' if we keep in mind that it is defined modulo $\Z^2$. 
\end{remark}

\begin{teorema}[\cite{mz}]     \label{mz.cc}
Let $\tl{f}:\T^2\ra\T^2$ be a homeomorphism, let $\Lambda\subset\T^2$ be a compact $\tl{f}$-invariant set, and let $f:\R^2\ra\R^2$ be a lift of $\tl{f}$. Then the rotation set set $\rho(f)$ is a compact and convex, and every extremal point of $\rho(f)$ is the rotation vector of some point.
\end{teorema}

Given $A\in\txt{GL}(2,\Z)$, we denote by $\tl{A}$ the homeomorphism of $\T^2$ lifted by $A$. If $\tl{h}\in\txt{Homeo}(\T^2)$, there is a unique $A\in\txt{GL}(2,\Z)$ such that for every lift $h$ of $\tl{h}$, the map $h-A$ is bounded (in fact, $\Z^2$-periodic). Then $\tl{h}$ is isotopic to $\tl{A}$.

\begin{lema}
Let $\tl{f}\in{Homeo}_*(\T^2)$, $A\in\txt{GL}(2,\Z)$ and $\tl{h}\in\txt{Homeo}(\T^2)$ isotopic to $A$. Let $f$ and $h$ be lifts of $\tl{f}$ and $\tl{h}$ to $\R^2$. Then
$$\rho(hfh^{-1})=A\rho(f).$$
In particular, $\rho(AfA^{-1})=A\rho(f)$. 
\end{lema}

For a proof of this lemma, see for example \cite{kk}. 

\begin{remark}      \label{r.mat}
If $\rho(f)$ is segment of rational slope, there exists $A\in\txt{GL}(2,\Z)$ such that $A\rho(f)$ is a vertical segment. Indeed, if $\rho(f)$ is a segment of slope $p/q$ (with $p$ and $q$ coprime integers), we can find $x,y\in\Z$ such that $px+qy=1$, and letting
$$A=
\left(
\begin{array}{cc}
p & -q \\
y & x
\end{array} 
\right)$$
we have that $\txt{det}(A)=1$, and since $A(q,p)=(0,1)$, $A\rho(f)$ is vertical.
\end{remark}

\subsubsection{The rotation set and periodic orbits.}

Recall that we say that a rational point $(p_1/q,p_2/q)\in\rho(f)$ (with gcd$(p_1,p_2,q)=1$) is realized by a periodic orbit if there exists $x\in\R^2$ such that 
$$f^q(x)=x+(p_1,p_2).$$
We mention the following realization results. 

\begin{teorema}[\cite{f1}]
If a rational point of $\rho(f)$ is extremal, then it is realized by a periodic orbit.
\end{teorema}

\begin{teorema}[\cite{f2}]
Any rational point in the interior of $\rho(f)$ is realized by a periodic orbit.
\end{teorema}

The following theorem is stated for diffemorphisms in \cite{lc3}, p. 106, but its proof remains valid for homeomorphisms using the results in \cite{lc2} (see p. 9 of that article).

\begin{teorema}      \label{patirrac}
If a rational point belongs to a line of irrational slope which bounds a closed half-plane that contains $\rho(f)$, then this point is realized by a periodic orbit. 
\end{teorema}

\subsubsection{The rotation set and invariant measures.}            \label{sec.crmed}

For a compact $\tl{f}$-invariant set $\Lambda\subset\T^2$, we denote by $\cl{M}_{\tl{f}}(\Lambda)$ the family of $\tl{f}$-invariant probability measures with support in $\Lambda$, and $\cl{M}_{\tl{f}}=\cl{M}_{\tl{f}}(\T^2)$. Define the \textit{displacement function} $\phi:\T^2\ra\R^2$ by
$$\phi(\tl{x})= f(x)-x, \ \ \ \txt{for $x\in\pi'^{-1}(\tl{x})$}.$$
This is well defined, as any two preimages of $\tl{x}$ by the projection $\pi':\R^2\ra\T^2$ differ by an element of $\Z^2$, and $f$ is $\Z^2$-periodic. Now, for $\mu\in\cl{M}_{\tl{f}}$, we define the \textit{rotation vector} of $\mu$ as
$$\rho(\mu,f)=\int \phi \, d\mu.$$
Then, we define the sets
$$\rho_{mes}(\Lambda,f)=\left\{ \rho(\mu,f)\, : \, \mu\in \cl{M}_{\tl{f}}(\Lambda)\right\},$$
and
$$\rho_{erg}(\Lambda,f)=\left\{ \rho(\mu) \, : \, \mu \txt{ is ergodic for $\tl{f}$ and $\txt{supp}(\mu)\subset\Lambda$}  \right\}.$$
When $\Lambda=\T^2$ we simply write $\rho_{mes}(f)$ and $\rho_{erg}(f)$.

\begin{proposicion}[\cite{mz}]
It holds the following:
$$\rho(f)=\rho_{mes}(f)=\txt{conv} (\rho_{erg}(f)).$$
\end{proposicion}

When $\Lambda$ is a proper (compact, invariant) subset of $\T^2$, the set $\rho(\Lambda,f)$ is not necessarily convex (nor even connected). However, we have the following result from the folklore, which follows from arguments analogous to those in \cite{mz}.

\begin{proposicion} \label{erg}
It holds
$$\txt{conv}\rho(\Lambda,f)=\rho_{mes}(\Lambda,f),$$
and therefore, if $v\in\R^2$ is an extremal point of $\txt{conv} \, \rho(\Lambda,f)$, there exists an ergodic measure $\mu$ for $\tl{f}$ with $\rho(\mu,f)=v$ and $\txt{supp}(\mu)\subset\Lambda$.
\end{proposicion}

\subsection{Brouwer Theory.}     \label{sec.brth}

In \cite{br}, Brouwer proved the following theorem for homeomorphisms of the plane, known as the Brouwer Translation Theorem:

\begin{teorema}    \label{teoremabrouwer}
Let $h$ be an orientation preserving homeomorphism of $\R^2$ without fixed points. Then:
\begin{enumerate}
\item For all point $x\in\R^2$ there exists a line $\ell$ passing through $x$ such that 
$$\ell \prec h(\ell) \ \ \ \txt{and} \ \ \ h^{-1}(\ell)\prec\ell.$$
\item There exists a cover of $\R^2$ by open invariant disks where $h$ is conjugate to a translation.  
\end{enumerate}
\end{teorema}

A line satisfying item $(1)$ is called a \textit{Brouwer line} for $h$. By item $(2)$ we have that $h$ has no periodic points, and moreover, every point is wandering for $h$. The proofs of this theorem use the Brouwer Translation Lemma, which states that if an orientation preserving homeomorphism of the plane has no fixed points, then it has no periodic points. In \cite{f1} Franks proved the following stronger property of non-recurrence:

\begin{teorema}[Franks' Lemma]
Let $h:\R^2\ra\R^2$ be an orientation preserving homeomorphism. If there exist a sequence $(U_i)_{i\in \Z/n\Z}$ of pairwise disjoint open disks and a sequence of positive integers $(m_i)_{i\in\Z/n\Z}$ such that: 
\begin{itemize}
\item $h(U_i)\cap U_i=\empt \ \ \ $  $\forall\, i\in\Z/n\Z$, and,
\item $h^{m_i}(U_i)\cap U_{i+1}\neq\emptyset \ \ \ $  $\forall\, i\in\Z/n\Z$,
\end{itemize}
then $h$ has a fixed point.
\end{teorema}

As a corollary one obtains the following.

\begin{teorema}[\cite{f2}]     \label{franks}
Let $h$ be an orientation preserving homeomorphism of the plane, without fixed points and which is the lift of a homeomorphism of $\T^2$. Then, there exists $\epsilon>0$ such that there are no periodic $\epsilon$-chains for $h$. 
\end{teorema}

In \cite{lc1}, Le Calvez showed the following remarkable and much stronger version of the Brouwer Translation Theorem. 

\begin{teorema}  \label{pat1}
Let $h$ be an orientation preserving homeomorphism of the plane without fixed points. There exists a topological oriented foliation $\cl{F}$ of the plane such that each leaf of $\cl{F}$ is a Brouwer line for $h$.
\end{teorema}

Then, in \cite{lc2} it is proved the following improvement of Theorem \ref{pat1}.

\begin{teorema} \label{fol.br}
Let $M$ be a surface and $(H_t)_{t\in[0,1]}$ an isotopy in $M$ joining the identity to a homeomorphism $f$. For all $z\in M$ we define the arc $\gamma_z:t\mapsto H_t(z)$. We suppose that $f$ does not have any contractible fixed point $z$, that is, a fixed point $z$ such that $\gamma_z$ is a closed curve homotopic to a point. Then there exists an oriented topological foliation $\cl{F}$ in $M$ and for all $z\in M$ an arc positively transverse to $\cl{F}$ joining $z$ to $f(z)$ that is homotopic with fixed endpoints to the arc $\gamma_z$.
\end{teorema}

As an application of this theorem, we have the following proposition, which follows essentially from Theorem 9.1 in \cite{lc2}. We include a sketch of the proof.


\begin{proposicion}   \label{teo.pat2}
Let $\tl{f}:\T^2\ra\T^2$ be a homeomorphism isotopic to the identity without contractible fixed points.
Fix an isotopy $(\tl{H}_t)_{t\in[0,1]}$ in $\T^2$ between $\tl{f}$ and the identity, and let $\cl{F}$ be the foliation of $\T^2$ transverse to $(\tl{H}_t)_{t\in[0,1]}$ given by Theorem \ref{fol.br}. Let $(H_t)_{t\in[0,1]}$ be the isotopy in $\R^2$ which is the lift of $(\tl{H}_t)$ and satisfies $H_0= \txt{Id}$, and let $f:\R^2\ra\R^2$ be the lift of $\tl{f}$ given by $f=H_1$. Let $\hat{\cl{F}}$ be the lift of $\cl{F}$ to $\R\times S^1$. 

There exists $\epsilon>0$ such that, if $\hat{x},\hat{y}\in\R\times S^1$, are points with lifts $x,y\in\R^2$ and:
\begin{itemize} 
\item there is an $\epsilon$-chain for $f$ from $x$ to $x+(0,m)$ for some $m\in\N$, and 
\item there is an $\epsilon$-chain for $f$ from $y$ to $y+(0,-n)$ for some $n\in\N$,
\end{itemize} 
then there exists a compact leaf $l\in \hat{\cl{F}}$ which is an essential curve that separates $\hat{x}$ from $\hat{y}$ (that is, $\hat{x}$ and $\hat{y}$ belong to different connected components of $\R\times S^1\setminus l$). In particular $\hat{x}\neq \hat{y}$.
\end{proposicion}

\begin{prueba}[Sketch of the Proof.]
Let $F:\R\times S^1 \ra\R\times S^1$ be the lift of $\tl{f}$ such that $F\circ \pi= \pi\circ f$. Let $(\hat{H}_t)_{t\in[0,1]}$ be the isotopy in $\R\times S^1$ between $F$ and the identity which is the lift of the isotopy $(\tl{H}_t)_{t\in[0,1]}$. By Theorem \ref{fol.br}, for every $\hat{x}\in\R\times S^1$ there exists an arc which is positively transverse to $\hat{\cl{F}}$, joins $\hat{x}$ to $F(\hat{x})$ and is homotopic with fixed extremes to the arc $\gamma_{\hat{x}}:\mapsto \hat{H}_t(\hat{x})$. By this, one can easily see that, for any $\hat{x}\in\R\times S^1$ there exists $\epsilon>0$ such that any point $\hat{z}$ in $B_{\epsilon}(\hat{x})$ can be joined to any point $\hat{z}'$ in $B_{\epsilon}(F(\hat{x}))$ by an arc which is positively transverse to $\hat{\cl{F}}$ and homotopic to an arc of the form $\gamma_{\hat{z}\hat{x}}\gamma_{\hat{x}}\gamma_{F(\hat{x})\hat{z}'}$, where $\gamma_{\hat{z}\hat{x}}$ joins $\hat{z}$ to $\hat{x}$ in $B_{\epsilon}(\hat{x})$ and $\gamma_{F(\hat{x})\hat{z}'}$ joins $F(\hat{x})$ to $\hat{z}'$ in $B_{\epsilon}(F(\hat{x}))$, and where the product of two arcs stands for their concatenation. 

As $F$ is the lift of the homeomorphism $\tl{f}:\T^2\ra\T^2$, and as $\T^2$ is compact, there exists $\eta>0$ such that for any point $\hat{x}\in\R\times S^1$, any point in $B_{\eta}(\hat{x})$ can be joined to any point in $B_{\eta}(F(\hat{x}))$ by an arc positively transverse to $\hat{\cl{F}}$ as above. Also, by the continuity of $F$, there is $0<\epsilon<\eta$ such that for any $\hat{x}\in\R\times S^1$, if $\{\hat{x}_i\}_{i=0}^n$ is a periodic $\epsilon$-chain for $F$ with $\hat{x}_0=\hat{x}_n=\hat{x}$, then $\hat{x}_{n-1}\in B_{\eta}(F^{-1}(\hat{x}))$.

Suppose then that there are $\hat{x},\hat{y}\in\R\times S^1$ with lifts $x,y\in\R^2$ such that there is an $\epsilon$-chain $\{x_i\}_{i=0}^{n_1}$ for $f$ with $x_0=x$ and $x_{n_1}=x+(0,m)$ for some $m\in\N$, and an $\epsilon$-chain $\{y_i\}_{i=0}^{n_2}$ for $f$ with $y_0=y$ and $y_{n_2}=y+(0,-n)$ for some $n\in\N$. Then, we can construct a sequence of arcs $(\gamma_n)_{n=1}^{n_1}$ positively transverse to $\cl{F}$, and such that:
\begin{itemize}
\item $\gamma_1$ joins $x_0$ to $f(x_0)$,
\item $\gamma_i$ joins $f(x_{i-2})$ to $f(x_{i-1})$ for $2\leq i \leq n_1-2$, 
\item $\gamma_{n_1-1}$ joins $f(x_{n_1-3})$ to $f^{-1}(x+(0,m))$, and
\item $\gamma_{n_1}$ joins $f^{-1}(x+(0,m))$ to $x+(0,m)$. 
\end{itemize}

Then, letting $\gamma= \Pi_{i=1}^{n_1}\gamma_i$, we have that $\gamma$ is an arc positively transverse to $\cl{F}$ joining $x$ to $x+(0,m)$. 

Analogously, we construct an arc $\beta$ positively transverse to $\cl{F}$ and joining $y$ to $y+(0,-n)$. In \cite{lc2} it is proved that $\gamma$ and $\beta$ project to disjoint (not necessarily simple) loops $\tl{\gamma}$, and $\tl{\beta}$ in $\R\times S^1$, and there is a connected component $U$ of $\R\times S^1\setminus (\tl{\gamma}\cup\tl{\beta})$ which is a topological essential annulus. As $\tl{\gamma}$ and $\tl{\beta}$ are positively transverse to $\hat{\cl{F}}$, then $\hat{\cl{F}}$ is transverse to the border of $U$, either inwards or outwards. By the Poincar\'e Bendixon theorem, there exists a closed essential leaf $l$ contained in $U$. As the points $\hat{x}$ and $\hat{y}$ belong to the border of $U$, $l$ separates $x$ from $y$.   
\end{prueba}



\subsection{Atkinson's Lemma.}    \label{sec.atk}

Let $(X,\mu)$ be a probabilty space, $T:X\ra X$ be an ergodic transformation with respect to $\mu$, and $\phi:X\ra\R$ a measurable map. We say that the pair $(T,\phi)$ is \textit{recurrent} if for any measurable set $A\subset X$ of positive measure, and every $\epsilon>0$ there is $n>0$ such that 
$$ \mu \left( A\cap T^{-n}(A) \cap \left\{ x\, : \, \sum_{i=0}^{n-1} |\phi(T^i(x))| < \epsilon \right\} \right) >0. $$

In \cite{atk} it is proved the following theorem. 

\begin{teorema}    
Let $(X,\mu)$ be a non-atomic probability space, $T:X\ra X$ an ergodic automorphism, and $\phi:X\ra \R$ an integrable function. Then, the pair $(T,\phi)$ is recurrent if and only if $\int \phi d\mu =0$.
\end{teorema}  

From this theorem, it is not difficult to obtain the following corollary, usually known as `Atkinson's Lemma'.


\begin{corolario}      \label{atk}
Let $X$ be a separable metric space and $\mu$ a non-atomic Borel probability measure in $X$ which is ergodic with respect to a measurable transformation $T:X\ra X$. Let $\phi:X\ra\R$ be an integrable function, with $\int \phi d\mu=0$. Then, for $\mu$-almost every $x\in X$, there is an increasing sequence of integers $n_i$ with 
$$ T^{n_i}(x) \ra x      \    \txt{ and }     \    \left| \sum_{j=0}^{n_i-1} \phi(T^j(x)) \right| \ra 0 \ \ \ \ \ \txt{ as }  i\ra \infty.$$
\end{corolario}

\section{Proof of Theorem C from Theorem A and Addendum B.}             \label{sec.teofromteo1}

If $\rho(f)$ has irrational slope and contains a rational point, then by Theorem \ref{patirrac} this point is realized by a periodic orbit. Then we are left then with the case that $\rho(f)$ has rational slope and contains rational points. 

We will prove now that if there is a rational point $v\in\rho(f)$ that is not realized by a periodic orbit, then there is a power of $\tl{f}$ that is topologically conjugate to a homeomorphism $\tl{g}:\T^2\ra\T^2$ satisfying the hypotheses of Theorem A; that is, $\tl{g}$ has a lift $g:\R^2\ra\R^2$ such that $\rho(g)$ is a vertical interval containing the origin in its interior, and such that $(0,0)\in\rho(g)$ is not realized by a periodic orbit. 

By Remark \ref{r.mat}, there is $A\in \txt{GL}(2,\Z)$ such that $\rho(AfA^{-1})=A\rho(f)$ is a vertical segment, containing the rational point $v'=(p_1/k,p_2/k)$ given by $v'=Av$. By Remark \ref{r.pot}, if $g_0=(AfA^{-1})^{k}$, then $\rho(g_0)=k\rho(AfA^{-1})$, and then $\rho(g_0)$ is a vertical interval containing the point $w=kv'=(p_1,p_2)\in\Z^2$. We know that $\rho(T_1^{-p_1}T_2^{-p_2}g_0)=T_1^{-p_1}T_2^{-p_2}\rho(g_0)$, and therefore, if $g=T_1^{-p_1}T_2^{-p_2}g_0$, $\rho(g)$ is a vertical interval containing the point $T_1^{-p_1}T_2^{-p_2}(w)=(0,0)$. Let $\tl{A}$, $\tl{g}$ and $\tl{g}_0$ be the homeomorphisms of $\T^2$ lifted by $A$, $g$ and $g_0$, respectively. Then $\tl{g}=\tl{g}_0$, and as $g_0=(AfA^{-1})^{k}=Af^{k}A^{-1}$ we have that $\tl{g}$ and $\tl{f}^{k}$ are conjugate by $\tl{A}$. 

It remains to see that $(0,0)\in\rho(g)$ is not realized by a periodic orbit for $\tl{g}$. As $v\in\rho(f)$ is not realized by a periodic point of $\tl{f}$, then $v'=Av\in\rho(AfA^{-1})$ is not realized by a periodic point for $\tl{A}\tl{f}\tl{A}^{-1}$, and then $w=kv'=(p_1,p_2)\in\rho(g_0)$ is not realized by a periodic point of $\tl{g}_0=\tl{A}\tl{f}^k\tl{A}^{-1}$. Therefore $(0,0)=T_1^{-p_1}T_2^{-p_2}(w)\in\rho(g)$ is not realized by a periodic point of $\tl{g}$, as we wanted. \\
\\
\textbf{There is an annular, essential, set $A\subset\T^2$ which is periodic for $\tl{f}$}.\\
As $\tl{f}^{k}$ is conjugate to a homeomorphism satisfying the hypotheses of Theorem A, there is an annular, essential set $A$ which is invariant for $\tl{f}^{k}$, and therefore periodic for $\tl{f}$. \\
\\
\textbf{There is a simple, closed, essential curve in $\T^2$ that is free forever for $\tl{f}$.}\\
We saw above that $\tl{f}$ is topologically conjugate to a homeomorphism $\tl{f}_0$ with a lift of the form $f_0=A f A^{-1}$, with $A\in\txt{GL}(2,\Z)$, and such that $\rho(f_0)$ is a vertical interval. Also, we saw that there is $k\in\N$ such that $\tl{g}=\tl{f}_0^k$ satisfies the hypotheses of Theorem A. 

Then, there is a finite family $\{l_i\}_{i=0}^{r-1}$ of vertical closed curves that are free forever for $\tl{f}_0^k$, and such that the maximal invariant set $C$ of one of the annuli $[l_i,l_{i+1}]$ for $\tl{f}_0^k$ is annular, essential and vertical. Up to taking a subfamily, we may assume that the cardinality $r$ of $\{l_i\}_{i=0}^r$ is minimal such that the conclusion of that theorem holds. Observe that as $\tl{f}$ is isotopic to the identity, $\tl{f}_0^i(C)$ is vertical, for every $i$.

\begin{lema}     \label{lemaBA}
Let $m\in\N$ be the minimal period of $C$ for $\tl{f}_0$. Then $\tl{f}_0^i(C)\cap C=\empt$ for every $0 < i < m$. 
\end{lema}

By the fact that $C$ is annular and essential, this lemma gives us that the complement of the orbit of $C$ for $\tl{f}_0$ is a union of open annuli. That is, we have
$$\left( \bigcup_{i=0}^{m-1} \tl{f}_0^i(C)\right)^c=\bigcup_{i=0}^m U_i,$$ 
where the $U_i$ are open, essential, vertical and pairwise disjoint annuli. The minimal period of each $U_i$ is $m$, and $\cup l_i \subset \cup U_i$. Using this, we will prove the following.

\begin{lema}   \label{lemaBA2}
There is a simple, closed, essential curve $\gamma$ contained in $U_{i_0}$, such that $\tl{f}_0^m(\gamma)\cap \gamma=\empt$.
\end{lema}

Lemma \ref{lemaBA2} implies that there exists an essential curve that is free forever for $\tl{f}$. To see this, let $\gamma$ be as in Lemma \ref{lemaBA2}. By the fact that $\tl{f}_0^m(U_{i_0})=U_{i_0}$ we get that 
$$\tl{f}_0^{im}(\gamma)\cap \gamma=\empt \ \ \ \ \txt{for all $i\in\Z$,} $$
and as the minimal period of $U_{i_0}$ is $m$, we then get that $\tl{f}_0^{im+j}(\gamma)\cap \tl{f}_0^{im}(\gamma)=\empt$ for all $i\in\Z$ and $j\in\{0,\ldots,m-1\}$. That is, $\gamma$ is free forever for $\tl{f}_0$. As $\tl{f}_0$ is conjugate to $\tl{f}$, there is also an essential curve that is free forever for $\tl{f}$, as we wanted.

We are left with the proofs of lemmas \ref{lemaBA} and \ref{lemaBA2}.

\begin{prueba}[Proof of Lemma \ref{lemaBA}]
Recall that $k\in\N$ is such that $\tl{f}_0^k$ satisfies the hypotheses of Theorem A, the family $\{l_i\}_{i=0}^{r-1}$ is the family of curves in $\T^2$ given by that theorem applied to $\tl{f}_0^k$, and the annular essential set $C$ is the maximal invariant set of one of the annuli $[l_i,l_{i+1}]$ for $\tl{f}_0^k$.

It is easy to verify that the family $\{\tl{f}_0(l_i)\}_{i=0}^{r-1}$ also satisfies the conclusion of Theorem A applied to $\tl{f}_0^k$. Therefore, if we denote
$$\Theta_i= \bigcap_{n\in\Z} \tl{f}_0^{kn} \left( [l_i,l_{i+1}] \right)\ \ \ \txt{for } i\in\Z/r\Z, \txt{ and}$$
$$\wt{\Theta}_i= \bigcap_{n\in\Z} \tl{f}_0^{kn} \left( [\tl{f}_0(l_i),\tl{f}_0(l_{i+1})] \right)\ \ \ \txt{for } i\in\Z/r\Z,$$
by Addendum B we have that $\{\Theta_i\}_{i=0}^{r-1} = \{\wt{\Theta}_i\}_{i=0}^{r-1}$. As $\wt{\Theta}_i = \tl{f}_0(\Theta_i)$, this gives us that $\tl{f}_0$ permutes the family $\{\Theta_i\}_{i=0}^{r-1}$. The set $C$ is one of the sets $\Theta_i$, and then we conclude that $\tl{f}_0(C)\cap C=\empt$ for $0 < i < m$, as desired.
\end{prueba}

In \cite{kk} it is proved the following lemma (Lemma 3.2 in that article). A \textit{vertical} line in $\R^2$ is the lift of a vertical simple closed curve in $\T^2$.

\begin{lema}    \label{art54}
Let $h:\R^2\ra\R^2$ be a lift of a torus homeomorphism which is homotopic to the identity. Let $n\in\N$ and suppose that $\ell\subset\R^2$ is a vertical Brouwer line for $h^n$. Then, there is a vertical Brouwer line $\ell'$ for $h$. Also, if $S\subset\R^2$ is an open set containing the curves $\ell,h(\ell),\ldots,h^n(\ell)$, then the curve $\ell'$ can be chosen to be contained in $S$.
\end{lema}

We will use Lemma \ref{art54} in the proof of Lemma \ref{lemaBA2}.

\begin{prueba}[Proof of Lemma \ref{lemaBA2}]
We recall that $k\in\N$ is such that $\tl{f}_0^k$ satisfies the hypotheses of Theorem A, and $m$ is the minimal period for $\tl{f}_0$ of the annular set $C$.

Fix one of the curves $l_i$ and denote it by $L$. If $\tl{f}_0^m(L)\cap L=\empt$, then setting $\gamma:= L$ the lemma follows. Otherwise, suppose that $\tl{f}_0^m(L)\cap L\neq\empt$. Recall that $U_{i_0}$ is the connected component of the complement of the orbit of $C$ that contains $L$, and that $\tl{f}^m(U_{i_0})=U_{i_0}$. Let $V$ be a connected component of $\pi'^{-1}(U_{i_0})$, and let $F$ be a lift of $\tl{f}_0^m$ such that $F(V)=V$. Let $\ell\subset\R^2$ be the lift of $L$ contained in $V$. As $m$ divides $k$, $\tl{f}_0^k(U_{i_0})=U_{i_0}$, and the curve $L$ was chosen such that it is free forever for $\tl{f}_0^k$. In particular, $\tl{f}_0^k(L)\cap L=\empt$. That is, $F^{k/m}(\ell)\cap\ell=\empt$. Since $F$ fixes $V$, and as the curves $\ell, F(\ell),\ldots, F^{k/m}(\ell)$ are contained in $V$, by Lemma \ref{art54} there is a vertical Brouwer line $\Gamma$ for $F$ contained in $V$. Let $\gamma=\pi'(\Gamma)\subset U_{i_0}$. Then $\gamma$ is a vertical circle, and $\tl{f}_0^m(\gamma)\cap\gamma=\empt$, which concludes the proof of the Lemma 
\end{prueba}

To finish the proof of Theorem C, it remains the following. \\
\\
\textbf{If $q$ is the minimal period of $A$, the sets $A$, $\tl{f}(A)$, \ldots, $\tl{f}^{q-1}(A)$ are pairwise disjoint.}\\
We showed above that $\tl{f}$ is topologically conjugate to a homeomorphism $\tl{f}_0$ such that there is an essential annular set $C$ which is vertical and periodic for $\tl{f}_0$. In Lemma \ref{lemaBA} it was proved that, if $m$ is the minimal period of $C$ for $\tl{f}_0$, then $\tl{f}_0^i(C)\cap C=\empt$ for $0 < i < m$. 

This implies that, for the homeomorphism $\tl{f}$, there is an essential, annular set $A$ which is periodic for $\tl{f}$, with minimal period $q=m$, and such that $A,\tl{f}(A),\ldots,\tl{f}^{q-1}(A)$ are pairwise disjoint, as desired.

\section{Proof of Theorem A, part I: construction of the curves $\tl{l}_i$, and items (2), (3) and (4).}            \label{sec.lema.prel}

We begin by mentioning a related result for homeomorphisms of the compact annulus by Le Calvez. For a homeomorphism $F:S^1\times [0,1]\ra S^1\times [0,1]$ isotopic to the identity, the rotation set of some lift $f:\R\times[0,1]\ra\R\times[0,1]$ is defined as the set of all accumulation points of sequences of the form 
\begin{equation}     \label{eqpat1}
\left\{ \frac{f^{m_i}(x_i)_1-(x_i)_1}{m_i} \right\}_{i\in\N}
\end{equation}
where $m_i\ra\infty$ and $x_i\in\R\times [0,1]$. In this case the rotation set $\rho(f)$ is a compact interval $I\subset\R$ (possibly degenerate). 

Also, if $\Lambda\subset S^1\times[0,1]$ is a compact invariant set, we can define the rotation set of $\Lambda$, denoted $\rho(\Lambda,f)$, as the set of all accumulation points of sequences of the form (\ref{eqpat1}), with $m_i\ra\infty$ and $x_i\in\Pi^{-1}(\Lambda)$, where $\Pi:\R\times [0,1]\ra S^1\times [0,1]$ is the canonical projection. 

The following theorem was proven for $C^1$ diffeomorphisms in \cite{lc3}, but by the results of \cite{lc2} it is also valid for homeomorphisms (see Theorem 9.1 in that article). 

\begin{figure}[h]        
\begin{center} 
\includegraphics{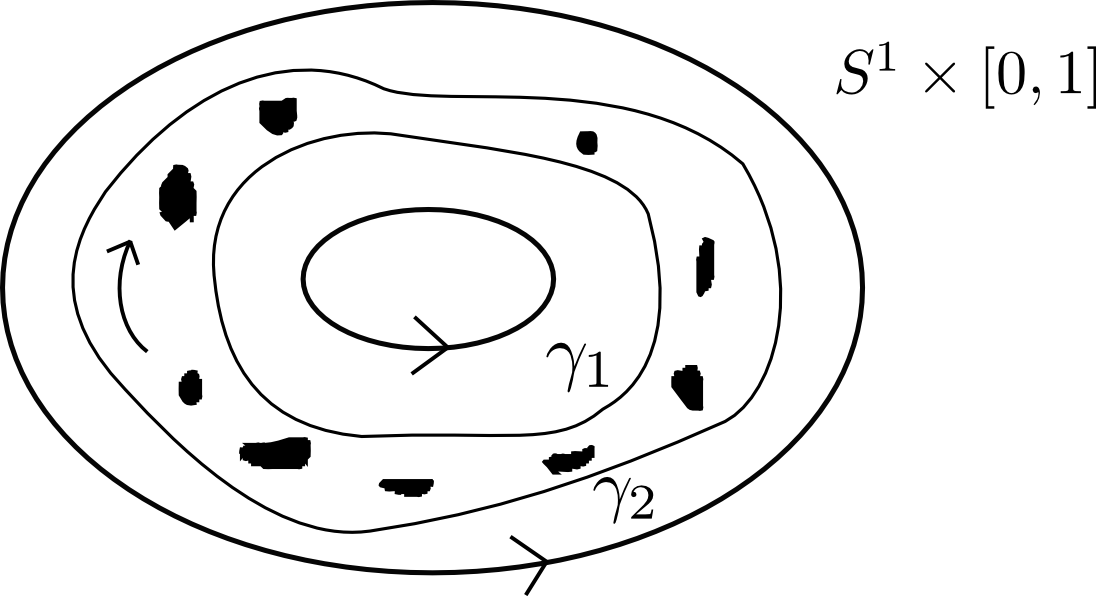}
\caption{Illustration for Theorem \ref{teo.intro.pat}. The free curves $\gamma_1$ and $\gamma_2$ are free for $F$.}
\label{fig.teopatintro}
\end{center}  
\end{figure}

\begin{teorema}      \label{teo.intro.pat}
Let $F:S^1\times[0,1]\ra S^1\times[0,1]$ be a homeomorphism isotopic to the identity with a lift $f:\R\times[0,1]\ra\R\times[0,1]$ that has no fixed points and whose rotation set is an interval containing $0$ in its interior.

Then, there exists a finite non-empty family $\{\gamma_i\}$ of essential, pairwise disjoint, free curves for $F$ such that, the maximal invariant set contained between two consecutive curves has rotation set contained either strictly to the right or strictly to the left of $0$ (see Fig. \ref{fig.teopatintro}).
\end{teorema}

In Section \ref{sec.curvasli} we will adapt this result for the torus case, and in this way we will construct a family of curves satisfying item (3) from Theorem A.

\subsection{Construction of the curves $\tl{l}_i$ satisfying item (3).}     \label{sec.curvasli}

We will prove the following.

\begin{proposicion}       \label{prop.prel}
Let $\tl{f}$ and $f$ be as in Theorem A. There exists a finite family $\{\tl{l}_i\}_{i=0}^{r-1}$ of pairwise disjoint, simple, closed, essential and vertical curves $\tl{l}_i\subset\T^2$ such that if $\Theta_i$ is the maximal invariant set of $[\tl{l}_i,\tl{l}_{i+1}]$, then $\Theta_i$ is non-empty, and $\rho(\Theta_i,f)$ is contained either in $\{0\}\times \R^+$ or in $\{0\}\times \R^-$.

The family $\{\tl{l}_i\}$ is minimal in the sense that, for any $i$, if $\rho(\Theta_i,f)\subset \{0\}\times \R^+$ then $\rho(\Theta_{i+1\mod r},f)\subset \{0\}\times \R^-$, and if $\rho(\Theta_i,f)\subset \{0\}\times \R^-$ then $\rho(\Theta_{i+1\mod r},f)\subset \{0\}\times \R^+$. 
\end{proposicion}

\begin{remark}     \label{remark.prop.prel}   
As the sets $\Theta_i$ are compact, the sets $\rho(\Theta_i,f)$ are also compact (see Remark \ref{remark.compact}). Therefore, the fact that $\rho(\Theta_i,f)$ is contained either in $\{0\}\times(0,\infty)$ or in $\{0\}\times (-\infty,0)$ means actually that $\rho(\Theta_i,f)$ is contained either in $\{0\}\times(\e,\infty)$ or in $\{0\}\times (-\infty,-\e)$, for some $\e>0$ and for any $i$. Therefore the family $\{\tl{l}_i\}$ given by Proposition \ref{prop.prel} satisfies item $(3)$ from Theorem A.
\end{remark}

\begin{remark}     \label{remark.nofix}
By the Brouwer Translation Lemma (see Section \ref{sec.brth}), the hypothesis in Theorem A that $(0,0)\in\rho(f)$ is not realized by a periodic orbit is equivalent to the fact that $f$ has no fixed points. Therefore all the results from Section \ref{sec.brth} apply to $f$.
\end{remark}

To prove Proposition \ref{prop.prel} it will be convenient to work on the lift $\R\times S^1$ of $\T^2$. Recall our notation for the canonical projections:
$$\R^2 \stackrel{\pi}{\ra} \R\times S^1 \stackrel{\pi''}{\ra} \T^2, \ \ \ \txt{and } \ \  \pi'= \pi''\circ \pi.$$
We will first prove the following.

\begin{figure}[h] 
\begin{center} 
\includegraphics{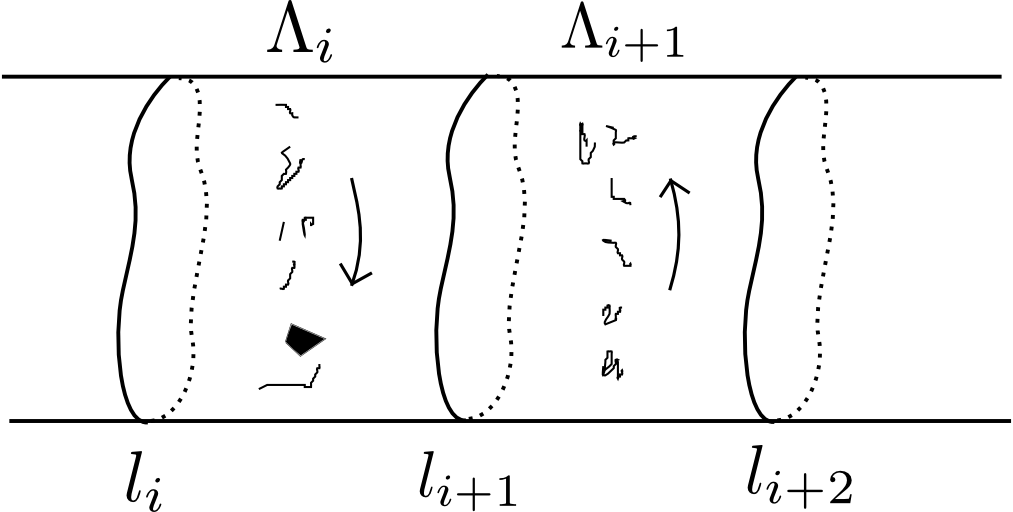}
\caption{The sets $\Lambda_i$ and the curves $l_i$.}
\label{fig.lambdas}  
\end{center}  
\end{figure}

\begin{lema} \label{lema.pat}
For $\tl{f}$ and $f$ as in Theorem A, let $F:\R\times S^1\ra\R\times S^1$ be the lift of $\tl{f}$ such that $F\circ \pi=\pi\circ f$. Then:
\begin{enumerate}
\item \label{lema.pat1}
The chain recurrent set $CR(F)$ is not empty, and $CR(F)=\Lambda^+\cup \Lambda^-$, where $\Lambda^+$ and $\Lambda^-$ are closed disjoint $F$-invariant sets such that, denoting $\wt{\Lambda}^{\pm} =\pi''(\Lambda^{\pm})\subset \T^2$, we have $\rho(\wt{\Lambda}^+,f)\subset \{0\}\times(\epsilon,\infty)$ and $\rho(\wt{\Lambda}^-,f)\subset \{0\}\times(-\infty,-\epsilon)$, for some $\epsilon>0$.
\item There exist simple, closed, essential curves $l_0\prec l_1\prec \cdots l_{r}=T_1(l_0)$ on $\R\times S^1$ which are free for $F$, and such that they `separate' $\Lambda^+$ from $\Lambda^-$, that is:
\begin{enumerate}
\item $CR(F)\bigcap \cup_{i=0}^r l_i=\emptyset$,
\item for $0\leq i < r$, the set $\Lambda_i:=CR(F)\cap(l_i,l_{i+1})$ is compact, non-empty and $F$-invariant,
\item for $0\leq i < r$, either $\Lambda_i\subset\Lambda^+$ or $\Lambda_i\subset\Lambda^-$, and
\item if $\Lambda_i\subset \Lambda^+$, then $\Lambda_{i+1}\subset \Lambda^-$, and if $\Lambda_i\subset \Lambda^-$ then $\Lambda_{i+1}\subset \Lambda^+$, for any $0\leq i< r-1$ (see Fig. \ref{fig.lambdas}).

\end{enumerate}   
\end{enumerate}
\end{lema}

\begin{prueba} 
First we observe the following elementary fact. There exists an isotopy $(\tl{H}_t)_{t\in[0,1]}$ between the identity and $\tl{f}$ with the property that if $(H_t)$ is the lift of $(\tl{H}_t)$ with $H_0=\txt{Id}$, then $H_1 = f$. To see this just observe that if $(\tl{H}_t')_{t\in[0,1]}$ is any isotopy between the identity and $\tl{f}$, and if $(H'_t)$ is the lift of $(\tl{H}_t')$ with $H_0'=\txt{Id}$, then $H_1'=f+(a,b)$, for some $a,b\in\Z$. Defining $H_t=H_t'+t(-a,-b)$, for $t\in [0,1]$, we have that $(H_t)$ is an isotopy between the identity and $f$ wich projects to an isotopy $(\tl{H}_t)$ on $\T^2$ between the identity and $\tl{f}$ with the desired properties. 

Now, let $\cl{F}$ be the Brouwer foliation of $\T^2$ transversal to $\tl{H}$ given by Theorem \ref{fol.br}. Let $\hat{\cl{F}}$ be the lift of $\cl{F}$ to $\R\times S^1$. \\
\\
\textit{Item 1.} \\
\textbf{$CR(F)$ is non-empty}. Let $\phi:\T^2\ra \R^2$ be given by 
$$\phi(\tl{x})=f(x)-x,$$
where $x\in\R^2$ is any lift of $\tl{x}$, and let $\phi_1=\txt{pr}_1\circ\phi$. Then 
$$\sum_{i=0}^{n-1}\phi_1(\tl{f}^i(\tl{x}))= f^{n}(x)_1-x_1.$$ 
Let $\mu$ be any ergodic measure for $f$. By hypothesis, $\int \phi_1 d\mu= \rho(\mu,f)_1=0$, and by Atkinson's Lemma \ref{atk} there exists a full $\mu$-measure set $X\subset \T^2$ such that for any $\tl{x}\in X$ and $x\in\pi'^{-1}(\tl{x})$ we have that there is a sequence of positive integers $n_i$ such that
$$\tl{f}^{n_i}(\tl{x})\ra \tl{x}  \ \ \txt{  and  } \ \ \left|\sum_{j=0}^{n_i-1}\phi_1(\tl{f}^{j}(\tl{x}))\right|=|f^{n_i}(x)_1-x_1| \ra 0, \ \ \ \txt{ as } i\ra\infty.$$
This means that $\pi(x)$ is recurrent for $F$, and in particular $CR(F)\neq\emptyset$. 

\textbf{Definition of $\Lambda^+$ and $\Lambda^-$}. As $f$ has no fixed points by hypothesis (see Remark \ref{remark.nofix}), by Theorem \ref{franks} there is $\epsilon_1>0$ such that $f$ has no periodic $\epsilon_1$-chains, and by Proposition \ref{teo.pat2} there is $\epsilon_2>0$ such that, if there are $\hat{x},\hat{y}\in \R\times S^1$ with lifts $x,y\in\R^2$ such that:
\begin{itemize}
\item there is an $\epsilon_2$-chain for $f$ from $x$ to $x+(0,m)$ for some $m\in\N$, and 
\item there is an $\epsilon_2$-chain for $f$ from $y$ to $y+(0,-n)$ for some $n\in\N$,
\end{itemize}
then there exists a compact leaf of $\hat{\cl{F}}$ that separates $\hat{x}$ from $\hat{y}$. Let $\epsilon_0=\min\{\epsilon_1,\epsilon_2\}$. 

We define $\Lambda^+\subset\R\times S^1$ as the set of points $\hat{x}\in CR(F)$ such that, if $x\in\pi^{-1}(\hat{x})$, there exists an $\epsilon_0$-chain $\{x_i\}_{i=0}^n$ for $f$, with $x_0=x$ and $x_n=x+(0,m)$ for some $m\in\N$. Analogously, we define $\Lambda^-\subset\R\times S^1$ as the set of points $\hat{x}\in CR(F)$ such that, if $x\in\pi^{-1}(\hat{x})$, there exists an $\epsilon_0$-chain $\{x_i\}_{i=0}^n$ for $f$, with $x_0=x$ and $x_n=x+(0,-m)$, for some $m\in\N$.

\textbf{$\Lambda^+$ and $\Lambda^-$ are non-empty}. We prove that $\Lambda^+$ is non-empty; the case of $\Lambda^-$ is similar. As $\rho(f)$ is a vertical interval containing the origin in its interior, by Proposition \ref{erg}, there exists an ergodic measure $\mu$ with respect to $\tl{f}$ with $\rho(\mu)_2>0$. By Birkhoff's Theorem, there exists a set $X\subset\T^2$ of full $\mu$-measure such that for $\tl{x}\in X$ and $x\in\pi'^{-1}(\tl{x})$, we have 
$$\rho(\tl{x},f)=\lim_n \frac{1}{n}\sum_{i=0}^n \phi(\tl{f}^i(\tl{x}))=\int \phi d\mu=\rho(\mu,f).$$ 
By Atkinson's Lemma \ref{atk}, there exists a full $\mu$-measure set $X'\subset\T^2$ such that if $\tl{x}\in X'$ and $x\in\pi'^{-1}(\tl{x})$, there is a sequence of positive integers $n_i$ such that
$$ \tl{f}^{n_i}(\tl{x}) \ra \tl{x} \ \ \txt{ and } \ \ \left| \sum_{j=0}^{n_i-1}\phi_1(\tl{f}^j(\tl{x})) \right| =|f^{n_i}(x)_1-x_1| \ra 0, \ \ \ \txt{ as } i\ra\infty.$$

Let $\tl{y}\in X\cap X'$ and $y\in\pi'^{-1}(\tl{y})$. Then there is an increasing sequence of naturals $\{r_n\}_n$ and a sequence of integers $\{s_n\}_n$ such that
$$|f^{r_n}(y)-y-(0,s_n)| \ra 0 \ \ \ \txt{ as } n\ra\infty,$$
with $s_n/r_n\ra \rho(\tl{y},f)_2=\rho(\mu,f)_2>0$. Therefore $\lim_n s_n= \infty$, and in particular $s_n>0$ for $n$ sufficiently large. If $\hat{y}\in \pi''^{-1}(\tl{y})$, then $F^{r_n}(\hat{y})) \ra \hat{y}$ as $n\ra\infty$, and we have that $\hat{y}$ is recurrent for $F$, and in particular $\hat{y}\in CR(F)$. Therefore $\hat{y}\in\Lambda^+$, and $\Lambda^+$ is non-empty. 

\textbf{It holds $CR(F)=\Lambda^+\cup\Lambda^-$}. Observe that by the definition of the sets $\Lambda^+$ and $\Lambda^-$, we have that $\Lambda^+\cup\Lambda^-\subset CR(F)$, and then we only need to prove that $CR(F)\subset \Lambda^+\cup\Lambda^-$. Suppose by contradiction that there is $\hat{x}\in CR(F)\setminus(\Lambda^+\cup\Lambda^-)$. Then, by definition of $\Lambda^+$ and $\Lambda^-$, there exists an $\epsilon_0$-chain for $f$ starting and ending in $x$, that is, a periodic $\epsilon_0$-chain for $f$. By definition of $\e_0$, we have $\e_0\leq \e_1$, where $\e_1$ is the constant given by Theorem \ref{franks}, and then that theorem implies there is a fixed point for $f$, a contradiction. Therefore we must have $CR(F)\subset\Lambda^+\cup\Lambda^-$ as we wanted.

\textbf{The sets $\Lambda^+$ and $\Lambda^-$ are disjoint and closed}. We will prove that $\ol{\Lambda^+}\cap \ol{\Lambda^-}=\emptyset$. As $CR(F)=\Lambda^+\cup\Lambda^-$ and $CR(F)$ is closed, this will imply that $\Lambda^+$ and $\Lambda^-$ are closed and disjoint. Suppose by contradiction that there is $\hat{x}\in\ol{\Lambda^+}\cap\ol{\Lambda^-}$. Let $\hat{y}\in\Lambda^+$ and $\hat{z}\in\Lambda^-$, be such that $d(\hat{y},\hat{x})<\epsilon_0/3$ and $d(\hat{z},\hat{x})<\epsilon_0/3$, and let $\{y_i\}_{i=0}^{n_1}$, and $\{z_i\}_{i=0}^{n_2}$ be $\epsilon_0/3$-chains for $f$ such that $y_0\in\pi^{-1}(\hat{y})$, $y_{n_1}=y_0+(0,m_1),$ $z_0\in\pi^{-1}(\hat{z})$, $|z_0-y_0|<2\epsilon_0/3$ and
$z_{n_2}=y_0+(0,-m_2)$, for some $m_1,m_2\in\N$. We now show that we can concatenate integer translates of these chains $\{y_i\}$ and $\{z_i\}$ to get a periodic chain for $f$. For each $0\leq i < m_2$ define the $\epsilon_0/3$-chain $\{y_l^i\}_{l=0}^{n_1}$ for $f$ as the translate of $\{y_l\}_{l=0}^{n_1}$ by $T_2^{im_1}$, that is,
$$y_l^i = T_2^{im_1} y_l, \ \ \txt{ for } 0\leq l < n_1,$$
and for each $0\leq j < m_1$, define the $\epsilon_0/3$-chain $\{z_k^j\}_{k=0}^{n_2}$ for $f$ as the translate of $\{z_k\}_{k=0}^{n_2}$ by $T_2^{m_1m_2-jm_2}$, that is,
$$z_k^j = T_2^{m_1m_2-jm_2} z_k, \ \ \txt{ for }0\leq k < n_2.$$
Define then the $\epsilon_0$-chain $\{w_i\}_{i=0}^{n_1m_2+n_2m_1}$ for $f$ as the concatenation of the chains defined above, given by
$$w_{in_1+l}=y_l^i, \ \ \txt{for $0\leq i < m_2$ and $0\leq l < n_1$,}$$
$$w_{m_2n_1+jn_2+k}= z_k^j \ \ \txt{for $0\leq j < m_1$ and $0\leq k < n_2$, and}$$
$$w_{n_1m_2+n_2m_1}=w_0.$$ 
Then, $\{w_i\}_{i=0}^{n_1m_2+n_2m_1}$ is a periodic $\epsilon_0$-chain for $f$. By Theorem \ref{franks}, $f$ has a fixed point, which is a contradiction. Therefore there cannot be $\hat{x}\in\ol{\Lambda^+}\cap\ol{\Lambda^-}$. As we mentioned, this implies that $\Lambda^+$ and $\Lambda^-$ are closed and disjoint.   

Before proving that the sets $\Lambda^+$ and $\Lambda^-$ are $F$-invariant and the last claim of Item $1$, we will prove Item $2$ of the lemma. \\
\\
\textit{Item 2}. \\
\textbf{Construction of the family $\{l_i\}_{i=0}^r$}. By Proposition \ref{teo.pat2} and by the definition of the sets $\Lambda^+$ and $\Lambda^-$, for each $x\in\Lambda^+$, $y\in\Lambda^-$, there exists a compact leaf $l\in\hat{\cl{F}}$ that separates $x$ from $y$. So, the set $\cl{F}_c$ of compact leaves of $\hat{\cl{F}}$ is not empty. The union of the compact leaves of a foliation of $\T^2$ is compact (see for ex. \cite{hf}), and as $\hat{\cl{F}}$ is a lift of a foliation of $\T^2$, the set $\cup \cl{F}_c$ is closed as a subset of $\R\times S^1$ ($\cup\cl{F}_c$ denotes the union of the elements of $\cl{F}_c$). Observe that, as the leaves of $\cl{F}$ are Brouwer lines for $f$, the elements of $\cl{F}_c$ are free curves for $F$. \\
\\
\textit{Claim: }$CR(F)\ \cap \ \cup\cl{F}_c=\emptyset$. \\
To see this, note that as the elements of $\cl{F}_c$ are free curves for $F$, the points in $\cup\cl{F}_c$ are wandering for $F$. As the set of wandering points does not intersect the set of chain recurrent points, the claim follows. 

This claim gives us that $CR(F)$ has an open cover $\cl{U}'$ whose elements are the connected components of $\R\times S^1\setminus \cup \cl{F}_c$, which are sets of the form $(l,l')$, with $l,l'\in\cl{F}_c$. By definition of the sets $\Lambda^+$ and $\Lambda^-$, and by Proposition \ref{teo.pat2} we have that for any element $(l,l')$ of $\cl{U}'$,
$$\txt{either } CR(F)\cap (l,l')\subset\Lambda^+, \txt{ or } CR(F)\cap(l,l')\subset\Lambda^-.$$
Now, fix $l_*\in \cl{F}_c$. The compact set $CR(F)\cap [l_*,T_1(l_*)]$ has a finite subcover $\cl{U}''\subset\cl{U}'$, of the form $\cl{U}''=\{(l_{2i},l_{2i+1})\}_{i=0}^{s-1}$. We reindex the curves $l_i$ in a way that $l_i\prec l_{i+1}$ for $0\leq i<2s-1$, and we extract from the family of compact leaves $\{l_i\}_{i=0}^{2s-1}$ a subfamily, which we denote simply by $\{l_i\}_{i=0}^{r-1}$, indexed in a similar way, and which is minimal with respect to the following property: if $l_r=T_1(l_0)$, then for each $0\leq i < r$
$$\txt{either } \emptyset\neq CR(F)\cap (l_i,l_{i+1})\subset\Lambda^+, \txt{ or } \emptyset\neq CR(F)\cap(l_i,l_{i+1})\subset\Lambda^-.$$
As a consequence we have that, if for $0\leq i < r$ we define
$$\Lambda_i= CR(F)\cap (l_i,l_{i+1}),$$
then, for $0\leq i < r$
\begin{itemize}
\item $\Lambda_i\neq\emptyset$, and 
\item if $\Lambda_i\subset\Lambda^+$ then $\Lambda_{i+1}\subset\Lambda^-$, and if $\Lambda_i\subset\Lambda^-$ then $\Lambda_{i+1}\subset\Lambda^+$.
\end{itemize} 
This concludes the construction of the family $\{l_i\}_{i=0}^{r-1}$ satisfying items $(a)$, $(c)$ and $(d)$ from Item 2 of the lemma. As the sets $\Lambda^+$ and $\Lambda^-$ are compact and disjoint from the curves $l_i$, the sets $\Lambda_i$ are also compact. Since $\Lambda_i\neq\emptyset$ for all $i$, to prove that $\{l_i\}_{i=0}^{r-1}$ satisfies item $(b)$ it remains to prove that $\Lambda_i$ is $F$-invariant, for each $0\leq i < r$.

\textbf{For any $0\leq i < r$, $\Lambda_i$ is $F$-invariant}. Fix $i\in\{0,\ldots,r-1\}$. First we prove that $F(\Lambda_i)\subset\Lambda_i$.  As $CR(F)$ is $F$-invariant and $\Lambda_i=CR(F)\cap(l_i,l_{i+1})$, to show that $F(\Lambda_i)\subset\Lambda_i$ it suffices to show that if $x\in\Lambda_i$ then $F(x)\in (l_{i},l_{i+1})$.

Suppose this is not true. Then there exists $x_0\in \Lambda_{i}$ such that $F(x_0)\notin (l_{i},l_{i+1})$. Without loss of generality suppose that $F(x_0)\in \ol{R}(l_{i+1})$. Then, as $l_{i+1}$ is free for $F$ we must have that $l_{i+1}\prec F(l_{i+1})$. Let $\delta_1:=d(l_{i+1},F(l_{i+1}))>0$. By the continuity of $F$ there is $\delta_2>0$ such that if $d(x,l_{i+1})<\delta_2$, then $F(x)\in R(l_{i+1})$ and $d(F(x),l_{i+1})>\delta_1/2$. Let $\delta=\min\{\delta_2,\delta_1/2\}$, and let $\{y_i\}_{i=0}^s$ be any $\delta$-chain for $F$ with $y_0=x_0$. Then, either $y_1\in R(l_{i+1})$ or $d(y_1,l_{i+1})<\delta_2$, and in both cases we have $F(y_1)\in R(l_{i+1})$ and $d(F(y_1),l_{i+1})>\delta_1/2$. Therefore, $y_2\in R(l_{i+1})$, and $F(y_2)\in R(F(l_{i+1}))$. Then $y_3\in R(l_{i+1})$. By induction, we get that $y_n\in R(l_{i+1})$ for all $n\geq 2$. As $\{y_i\}_{i=0}^s$ was an arbitrary $\delta$-chain with $y_0=x_0$, we then have that $x_0$ is not $\delta$-chain recurrent, which contradicts that $x_0\in \Lambda_{i}\subset CR(F)$. This contradiction gives us that $F(x_0)$ must be contained in $(l_i, l_{i+1})$, and therefore $F(\Lambda_i)\subset\Lambda_i$.

Now we prove that $F^{-1}(\Lambda_i)\subset\Lambda_i$. Applying the arguments in last paragraph to $F^{-1}$ we get that $F^{-1}(CR(F^{-1})\cap(l_i,l_{i+1})) \subset CR(F^{-1}) \cap(l_i,l_{i+1})$, and as $CR(F)=CR(F^{-1})$ we get that $F^{-1}(\Lambda_i)\subset\Lambda_i$.

As the choice of $i$ was arbitrary, we conclude that for any $i$, $\Lambda_i$ is $F$-invariant, as we wanted. This finishes the proof of Item $2$ of the lemma.\\

Now we proceed to the remaining part of the proof of Item $1$. 

\textbf{The sets $\Lambda^+$ and $\Lambda^-$ are $F$-invariant.} We proved that for each $i$ the set $\Lambda_i$ is $F$-invariant, and contained either in $\Lambda^+$ or in $\Lambda^-$. As both $\Lambda^+$ and $\Lambda^-$ are contained in $\cup_{n,i} T_1^n(\Lambda_i)$, we conclude that $\Lambda^+$ and $\Lambda^-$ are $F$-invariant. 

\textbf{There is $\epsilon>0$ such that $\rho(\wt{\Lambda}^+,f)\subset \{0\}\times (\epsilon,\infty)$, and $\rho(\wt{\Lambda}^-,f)\subset \{0\}\times (-\infty,-\epsilon)$}. We will deal only with the case of $\rho(\wt{\Lambda}^+,f)$; the case of $\rho(\wt{\Lambda}^-,f)$ is analogous. As $\Lambda^+$ is closed and $F$-invariant, $\wt{\Lambda}^{\pm}$ is a compact $\tl{f}$-invariant set. 
Let $v^-$ the lower endpoint of $\txt{conv}\,\rho(\wt{\Lambda}^+,f)$. It suffices to prove that $(v^-)_2 >0$. 

Suppose by contradiction that $(v^-)_2 \leq 0$. By Proposition \ref{erg}, as $v^-$ is an extremal point of $\txt{conv}\,\rho(\wt{\Lambda}^+,f)$, there exists an ergodic measure $\mu$ for $\tl{f}$ with $\rho(\mu,f)=v^-$ and $\txt{supp}(\mu)\subset\wt{\Lambda}^+$.
Let $\e_0$ be as in the definition of $\Lambda^+$ and $\Lambda^-$. As $(v^-)_2\leq 0$, proceeding as in the proof that $\Lambda^+$ and $\Lambda^-$ are non-empty, with the aid of Birkhoff's Theorem and Atkinson's Lemma \ref{atk} we find a point $\tl{x}\in\txt{supp}(\mu)$ such that
$$| f^n(x)-x-(0,-m) | <\epsilon_0$$
for $x\in\pi'^{-1}(\tl{x})$ and for some $n\in\N$ and $m \in \N_0$. If $m>0$ this means that $\pi(x)\in\Lambda^-$,
and therefore $\Lambda^+\cap\Lambda^-\neq\emptyset$. This is a contradiction, and then we cannot have that $m>0$. If $m=0$ we have that $x$ is $\epsilon_0$-chain recurrent for $f$, which by the definition of $\epsilon_0$ and by Theorem \ref{franks} is also a contradiction. Therefore we cannot have that $(v^-)_2 \leq 0$, as we wanted. This finishes the proof of Item $1$, and of Lemma \ref{lema.pat}.
\end{prueba}

\begin{remark}    \label{remark.maxinv}
The fact that the sets $\Lambda_i\subset (l_i,l_{i+1})$ are non-empty and $F$-invariant implies the following. If $\ell\subset\R^2$ is a lift of $l_i$ and if $\ell'\subset\R^2$ is the lift of $l_{i+1}$ such that $\ell \prec \ell' \prec T_1(\ell)$, then
$$\bigcap_{n\in\Z}f^n\left((\ell,\ell')\right)=\pi^{-1}\left(\bigcap_{n\in\Z}F^n\left((l_i,l_{i+1})\right)\right)\supset \pi^{-1}(\Lambda_i)\neq\empt.$$ 
\end{remark}

We now are ready to prove Proposition \ref{prop.prel}.

\begin{prueba}[Proof of Proposition \ref{prop.prel}.]            
By construction,
the curves $\tl{l}_i=\pi''(l_i)\subset\T^2$ are compact leaves from a foliation of $\T^2$, and therefore pairwise disjoint. As the curves $l_i\subset\R\times S^1$ are essential, and as $\pi'':\R\times S^1\ra \T^2$ is a covering map, the curves $\tl{l}_i\subset\T^2$ are also essential, and by the definition of $\pi''$ it is easy to see that they are vertical. For any $0\leq i < r$, let $\Theta_i\subset\T^2$ be as in Theorem A. Then, we observe that $\Theta_i$ is non-empty: if $F$ and $\Lambda_i\subset\R\times S^1$ are as in Lemma \ref{lema.pat}, then by that lemma $\Lambda_i$ is non-empty and $F$-invariant, and as $F$ lifts $\tl{f}$, the set $\pi''(\Lambda_i)\subset (\tl{l}_i,\tl{l}_{i+1})$ is non-empty and $\tl{f}$-invariant. As $\Theta_i$ is the maximal invariant set in $(\tl{l}_i,\tl{l}_{i+1})$ for $\tl{f}$, $\empt\neq \pi''(\Lambda_i)\subset\Theta_i$. 

We now prove that $\rho(\Theta_i,f)$ is contained either in $\{0\}\times(0,\infty)$ or in $\{0\}\times (-\infty,0)$ for each $i$. 

Let $i\in\{0,1,\ldots,r-1\}$, and suppose first that $\Lambda_i\subset\Lambda^+$, where $\Lambda^+$ is as in Lemma \ref{lema.pat}. Then we will prove that $\rho(\Theta_i,f) \subset \{0\}\times(0,\infty)$. Let $v^-$ be the lower endpoint of the interval $\txt{conv}(\rho(\Theta_i,f))$. To prove that $\rho(\Theta_i,f) \subset \{0\}\times(0,\infty)$, it suffices then to prove that $(v^-)_2>0$. By contradiction, suppose that $(v^-)_2\leq 0$.

By Proposition \ref{erg}, we can find an ergodic measure $\mu$ with support contained in $\Theta_i$ and with $\rho(\mu,f) = v^-$. As in the proof in Lemma \ref{lema.pat} that the sets $\Lambda^+$ and $\Lambda^-$ are non-empty, with the use of Atkinson's Lemma we can find a point $\tl{x}\in\txt{supp}(\mu)$ such that for any $\e>0$ and $x\in\pi'^{-1}(\tl{x})$, there is $n>0$ and $m\leq 0$ such that 
\begin{equation}      \label{z200}
|f^n(x)-x-(0,m)|<\e.
\end{equation}
Therefore if $\hat{x}\in\pi''^{-1}(\tl{x})\subset\R\times S^1$, $\hat{x}$ is recurrent for $F$, and in particular $\hat{x}\in CR(F)$. As $\tl{x}\in\txt{supp}(\mu)\subset [\tl{l}_i,\tl{l}_{i+1}]$, $\hat{x}$ belongs to some integer translate of $[l_i,l_{i+1}]$. Therefore, there is an integer translate $\hat{x}'$ of $\hat{x}$ such that $\hat{x}'\in CR(F)\cap [l_i,l_{i+1}]=\Lambda_i$. As in (\ref{z200}) $\e>0$ can be taken arbitrarily and $m\leq 0$, we have that $\hat{x}'\in\Lambda^-$, and then $\Lambda_i\cap\Lambda^-\neq\empt$, which contradicts our assumption that $\Lambda_i\subset\Lambda^+$. Therefore we must have $(v^-)_2>0$, and this proves that $\rho(\Theta_i,f) \subset \{0\}\times(0,\infty)$, as we wanted.

Similarily, if $\Lambda_i\subset\Lambda^-$ we prove that $\rho(\Theta_i,f)$ is contained in $\{0\}\times (-\infty,0)$. The choice of $i\in\{0,\dots,r-1\}$ was arbitrary, and then by Remark \ref{remark.prop.prel} we have that, for each $i$, either $\rho(\Theta_i,f)\subset\{0\}\times (-\infty,-\e)$ or $\rho(\Theta_i,f)\subset\{0\}\times (\e,\infty)$ for some $\e>0$, and we conclude that	for the family $\{\tl{l}_i\}_{i=0}^{r-1}$ it holds item $(3)$ from Theorem A 
\end{prueba}

\subsection{Item $(1)$ implies items $(2)$ and $(4)$.}      \label{sec.24from1}

We recall that, in Theorem A, the sets $\Theta_i$ are the maximal invariant sets in $[\tl{l}_i,\tl{l}_{i+1}]$ for $\tl{f}$. Also, we recall items $(1)$, $(2)$ and $(4)$ from that theorem.

\begin{enumerate}

\item[\textit{(1)}] One of the sets $\Theta_i$ is annular, essential, vertical and a semi-attractor for $\tl{f}$.
\item[\textit{(2)}] The curves $\tl{l}_i$ are free forever for $\tl{f}$.
\item[\textit{(4)}] $\Omega(\tl{f})\subset \cup_i \Theta_i$.

\end{enumerate}

We will prove items (2) and (4) assuming it holds the following weaker version of Item $1$: 
\begin{enumerate}
\item [\textit{(1$^*$)}] One of the sets $\Theta_i$ is essential.
\end{enumerate}

\begin{remark}    \label{remark.essential}
By the construction of the curves $\tl{l}_i$, if one of the sets $\Theta_i$ is essential then it is necessarily annular and vertical. 
\end{remark}

$\mathbf{(1^*)\Rightarrow (2).}$ Let $i_0$ such that $\Theta_{i_0}$ is essential, and therefore annular and vertical. Let $C\subset\R^2$ be a connected component of $\pi'^{-1}(\Theta_{i_0})$. Then $T_2^n(C)=C$ for all $n$, and $\pi'^{-1}(\Theta_{i_0})=\cup_n T_1^n(C)$. By the fact that $\Theta_{i_0}$ is annular, essential and vertical, $C$ is a connected set such that $\R^2\setminus C$ has two unbounded connected components. By our hypothesis that $\rho(f)=\{0\}\times I$ we can easily deduce that for each $n$, $T_1^n(C)$ is $f$-invariant. 

Let $S\subset\R^2$ be the $f$-invariant open strip bounded by $C$ and $T_1(C)$. In the construction of the curves $\tl{l}_i$, we saw that $\tl{l}_i$ is disjoint from $\Theta_j$ for every $i$ and $j$. By this, and since the curves $\tl{l}_i$ are essential, we have that for any $i$ there is exactly one lift $\ell_i\subset\R^2$ of the curve $\tl{l}_i$ which intersects $S$, and actually $\ell_i\subset S$. By the invariance of $S$ we have that $f^n(\ell_i)\subset S$ for all $n$ and $i$, and then $f^n(\ell_i)\cap T_1^m(\ell_i)=\empt$ for every $i$ and for all $n\in\Z$ and $m\neq 0$. Also, by the construction of $\tl{l}_i$, the lifts $\ell_i$ are Brouwer curves for $f$, and then $f^n(\ell_i)\cap\ell_i=\empt$ for every $i$ and all $n\in\Z$. We conclude that $f^n(\ell_i)\cap T_1^m(\ell_i)=\empt$ for every $i$, and all $n,m\in\Z$. This implies that the curves $\tl{l}_i$ are free forever for $\tl{f}$, and it holds item $(2)$ from Theorem A.\\
\\
$\mathbf{(1^*)\Rightarrow (4).}$ To prove it holds item $(4)$ we let $\tl{x}\in\T^2\setminus \cup_i \Theta_i$, and we will prove that $\tl{x}$ is wandering. Let $x\in\R^2$ be a lift of $\tl{x}$ contained in the open strip $S$ defined above. By the invariance of $S$, we have that if $B\subset S$ is a ball containing $x$, then the iterates of $B$ do not meet any integer translate of $B$ that is outside $S$, that is, $f^n(B)\cap T_1^r T_2^s(B)=\empt$ for all $r\neq 0$ and $n,s\in\Z$. Then, to get that $\tl{x}$ is wandering for $\tl{f}$ it suffices to show that there is an open ball $B'\subset S$ containing $x$ such that $f^n(B')\cap T_2^r(B')=\empt$ for all $n > 0$ and $r\in\Z$.

Let $i_1$ be an integer such that $x\in [\ell_{i_1},\ell_{i_1+1}]$. By hypothesis, $\tl{x}=\pi'(x)\notin\Theta_{i_1}$, and then there exists $n_0\in\Z$ such that $\tl{x}$ is contained in the open annulus bounded by $\tl{f}^{n_0}(l_{i_1})$ and $\tl{f}^{n_0+2}(l_{i_1})$. Then, $x$ is contained in the open strip $S'\subset\R^2$ bounded by $f^{n_0}(\ell_{i_1})$ and $f^{n_0+2}(\ell_{i_1})$. As the curve $\ell_{i_1}$ is a Brouwer curve for $f$ (by construction of the curves $\tl{l}_i$), then 
\begin{equation} \label{art1}
f^n(S')\cap S'=\empt  \ \ \ \txt{ for all $n\geq 2$}.
\end{equation}
Also, as the curve $\ell_{i_1}$ is a lift of a vertical curve in $\T^2$, we have that if $B'\subset S'$ is an open ball containing $x$ then $T_2^j(B')\subset S'$ for all $j\in\Z$. By this and by (\ref{art1}) we have that $f^n(B')\cap T_2^j(B')=\empt$ for all $n\geq 2$ and $j\in\Z$, and therefore there is a ball $B''\subset B'$ containing $x$ such that $f^n(B'')\cap T_2^j(B'')=\empt$ for all $n > 0$, which by last paragraph implies that $\tl{x}=\pi'(x)$ is wandering for $\tl{f}$. As $\tl{x}\in\T^2\minus\cup_i \Theta_i$ was arbitrary, we get that $\Omega(\tl{f})\subset\cup_i\Theta_i$, and it holds item $(4)$ from Theorem A, as desired.

\section{Proof of Theorem A, part II: proof of item (1).}         \label{sec.teoa2}

For the sake of completeness, we begin with the following remark about Theorem A. The fact that the rotation set is a \textit{non-degenerate} interval is essential in order to obtain annular dynamics (i. e., an invariant essential annular set). To illustrate this, we note that one can construct an example of a homeomorphism $\tl{f}$ of $\T^2$ with a lift $f$ to $\R^2$ whose rotation set consists of a point, and such that $f$ has no invariant essential annular sets. We briefly outline such a construction (for a more detailed description, see for ex. \cite{handel}). Let $v\in\R^2$ be a vector with irrational slope, and let $\chi$ denote the constant vector field $\chi\equiv v$ in $\T^2$. For $p\in\T^2$, let $\psi:\T^2\ra\R$ be a continuous function such that $\psi \geq 0$ and $\psi(x)=0$ if and only if $x=p$. Now, let $\tl{f}:\T^2\ra\T^2$ be the time-$1$ map of the flow given by the vector field $\psi\chi$. We have that $\txt{Fix}(\tl{f})=\{p\}$, and that the future orbit of every point passes arbitrarily close of $p$. If $\psi$ is chosen also to be close enough to zero near $p$, one can prove that every point has rotation vector equal to $(0,0)$, and then $\rho(f)=(0,0)$, for a lift $f:\R^2\ra\R^2$ of $\tl{f}$. Such a lift $f$ also has unbounded orbits both in the horizontal and vertical direction, which impedes the existence of any invariant annular essential set for $\tl{f}$. Other examples of homeomorphisms whose rotation set reduces to a vector with rational slope and have dynamics far from being annular can be found in \cite{korotal}.         

We now return the proof of Theorem A.

\subsection{It suffices to prove $(1^*)$.}          \label{sec.1b}

We recall that in Theorem A, for each $i$, the set $\Theta_i$ is the maximal invariant set of $[\tl{l}_i,\tl{l}_{i+1}]$ for $\tl{f}$. Also, we recall item (1) from that theorem:
\begin{enumerate}
\item[\textit{(1)}] At least one of the sets $\Theta_i$ is annular, essential, vertical and a semi-attractor for $\tl{f}$.
\end{enumerate}
In Section \ref{sec.24from1} we proved that the curves $\tl{l}_i$ from Theorem A are free forever for $\tl{f}$, if it holds the following:
\begin{enumerate}
\item[\textit{(1$^*$)}] At least one of the sets $\Theta_i$ is essential.
\end{enumerate}
As we mentioned in Remark \ref{remark.essential}, if one of the sets $\Theta_i$ is essential, then it is also annular and vertical. Now, observe that if the curves $\tl{l}_i$ are free forever for $\tl{f}$, and if $\Theta_{i_0}$ is essential, annular and vertical for some $i_0$, we easily get that $\Theta_{i_0}$ must be a semi-attractor for $\tl{f}$. 

Therefore, to prove item $(1)$ from Theorem A it suffices to prove it holds $(1^*)$. The rest of this section will be devoted to the proof of $(1^*)$.

\subsection{Strategy and outline of the proof of $(1^*)$.}         \label{sec.outline}

To prove $(1^*)$ we will work in $\R^2$ with lifts of the curves $\tl{l}_i$, so we start by fixing a family of such lifts.

\begin{figure}[h] 
\begin{center} 
\includegraphics{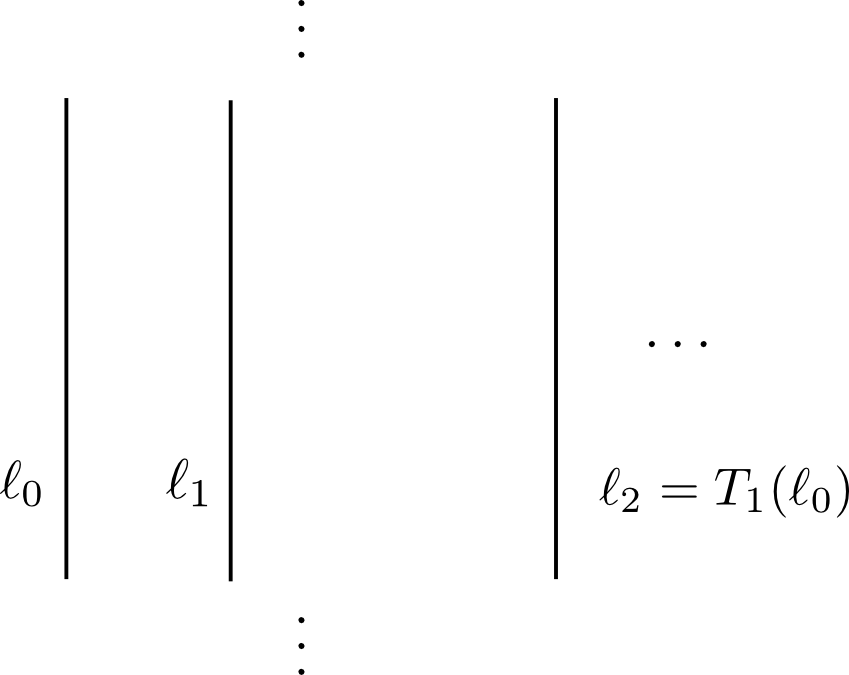}
\caption{The curves $\ell_i$.}
\label{fig.eles}  
\end{center}  
\end{figure}

\begin{definicion}[The curves $\ell_i$]              \label{def.eles}
For $i\in\N_0$ we define a lift $\ell_i\subset\R^2$ of the curve $\tl{l}_{i \mod 2}$ in the following way. First define $\ell_0\subset\R^2$ as any lift of $\tl{l}_0$. Then, define $\ell_1$ as the lift of $\tl{l}_1$ such that $\ell_{0}\prec\ell_1\prec T_1(\ell_0)$.
Then, for $i=0,1$, and $j\in\N$ we define $\ell_{2j+i}= T_1^{j}\ell_{i}$ (see Fig. \ref{fig.eles}).
\end{definicion}

In order to simplify some of the subsequent arguments, without loss of generality we make the following assumption. 

\begin{asuncion}     \label{asuncion2}
For every $i$, the curve $l_i$ a straight, vertical line. 
\end{asuncion}

\textbf{Strategy.} The proof of $(1^*)$ will be by contradiction. We will suppose that none of the sets $\Theta_i$ is essential, and we will obtain that $\max|\txt{pr}_1(\rho(f))|>0$, which contradicts our hypothesis that $\rho(f)$ is an interval of the form $\{0\}\times I$. 

We start by noting the following.

\begin{claim} \label{art10}
If none of the sets $\Theta_i$ is essential, then there is $n_0\in\N$ such that either 
\begin{enumerate}
\item $f^{n_0}(\ell_i)\cap\ell_{i+1}\neq\empt$ for all $i\geq 0$, or
\item $f^{n_0}(\ell_i)\cap\ell_{i-1}\neq\empt$ for all $i>0$.
\end{enumerate}
\end{claim}
\begin{prueba}
If there was $i_0$ such that $f^n(\ell_{i_0})\cap\ell_{i_0+1}=\empt$ for all $n\in\Z$, then the maximal invariant set of $(\ell_{i_0},\ell_{i_0+1})$, which is not empty (see Remark \ref{remark.maxinv}), would be a connected set $A$ such that $\pi'(A)=\Theta_{i_0}$ is essential. Details are left to the reader. 
\end{prueba}

Therefore, if we assume that none of the sets $\Theta_i$ is essential, we are in case 1 or 2 from Claim \ref{art10}. From now on, without loss of generality we will assume that we are in case 1. Moreover, we may assume also that $n_0=1$ (indeed, if we get positive horizontal speed for $f^{n_0}$ we also get this for $f$). We state this explicitly:

\begin{asuncion}             \label{asuncion3}
It holds $f(\ell_i)\cap\ell_{i+1}\neq\empt$ for all $i\geq 0$. In particular $\ell_i\prec f(\ell_i)$ for all $i\geq 0$.
\end{asuncion}

\textbf{Outline of the proof of $(1^*)$.}

We will see that in each strip $[\ell_i,\ell_{i+1}]$ there exists a family $\cl{F}$ of compact, connected sets, with the property that 
$$f^n(F)\subset[\ell_i,\ell_{i+1}] \ \ \ \ \ \ \ \ \forall\, F\in \cl{F}, \ \forall\, n\geq 0,$$ 
(cf. Definition \ref{def.ridesanchors} and Fig. \ref{fig.anchors}). The sets $F\in\cl{F}$ will be called \textit{anchors}.

Under Assumption \ref{asuncion3} we will prove by induction the following:
\begin{itemize} 
\item[$\star$] There exists an integer $N_*>0$ such that, for all $n\geq 0$, $f^{N_*(n+1)}(\ell_0)$ has `good' intersection with the anchors of the $n$-th strip $[\ell_{n},\ell_{n+1}]$ and with the curve $\ell_{n+1}$. 
\end{itemize}
By `good' intersection we mean roughly that the curve $f^{N_*(n+1)}(\ell_0)$ contains an arc $\gamma$ such that $\gamma\subset [\ell_n,\ell_{n+1}]$, $\gamma(0)$ lies in an anchor, and $\gamma(1)\in\ell_{n+1}$ (see Definition \ref{def.good2}, and Fig. \ref{fig.goodintersection}). 

In this way, $(\star)$ gives us in particular that for all $n>0$, $f^{N_* n}(\ell_0)\cap \ell_n\neq\empt$, and then $\max\txt{pr}_1(\rho(f))>0$, which is the desired contradiction. We now give an idea of the proof of $(\star)$. 

Assume that an iterate of $\ell_0$ contains an arc $\gamma_{n-1}$ which has `good' intersection with an anchor of $[\ell_{n-1},\ell_n]$ and with the curve $\ell_n$, as described above. We first sketch how we will prove there is $N'_*>0$ such that $f^{N'_*}(\gamma_{n-1})$ intersects $\ell_{n+1}$.

\begin{figure}[h] 
\begin{center} 
\includegraphics{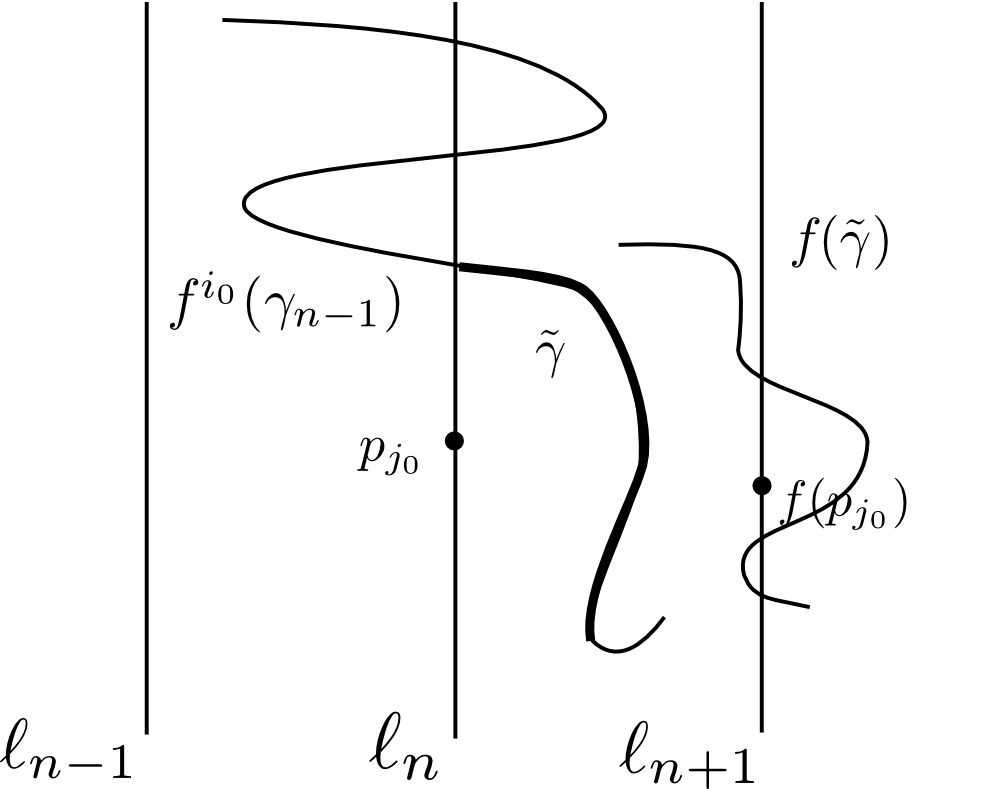}
\caption{The iterates of $\gamma_{n-1}$ get `stretched' in $R(\ell_n)$, and they get `pushed to the right' by $p_{j_0}$.}
\label{fig.outline}  
\end{center}  
\end{figure}

By Proposition \ref{prop.prel} and by Remark \ref{remark.prop.prel} we have that the maximal invariant sets of the strips $[\ell_j,\ell_{j+1}]$ move either upwards or downwards (these maximal invariant sets are non-empty, see Remark \ref{remark.maxinv}). Suppose that the maximal invariant set of $[\ell_{n-1},\ell_n]$ moves upwards, and that the maximal invariant set of $[\ell_n,\ell_{n+1}]$ moves downwards. As the arc $\gamma_{n-1}$ is such that $\gamma_{n-1}(0)$ belongs to an anchor, and $\gamma_{n-1}(1)\in\ell_n$, we have that $f^i(\gamma_{n-1}(0))\in L(\ell_n)$ and $f^i(\gamma_{n-1}(1))\in R(\ell_n)$, for all $i > 0$. We will see that, as long as the iterates $f^i(\gamma_{n-1})$ of $\gamma_{n-1}$ do not intersect $\ell_{n+1}$, they will get `stretched' vertically in $[\ell_{n-1},\ell_{n+1}]$, that is, the vertical diameter of the iterates $f^i(\gamma_{n-1})$ will grow. Moreover, we will see that this `stretching' occurs in $R(\ell_n)$, that is, there will be an iterate $f^{i_0}(\gamma_{n-1})$ and a subarc $\tl{\gamma}$ of it contained in $(\ell_n,\ell_{n+1})$ with `large' vertical diameter (see Fig. \ref{fig.outline}). On the other hand, by Assumption \ref{asuncion3} there are points $p_j$ in $\ell_n$ such that $f(p_j)\in \ell_{n+1}$. If the vertical diameter of the arc $\tl{\gamma}$ is large enough, we will see that $\tl{\gamma}$ is `to the right' of one of the points $p_j$, say $p_{j_0}$, and that $\tl{\gamma}$ will get `pushed' to the right by $p_{j_0}$; that is, $f(\tl{\gamma})\subset f^{i_0+1}(\gamma_{n-1})$ will intersect $\ell_{n+1}$ (see Fig. \ref{fig.outline}). A precise meaning of $\tl{\gamma}$ being `to the right' of a point $p_j$ and getting `pushed' to the right by it will be given in Section \ref{sec.prpl}.

\begin{figure}[h] 
\begin{center} 
\includegraphics{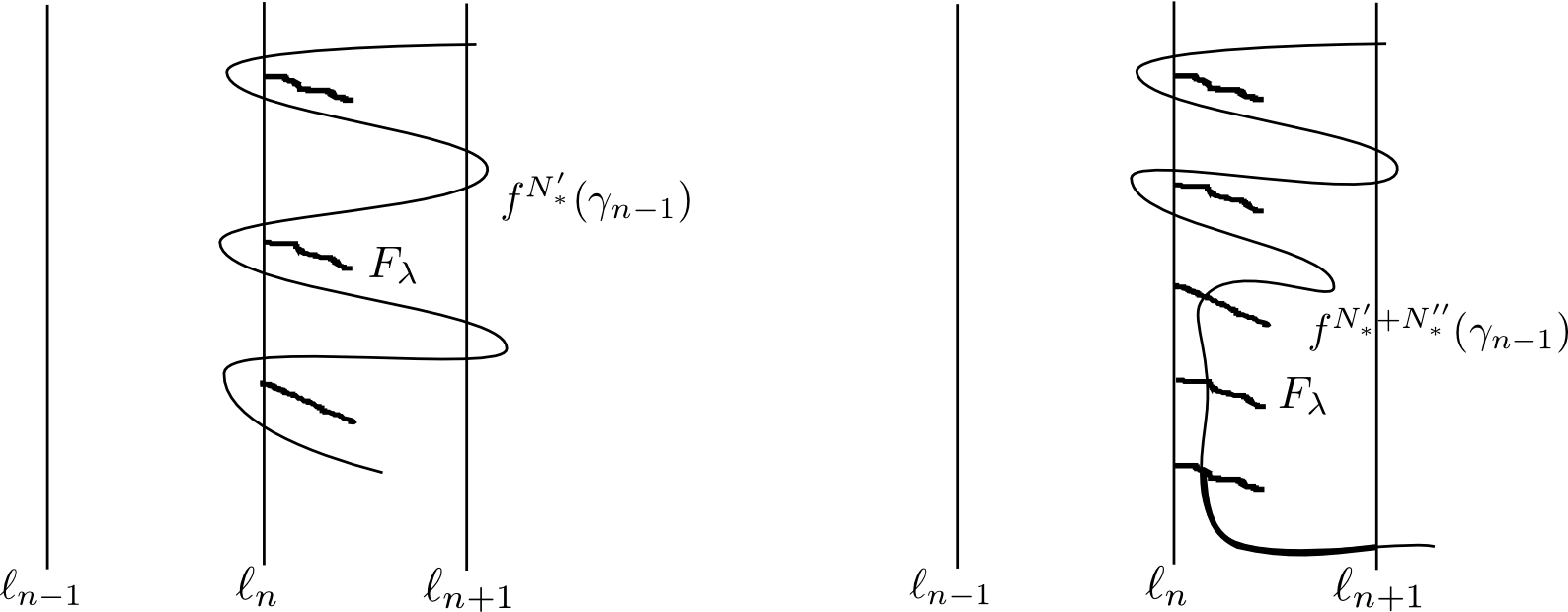}
\caption{$f^{N'_*}(\gamma_{n-1})$ might not intersect the anchors $F_{\lambda}$ of $[\ell_n,\ell_{n+1}]$, but $f^{N'_*+N''_*}(\gamma_{n-1})$ does have good intersection with the anchors and with $\ell_{i+2}$.}
\label{fig.outline2}  
\end{center}  
\end{figure}

In this way we get that, if $N'_*=i_0+1$, $f^{N'_*}(\gamma_{n-1})$ intersects $\ell_{n+1}$. However, the arc $f^{N'_*}(\gamma_{n-1})$ could have no good intersection with the anchors and $\ell_{n+1}$ (the arc $f^{N'_*}(\gamma_{n-1})$ might not even intersect the anchors of $[\ell_n,\ell_{n+1}]$, see Fig. \ref{fig.outline2}). With more work, using the fact that the iterates of $f^{N'_*}(\gamma_{n-1})$ will get `stretched' in $R(\ell_n)$, we will see that they cannot always avoid the anchors of $[\ell_n,\ell_{n+1}]$, and that there is $N''_*>0$ such that $f^{N'_*+N''_*}(\gamma_{n-1})$ does have good intersection with the anchors of $[\ell_n,\ell_{n+1}]$ and with $\ell_{n+1}$. The numbers $N'_*$ and $N''_*$ will be shown to be independent of $n$, and in this way we will get an integer $N_*=N'_*+N''_*$, independent of $n$, such that $f^{N_*}(\gamma_{n-1})$ has good intersection with the anchors of $[\ell_n,\ell_{n+1}]$ and with $\ell_{n+1}$, which will give us the induction step of $(\star)$.

\subsection{Being \textit{to the right} and being \textit{pushed}.}          \label{sec.prpl}


In this section we make more precise the idea of an arc, or a continuum in general, being `to the right' (or `to the left') of a point in the plane. Also, we explain what we mean by a continuum being `pushed' to the right (or to the left) by a point. 



The idea of a continuum being `to the right' or `to the left' of a point is represented by the following properties $\mathbf{PR}$ and $\mathbf{PL}$.

\begin{figure}[h] 
\begin{center} 
\includegraphics{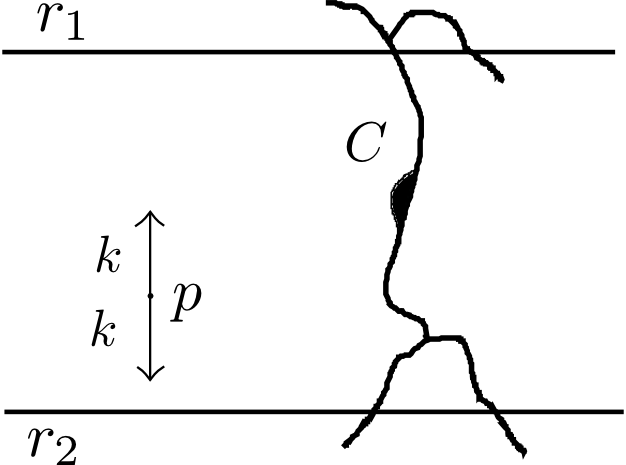}
\includegraphics{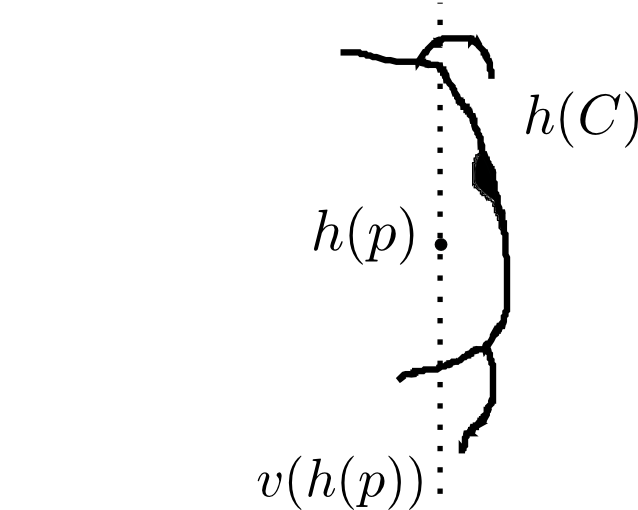}
\caption{Left: a set $C$ satisfying property $\mathbf{PL}(k,p)$. Right: $h(C)\cap R(v(h(p)))\neq\emptyset$.}
\label{fig.propp}  
\end{center}  
\end{figure}

\begin{definicion}
Let $\tl{h}:\T^2\ra\T^2$ be a homeomorphism isotopic to the identity, and let $h:\R^2\ra\R^2$ be a lift of  $\tl{h}$. Let $C\subset \R^2$ be a continuum, $k\in\R^+$ and $p\in\R^2$. We say that $C$ satisfies the property $\mathbf{PL}(k,p)$ if the following hold (see Fig. \ref{fig.propp}):
\begin{enumerate}
\item There exist horizontal (disjoint) straight lines $r_1 \prec r_2$, (oriented as going to the right) such that $r_1\cap C\neq\emptyset$, $r_2\cap C\neq\emptyset$, and such that the strip $(r_1,r_2)\subset\R^2$ contains a ball of radius $k$ centered in $p$.
\item The point $p$ belongs to the (unique) connected component of $(r_1,r_2)\setminus C$ which is unbounded to the left. 
\end{enumerate}
Analogously, we say that $C$ satisfies the property $\mathbf{PR}(k,p)$ if it holds item (1) from property $\mathbf{PL}(k,p)$ and $p$ belongs to the (unique) connected component of $(r_1,r_2)\setminus C$ which is unbounded to the right. 
\end{definicion}

The following lemma will be an important tool in the proof of Theorem A. It gives a meaning to the following idea: if $h:\R^2\ra\R^2$ is a lift of a torus homeomorphism isotopic to the identity, if $C$ is a continuum which is `to the right' (or `to the left') of a point $p\in\R^2$, and if $C$ is `large enough', then, applying $h$, $p$ `pushes $C$ to the right (or to the left) of $h(p)$'.

\begin{lema} \label{propp}
Let $\tl{h}:\T^2\ra\T^2$ be a homeomorphism isotopic to the identity, let $h:\R^2\ra\R^2$ be a lift of $\tl{h}$, and for $x\in\R^2$, denote by $v(x)\subset\R^2$ the vertical straight line that passes through $x$. There exists $k>0$ such that if a compact connected set $C\subset\R^2$ satisfies $\mathbf{PL}(k,p)$ (resp. $\mathbf{PR}(k,p)$) for some $p\in\R^2$, then $h(C)\cap R(v(h(p)))\neq\emptyset$ (resp. $h(C)\cap L(v(h(p)))\neq\emptyset$, see Fig. \ref{fig.propp}). 
\end{lema}
\begin{prueba}
First observe that as $h$ is the lift of a homeomorphism of $\T^2$, $\left\| h-\txt{Id}\right\|_{0}<\infty$. Define $k=2\left\| h-\txt{Id}\right\|_0+1$. Suppose that $C$ satisfies the property $\mathbf{PL}(k,p)$ for some $p\in\R^2$ (the case of $\mathbf{PR}(k,p)$ is similar). Then there are two horizontal straight lines $r_1\prec r_2$ intersecting $C$ and such that $(r_1,r_2)$ contains a ball of radius $k$ centered in $p$, and $p$ belongs to the connected component $U_L$ of $(r_1,r_2)\setminus C$ which is unbounded to the left. Observe that by the definition of $k$, $\min \txt{pr}_2(h(r_1)) > h(p)_2 >\max \txt{pr}_2(h(r_2))$, and then if $w$ is the horizontal straight line passing through $h(p)$, we have 
\begin{equation}  \label{j3}
w\subset (h(r_1),h(r_2)).
\end{equation} 
Let $U_R$ be the connected component of $(r_1,r_2)\setminus C$ which is unbounded to the right. As $\left\| h- \txt{Id}\right\|_{0}<\infty$, $h(U_L)$ is unbounded to the left and bounded to the right, and also $h(U_R)$ is unbounded to the right and bounded to the left.

We claim that for this choice of $k$, we have $h(C)\cap R(v(h(p)))\neq\emptyset$. If this was not the case, then we would have that $C\cap w_+=\emptyset$, where $w_+=w\cap R(v(h(p))$. By (\ref{j3}), $w_+\subset (h(r_1),h(r_2))$, and therefore $w_+$ is contained in $h(U_R)$. Then $h(p)$ belongs to $h(U_R)$, which is unbounded to the right, which contradicts the fact that $p$ belongs to the connected component $U_L$ of $(r_1,r_2)\setminus C$ which is bounded to the right. We must have then that $h(C)\cap R(v(h(p)))\neq\emptyset$, and this proves the lemma.     
\end{prueba}

The following lemma relates the properties $\mathbf{PL}$ and $\mathbf{PR}$ and the curves $\ell_i$.

\begin{lema} \label{isotopia2} 
Let $i,j\in\N$ and suppose that $n\in\Z$ is such that $f^{n}(\ell_{i})\cap \ell_{j}\neq\emptyset$. Then, there exists a constant $K>0$ such that, if $C\subset\R^2$ is a continuum contained in the open strip bounded by $\ell_i$ and $\ell_j$ and such that $\txt{diam}_2(C)\geq K$, then:

\begin{itemize}
\item If $i<j$, then $f^n(C)\cap R(\ell_{j})\neq\emptyset$.
\item If $j<i$, then $f^n(C)\cap L(\ell_{j})\neq\emptyset$.
\end{itemize}
\end{lema}

\begin{prueba}
Without loss of generality suppose that $n>0$. By Lemma \ref{propp} applied to $f^n$ there is a constant $k>0$ such that, if $C$ is a continuum that satisfies the property $\mathbf{PL}(k,p)$ (or $\mathbf{PR}(k,p)$) for some $p\in\R^2$, then $f^n(C)\cap R(v(f^n(p)))\neq\empt$ (resp., $f^n(C)\cap L(v(f^n(p)))\neq\empt$). 

We treat the case $i<j$, the case $i>j$ being similar. By hypothesis $f^n(\ell_i)\cap \ell_j\neq\empt$. Take $x\in f^{-n}(\ell_j)\cap\ell_i$ and define $K=k+1$. Suppose that $C$ is a continuum contained in $(\ell_i,\ell_j)$ and such that $\txt{diam}_2(C)\geq K$. Then there is $s\in\Z$ such that 
$$((T_2^s (x))_2-k,(T_2^s(x))_2+k) \subset \txt{pr}_2(C).$$ 
As we assumed that the lines $\ell_i$ are straight, and as $C\subset R(\ell_i)$, it is easy to see that $C$ satisfies property $\mathbf{PL}(k,x)$. Therefore as $f^n(T_2^s(x))\in \ell_j$, Lemma \ref{propp} gives us that $f^n(C)\cap R(\ell_j) = f^n(C)\cap R(v(f^n(T_2^s(x))))\neq\emptyset$, as we wanted. 
\end{prueba}

\subsection{The anchors and good intersection.}

We start with some preliminary constructions.

\subsubsection{The sets $L_{\infty}^i, R_{\infty}^i$, and $X_i$.}        \label{sec.linf}

\begin{figure}[h] 
\begin{center} 
\includegraphics{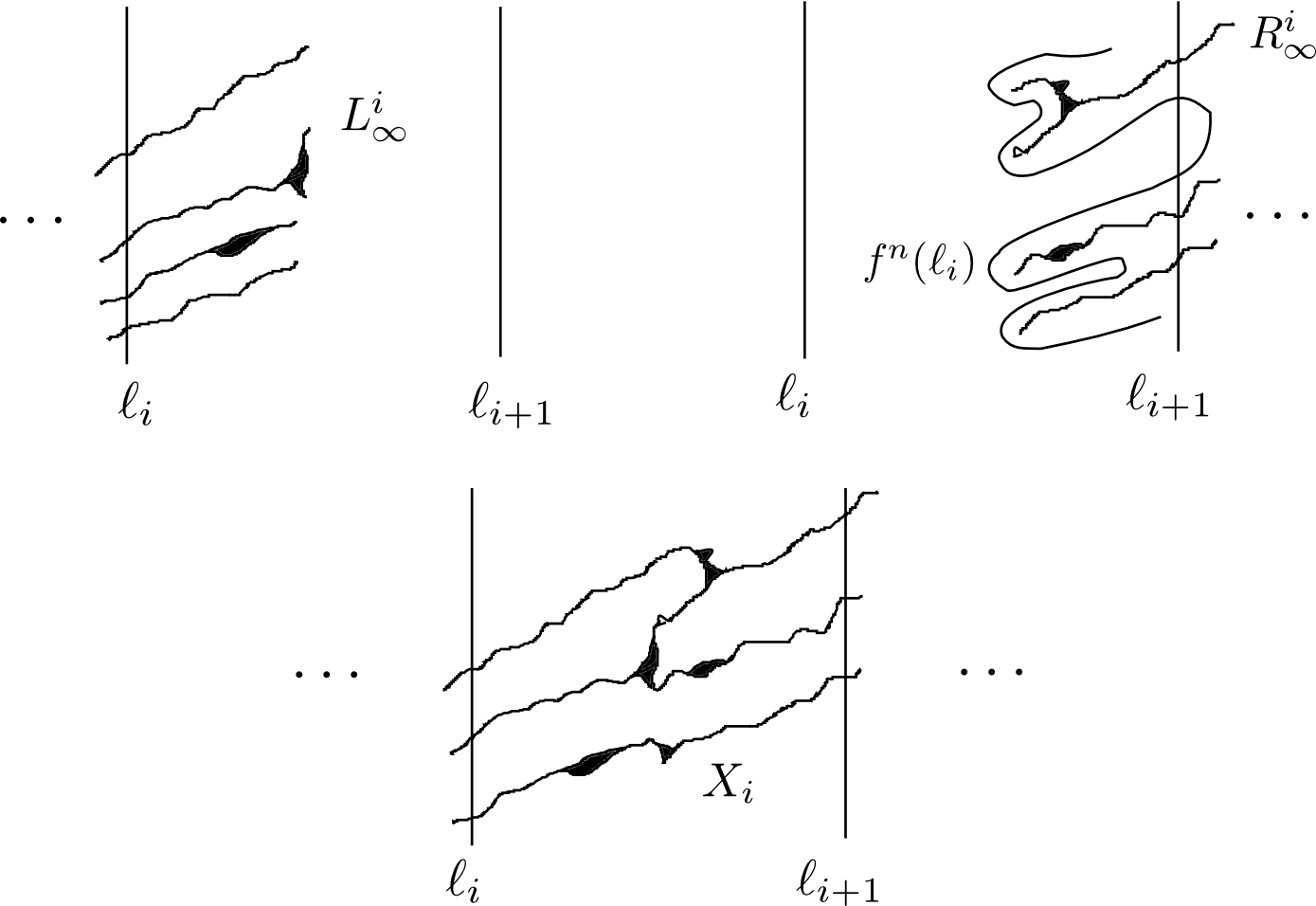}
\caption{Some examples of the sets $L_{\infty}^i$, $R_{\infty}^i$ and $X_i$.}
\label{fig.xi}  
\end{center}  
\end{figure}

For each $i\in\N$, we define the sets $L_{\infty}^i$ and $R_{\infty}^i$, which in some sense are the `stable' and `unstable' sets (resp.) of the maximal invariant set in $[\ell_i,\ell_{i+1}]$ for $f$. Let
$$R_{\infty}^i= \bigcap_{n\in\Z} R \left(f^n(\ell_i)\right), \ \ \    L_{\infty}^i=\bigcap_{n\in\Z} L (f^{-n}(\ell_{i+1})), \ \ \txt{and} \ \ \ X_i=L_{\infty}^i\cup R_{\infty}^i$$
(see Fig. \ref{fig.xi}.) 

By definition, the sets $R_{\infty}^i, L_{\infty}^i$ and $X_i$ are $f$-invariant. As we are under Assumption \ref{asuncion3}, $\ell_i\prec f(\ell_i)$ for all $i$, and therefore we have that 
$$R_{\infty}^i= \{x\in \R^2\ : \  f^{-n}(x)\in R(\ell_i) \ \forall n\geq 0\},$$
and
$$L_{\infty}^i= \{x\in \R^2\ : \  f^{n}(x)\in L(\ell_{i+1}) \ \forall n\geq 0\},$$
Therefore, for each $i$, the set $R_{\infty}^i\cap L_{\infty}^i$ is the maximal invariant set of $[\ell_i,\ell_{i+1}]$ for $f$ (which is non-empty by Remark \ref{remark.maxinv}), and then
$$R_{\infty}^i\cap[\ell_i,\ell_{i+1}]= \{x\in[\ell_i,\ell_{i+1}]\ : \ d(f^{-n}(x),L_{\infty}^i\cap R_{\infty}^i)\ra 0 \txt{ as } n\ra\infty\},$$
and
$$L_{\infty}^i\cap[\ell_i,\ell_{i+1}]=\{x\in [\ell_i,\ell_{i+1}] \ : \ d(f^{n}(x),L_{\infty}^i\cap R_{\infty}^i)\ra 0 \txt{ as } n\ra\infty\}.$$
That is, the set $L_{\infty}^i\cap[\ell_i,\ell_{i+1}]$ can be thought as the `local stable set' of $R_{\infty}^i\cap L_{\infty}^i$, and $R_{\infty}^i\cap[\ell_i,\ell_{i+1}]$ can be thought as the `local unstable set' of $R_{\infty}^i\cap L_{\infty}^i$. The following lemmas study some properties of these sets.

\begin{lema}\label{discos}
For every $i\geq 0$:
\begin{enumerate}
\item if $C$ is a connected component of $R_{\infty}^i$, then $\sup \txt{pr}_1(C)=+\infty$,
\item if $C'$ is a connected component of $L_{\infty}^i$, then $\inf \txt{pr}_1(C')=-\infty$, and
\item the connected components of $\R^2\setminus X_i$ are simply connected.
\end{enumerate}
\end{lema} 

\begin{prueba}
Let $S=\R\times S^1\cup\{\infty\}\cup\{-\infty\}$ be the two-point compactification of $\R\times S^1$, which is homeomorphic to $S^2$, and let $j:\R\times S^1\hookrightarrow S$ be the inclusion. The curves $\pi(\ell_i)\subset \R\times S^1$ are vertical circles, and then the sets $D_n:=j(\pi(\ol{R}(f^n(\ell_i))))\cup\{\infty\}$, and $D_n':=j(\pi(\ol{L}(f^{-n}(\ell_i))))\cup\{-\infty\}$ are topological closed discs in $S$, for any $n$ and $i$. Observe that
$$\wh{L}_i:=j(\pi(L_{\infty}^i))\cup \{-\infty\}=\bigcap_{n\in\N} D_n',$$ 
and 
$$\wh{R}_i:=j(\pi(R_{\infty}^i))\cup\{\infty\}=\bigcap_{n\in\N} D_n.$$ 
As we are under Assumption \ref{asuncion3}, $\ell_i\prec f(\ell_{i})$ for any $i$, and then $D_{n+1}\subset D_n$ for all $n$. Therefore the sets $\wh{L}_i$ and $\wh{R}_i$ are nested intersections of topological closed discs, and thus they are compact and connected. 

Observe that, for every $i$, $L_{\infty}^i\cap R_{\infty}^i= \cap_{n\in\Z}f^n((\ell_i,\ell_{i+1}))$. By Remark \ref{remark.maxinv}, $L_{\infty}^i\cap R_{\infty}^i\neq\empt$, and then $\wh{L}_i\cap\wh{R}_i\neq\empt$. Therefore, $\wh{X}_i:=j(X_i)\cup\{\infty\}\cup\{-\infty\}$ is compact and connected, as it is the union of $\wh{L}_i$ and $\wh{R}_i$, which are connected sets with nonempty intersection. \\
$(1)$. It suffices to show that, for $x\in \wh{R}_i\setminus\{\infty\}$, if $C_x$ is the connected component of $\wh{R}_i\setminus\{\infty\}$ containing $x$, then $\infty\in \ol{C_x}$. This result from plain topology is known as the `boundary bumping theorem', and can be found for example in \cite{kr} (Theorem 2, p. 172).\\
$(2)$. The proof is analogous to $(1)$.\\
$(3)$. First note that, as $\wh{X}_i$ is connected, then the connected components of $S\setminus \wh{X}_i$ are simply connected. Now, let $V$ be any connected component of $X_i^c$. Then $j(V)\subset S$ is a connected component of $S\setminus\wh{X}_i$, and therefore simply connected. As $j:\R^2\ra S\setminus \{\infty\}$ is a homeomorphism, $V$ must be also simply connected.
\end{prueba}

\begin{corolario}  \label{uisc}
For each $i\geq 0$, the connected components of $X_i^c\cap(\ell_i,\ell_{i+1})$ are simply connected.
\end{corolario}
\begin{prueba}
This is an easy consequence of the fact that the connected components of the intersection of two simply connected sets in the plane are simply connected.
\end{prueba}

The following lemma is an application of Lemma \ref{isotopia2}.

\begin{lema} \label{propp2}
There exists a constant $M_0$ such that, for any $i\geq 0$, any connected component of $R_{\infty}^i\cap L(\ell_{i+1})$ has vertical diameter less than $M_0$, and also any connected component of $L_{\infty}^i \cap R(\ell_i)$ has vertical diameter less than $M_0$. 
\end{lema}
\begin{prueba}
First we treat the case of $R_{\infty}^i\cap L(\ell_{i+1})$. As we are under Assumption \ref{asuncion3}, we have that $f^{-1}(\ell_{i+1})\cap \ell_i\neq\emptyset$ for all $i\geq 0$. By Lemma \ref{isotopia2}, there exists a constant $K_0>0$ such that if $C\subset\R^2$ is a continuum contained in $(\ell_i,\ell_{i+1})$ with $\txt{diam}_2(C)>K_0$, then $f^{-1}(C)\cap L(\ell_i)\neq\emptyset$. Therefore, for any $i\geq 0$, any connected component $C_0$ of $R_{\infty}^i\cap L(\ell_{i+1})$ must have vertical diameter less than $K_0$, because otherwise $f^{-1}(C_0)$ would intersect $L(\ell_i)$, which contradicts the definition of $R_{\infty}^i$. 

From a symmetric argument, we obtain that any connected component of $L_{\infty}^i\cap R(\ell_i)$ must have vertical diameter less that $K_0$. Setting $M_0=K_0$, the lemma follows. 
\end{prueba}

\begin{lema} \label{limitadas}
There exists $M_1>0$ such that for any $i\geq 0$, any connected component of $X_i^c\cap(\ell_i,\ell_{i+1})$ has vertical diameter less than $M_1$.
\end{lema}
\begin{prueba}
Let $i\geq 0$, and let $x\in L_{\infty}^i\cap R_{\infty}^i$. Let $C_1$ and $C_2$ be the connected components of $R_{\infty}^i\cap L(\ell_{i+1})$ and $L_{\infty}^i\cap R(\ell_i)$, respectively, that contain $x$. By Lemma \ref{discos}, the connected component of $R_{\infty}^i$ that contains $C_1$ is unbounded to the right and the connected component of $L_{\infty}^{i}$ that contains $C_2$ is unbounded to the left. By this reason, the connected set $C=C_1\cup C_2$ separates $(\ell_i,\ell_{i+1})$, that is, $(\ell_i,\ell_{i+1})\minus C$ is not connected. Also, by Lemma \ref{propp2}, there is a constant $M_0$ such that $\txt{diam}_2(C_i)\leq M_0$ for $i=1,2$, and then $\txt{diam}_2(C)\leq 2M_0$. Thus, $C\cap T_2^{3M_0}(C)=C\cap T_2^{-3M_0}(C)=\empt$.

Now, consider the set 
$$A=\bigcup_{n\in\Z} T_2^{3M_0 n} (C).$$
The connected components of $(\ell_i,\ell_{i+1})\minus A$ have then vertical diameter less than $\txt{diam}_2(C)+ 3M_0\leq 5M_0$. As $A\subset X_i$, any connected component of $X_i^c\cap(\ell_i,\ell_{i+1})$ is contained in a connected component of $(\ell_i,\ell_{i+1})\minus A$, and therefore has diameter less than $5M_0$. Therefore, making $M_1:=5M_0$, the lemma follows. 
\end{prueba}

\subsubsection{Definition of the anchors and good intersection.}

\begin{figure}[h] 
\begin{center} 
\includegraphics{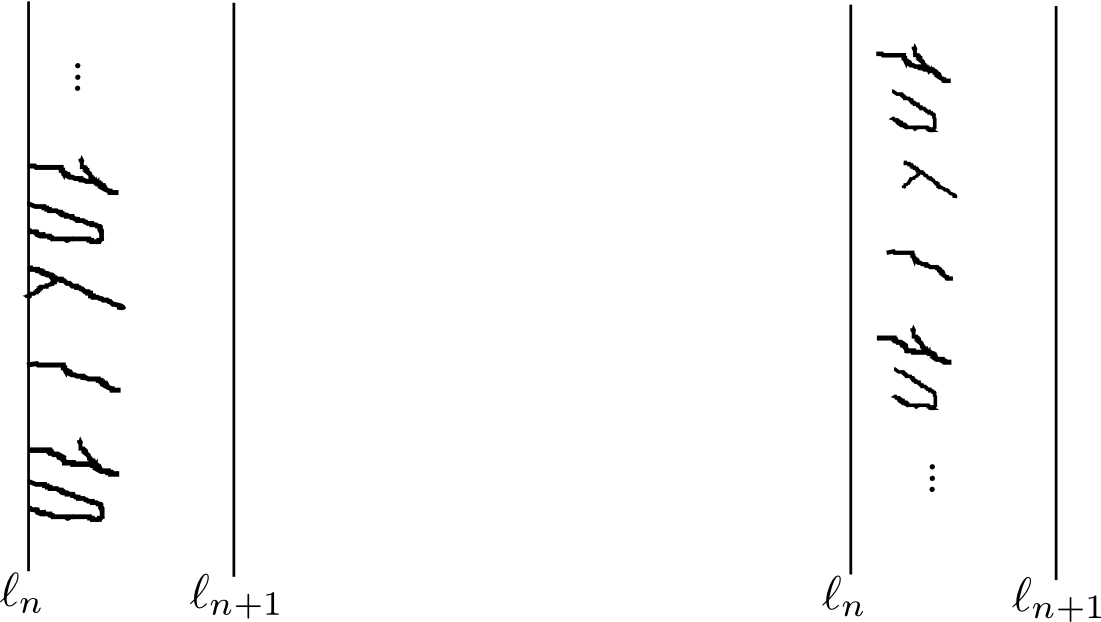}
\caption{The anchors and their images by $f$.}
\label{fig.anchors}  
\end{center}  
\end{figure}

We now are ready to make the main definitions of this section:

\begin{definicion} \label{def.ridesanchors}
For any $i\geq 0$, the connected components of the set
$$L_{\infty}^i\cap [\ell_i,\ell_{i+1}]$$
will be called \textbf{anchors}. For a fixed $i$, the set $L_{\infty}^i\cap [\ell_i,\ell_{i+1}]$ will be refered to as \textit{the anchors of the strip} $[\ell_i,\ell_{i+1}]$ (see Fig. \ref{fig.anchors}).
\end{definicion}

By definition, for any $i$, the set $L_{\infty}^i\cap [\ell_i,\ell_{i+1}]$ is forward $f$-invariant, and if a point belongs to an anchor of the strip $[\ell_i,\ell_{i+1}]$, all its future iterates by $f$ remain in $[\ell_i,\ell_{i+1}]$. We now define good intersection of a curve with the anchors and with one of the curves $\ell_i$.

\begin{figure}[h] 
\begin{center} 
\includegraphics{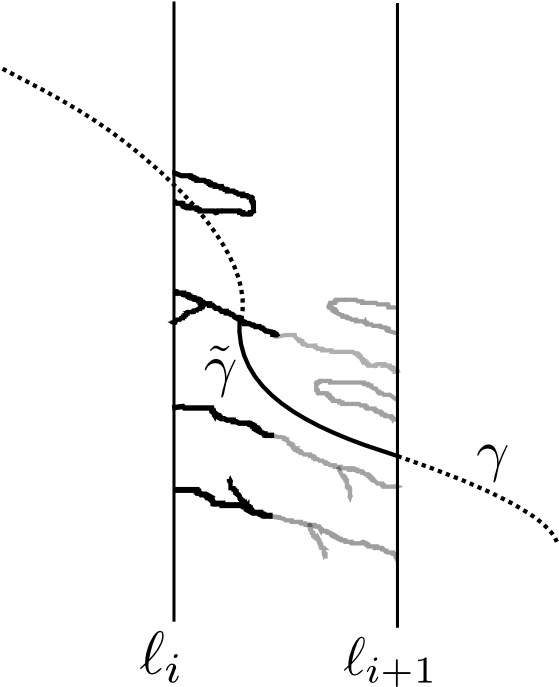}
\caption{Good intersection. In black are the anchors, and in gray appears the set $R_{\infty}^i\cap[\ell_i,\ell_{i+1}]$.}
\label{fig.goodintersection}  
\end{center}  
\end{figure}

\begin{definicion}             \label{def.good2}
Let $i\geq 0$. We say that a curve $\gamma$ has \textit{good intersection with the anchors of the strip $[\ell_i,\ell_{i+1}]$ and with $\ell_{i+1}$} if there is an arc $\tl{\gamma}\subset\gamma$ such that:
\begin{itemize}
\item $\tl{\gamma}(0)$ lies in an anchor of $[\ell_i,\ell_{i+1}]$,
\item $\tl{\gamma}(1)$ lies in the curve $\ell_{i+1}$, and
\item $\tl{\gamma}(t)$ belongs to $(X_i)^c \cap [\ell_i,\ell_{i+1}]$, for all $0<t<1$ (see Fig. \ref{fig.goodintersection}).
\end{itemize}
In this situation, we say that $\tl{\gamma}$ \textit{realizes the good intersection}.
\end{definicion}

\subsection{Main Lemma.}     \label{sec.mainlemma}

As we mentioned at the begining of this section, our proof of item $(1^*)$ is by contradiction. If item $(1^*)$ is not true then we can assume that $f(\ell_i)\cap\ell_{i+1}\neq\empt$ for all $i\geq 0$ (Assumption \ref{asuncion3}), and from this we want to prove that there is horizontal speed, namely, $\max \txt{pr}_1(\rho(f)) >0$, which contradicts our hypothesis that $\rho(f)$ is an interval of the form $\{0\}\times I$.

\begin{figure}[h] 
\begin{center} 
\includegraphics{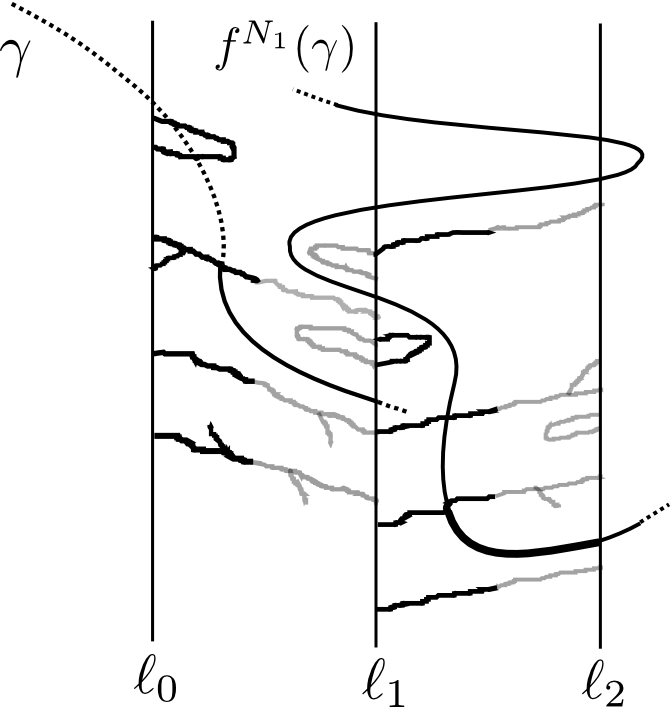}
\caption{Illustration for the Main Lemma.}
\label{fig.mainlemma}  
\end{center}  
\end{figure}

As we also explained, to find horizontal speed we will see that:
\begin{itemize} 
\item[$\star$] There exists an integer $N_*>0$ such that, for all $n\geq 0$, $f^{N_*(n+1)}(\ell_0)$ has good intersection with the anchors of the $n$-th strip $[\ell_{n},\ell_{n+1}]$ and with the curve $\ell_{n+1}$. 
\end{itemize}
In particular, $f^{N_*(n+1)}(\ell_0)$ intersects $\ell_{n+1}$ for all $n\geq 0$, and $\max \txt{pr}_2(\rho(f)) >0$. The proof of ($\star$) will be given by our Main Lemma:

\begin{lema}[Main Lemma] \label{N6}
There exists $N_1>0$ such that, if a curve $\gamma$ has good intersection with the anchors of the strip $[\ell_0,\ell_1]$ and with the curve $\ell_1$, then $f^{N_1}(\gamma)$ has good intersection with the anchors of the next strip $[\ell_1,\ell_2]$ and with the curve $\ell_2$ (see Fig. \ref{fig.mainlemma}).
\end{lema}

\subsection{The Main Lemma implies $(1^*)$}

We begin by making the following remark. The arguments in the proof of the Main Lemma \ref{N6} will automatically apply to every strip $[\ell_i,\ell_{i+1}]$, and this will give us that for any $i\geq 0$, there is a constant $N_1(i)>0$ such that if an arc $(\gamma)$ has good intersection with the anchors of $[\ell_i,\ell_{i+1}]$ and with $\ell_{i+1}$, then $f^{N_1(i)}(\gamma)$ has good intersection with the anchors of $[\ell_{i+1},\ell_{i+2}]$ and with $\ell_{i+2}$. As the curves $\ell_i$ come from a finite family of closed curves in $\T^2$, the numbers $N_1(i)$ will have a maximum $N_1^*$.  

With this in mind, to see that the Main Lemma implies ($\star$), and therefore $(1^*)$, it suffices to show the following:

\begin{claim}      \label{claim.star}
The curve $f(\ell_0)$ has good intersection with the anchors of $[\ell_0,\ell_1]$ and with $\ell_1$.
\end{claim}

Indeed, with this claim and the remark we just made, we get by induction that for each $n\geq 0$, $f^{N_1^*(n+1)}(\ell_0)$ has good intersection with the anchors of $[\ell_n,\ell_{n+1}]$ and with $\ell_{n+1}$. Setting $N_*= N_1^*$, this proves ($\star$). 

Therefore, the proof of Claim \ref{claim.star} and of the Main Lemma will conclude with the proof by contradiction of $(1^*)$. We proceed now to the easy proof of the claim.  

\begin{prueba}[Proof of Claim \ref{claim.star}]
We know by Lemma \ref{discos} that the connected components of $L_{\infty}^0$ are unbounded to the left, and in particular they intersect $\ell_0$. That is, there exists a point $x\in\ell_0$ which also belongs to an anchor of $[\ell_0,\ell_1]$. On the other hand, by Assumption \ref{asuncion3} we know there is $y\in\ell_0$ such that $f(y)\in\ell_1$. We claim that $f(\ell_0)$ has good intersection with the anchors of $[\ell_0,\ell_1]$ and with $\ell_1$. To see this, let $\gamma_1\subset\ell_0$ be the arc from $x$ to $y$. Then, as the union of the anchors of $[\ell_0,\ell_1]$ is forward $f$-invariant, $f(\gamma_1)$ is an arc in $R(\ell_0)$ which goes from an anchor of $[\ell_0,\ell_1]$ to $\ell_1$. Observe that as $\gamma_1\subset\ell_0$ and by the definition of $R_{\infty}^0$,
\begin{equation}    \label{eq1}
f(\gamma_1)\cap R_{\infty}^0=\empt.
\end{equation}
Let $\gamma_2\subset\gamma_1$ be an arc which is minimal with respect to the property of having one endpoint in an anchor of $[\ell_0,\ell_1]$ and the other endpoint in $\ell_1$. Then:
\begin{itemize}
\item $\gamma_2$ is a non-degenerate arc, as $L_{\infty}^0\cap\ell_1=\empt$ (by definition of $L_{\infty}^0)$,
\item $\gamma_2\subset [\ell_0,\ell_1]$ (by the minimality of $\gamma_2$), and
\item $\gamma_2(t)\notin X^0$ for $0<t<1$ (by \ref{eq1}). 
\end{itemize}
This means that $f(\gamma_1)\subset f(\ell_0)$ has good intersection with the anchors of $[\ell_0,\ell_1]$ and with $\ell_1$. This proves our claim.
\end{prueba}

\subsection{Proof of Main Lemma.}    \label{pruebamainlemma}

We will prove first two preliminary results, namely lemmas \ref{suben3} and \ref{pasa}.

Lemma \ref{suben3} will be a key step in proving there is uniformity in the advance to the right of the iterates of $\ell_0$. It tells us that the points that remain under iteration by $f$ in a strip $(\ell_i,\ell_{i+1})$, must go either upwards or downwards uniformly. We recall that in Section \ref{sec.curvasli} we proved that for each of the sets $\Theta_i$ from Theorem A, we have either $\rho(\Theta_i,f)\subset\{0\}\times\R^+$ or $\rho(\Theta_i,f)\subset\{0\}\times \R^-$. From now on we make the following assumption:

\begin{asuncion}   \label{asuncion.1a}
It holds $\rho(\Theta_0,f)\subset \{0\}\times \R^+$ and $\rho(\Theta_1,f)\subset \{0\}\times \R^-$.
\end{asuncion}

In the subsequent proofs, the complementary case $\rho(\Theta_0,f)\subset \{0\}\times\R^-$ will be obviously symmetric.


\begin{lema}[Uniformity Lemma] \label{suben3}
Given $m>0$ there exists $N\in\N$ such that, if:
\begin{itemize}
\item $i\in\{0,1\}$,
\item $n\in\Z$, $|n|\geq N$,
\item $x\in(\ell_i,\ell_{i+1})$ and $f^n(x)\in(\ell_i,\ell_{i+1})$,
\end{itemize} 
then:
\begin{itemize}
\item If $i=0$, then $f^n(x)_2-x_2>m$ if $n>0$ and $x_2-f^n(x)_2>m$ if $n<0$.
\item If $i=1$, then $x_2-f^n(x)_2>m$ if $n>0$ and $f^n(x)_2-x_2>m$ if $n<0$. 
\end{itemize}
\end{lema}

\begin{prueba}
We suppose the lemma does not hold, and we will find a contradiction. We treat the case $n>0$ and $i=0$, the cases $n<0$ and $i=1$ being completely analogous. We have then that there exists $m_0>0$, and sequences $\{x_j\}_j\subset (\ell_{0},\ell_{1})$, $\{s_j\}_j\subset\N$, such that:
\begin{itemize}
\item $s_j\ra\infty$ as $j\ra\infty$,
\item $f^{s_j}(x_j)\in(\ell_{0},\ell_{1})$ for all $j\in\N$, and
\item $f^{s_j}(x_j)_2-(x_j)_2<m_0$ for all $j\in\N$
\end{itemize}
Then,  
$$\limsup_j (f^{s_j}(x_j)_2-(x_j)_2)/s_j\leq 0,$$ 
and there is a subsequence of $\{x_j\}_j$, that we denote also $\{x_n\}_n$, such that 
$$\lim_j (f^{s_j}(x_j)-(x_j))/s_j = (0,a)$$ 
for some $a\leq 0$. Define the sequence of probability measures $\{\delta_j\}_j$ in $\T^2$ by
$$\delta_j = \frac{\delta_{\pi'(x_j)} + \delta_{\pi'(f(x_j))} + \cdots + \delta_{\pi'(f^{s_j-1}(x_j))} }{s_j},$$
and let $\delta$ be an accumulation point of $\{\delta_j\}_j$ in $\cl{M}_{\tl{f}}(\T^2)$. Then $\delta$ is $\tl{f}$-invariant, and 
$$\rho(\delta,f)=\int \phi \, d\delta= \lim_j \int \phi \, d(\delta_j) = \lim_j \frac{1}{s_j} (f^{s_j}(x_j)-x_j )=(0,a)$$
where $\phi:\T^2\ra\R^2$ is the displacement function defined in section \ref{sec.crmed}. Also, as $\txt{supp}(\delta)$ is $\tl{f}$-invariant and is contained in $[\tl{l}_{0},\tl{l}_{1}]$, $\txt{supp}(\delta)$ must be contained in $\Theta_{0}$ (because $\Theta_{0}$ is the maximal invariant set of $[\tl{l}_{0},\tl{l}_{1}]$). This means that $(0,a)\in\rho(\Theta_{0},f)$, and this is a contradiction, as $a\leq 0$ and by Assumption \ref{asuncion.1a} $\rho(\Theta_{0},f)\subset\{0\}\times\R^+$. This concludes the proof of the lemma.
\end{prueba}

As a corollary, we get that there is a maximum amount of displacement downwards for points that remain in $(\ell_0,\ell_{1})$ under iteration by $f$.

\begin{figure}[h] 
\begin{center} 
\includegraphics{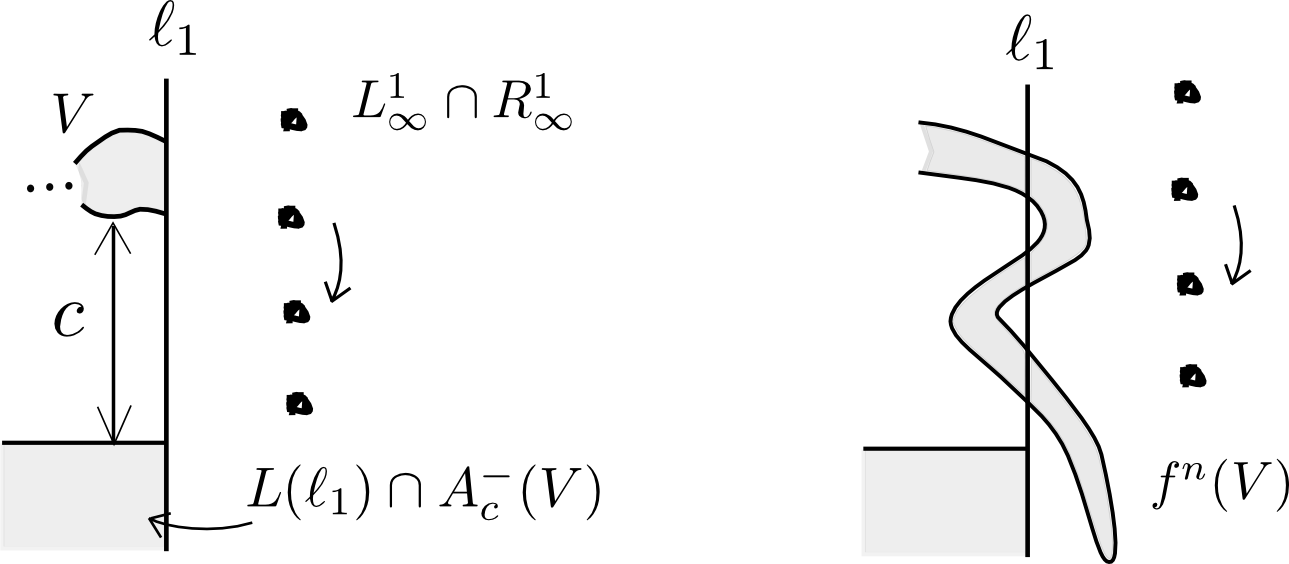}
\caption{Illustration of Corollary \ref{c}. Left: The sets $V$ and $A_c^-(V)$. Right: $f^n(V) \cap A_c^-(V) \subset R(\ell_{1})$ for all $n\geq 0$.}
\label{fig.lemac}  
\end{center}  
\end{figure}

\begin{corolario} \label{c}
There exists $c>0$ such that for any connected component $V$ of $X_0^c\cap(\ell_0,\ell_{1})$, we have that:
$$f^n(V)\cap  A_c^-(V) \subset R(\ell_{1}) \ \ \ \ \txt{for all $n\geq 0$},$$
where $A_c^-(V)=\{x\in\R^2\ : \ y_2-x_2>c \txt{ for all }y\in V\}$ (see Fig. \ref{fig.lemac}).  
\end{corolario}

\begin{prueba}
By Lemma \ref{suben3} there exists $N_0>0$ such that if:
\begin{itemize}
\item $n>N_0$, and
\item $x\in(\ell_0,\ell_{1})$ and $f^n(x)\in(\ell_0,\ell_{1})$,
\end{itemize}
then $f^n(x)_2 - x_2 > 0$. Let $V$ be any connected component of $X_0^c\cap(\ell_0,\ell_{1})$. Then we have that for any $x\in V$, either $f^{n}(x)\in R(\ell_{2})$ or $f^{n}(x)_2>x_2$, for any $n\geq N_0$. Making $c= N_0 \| f - \txt{Id} \|_0$, the lemma follows. 
\end{prueba}

Now we give our second preliminary result, which is a preliminary version of the Main Lemma \ref{N6}. 

\begin{lema} \label{pasa}
There exists $N_2>0$ such that for any curve $\gamma$ having good intersection with the anchors of the strip $[\ell_0,\ell_{1}]$ and with $\ell_{1}$, we have that $f^{N_2}(\gamma)$ intersects $\ell_{2}$.
\end{lema}

\begin{prueba}
Let $\gamma$ be a curve with good intersection with the anchors of $[\ell_0,\ell_{1}]$ and with $\ell_{1}$. By Lemma \ref{isotopia2} there is a constant $K_0>0$ such that if $C$ is a continuum contained in $(\ell_{1},\ell_{2})$ with $\txt{diam}_2\geq K_0$ then $f(C)\cap R(\ell_{2})\neq\empt$. 

By Lemma \ref{limitadas} there is $M_1>0$ such that for any connected component $V$ of $X_0^c\cap(\ell_0,\ell_{1})$, we have diam$_2(V)\leq M_1$. By Corollary \ref{c}, there is $c>0$ such that for any connected component $V$ of $X_0^c\cap(\ell_0,\ell_{1})$, and for any $n\geq 0$ we have that 
\begin{equation}    \label{art11}
f^n(V)\cap A_c^-(V)\subset R(\ell_{1}),
\end{equation}
where $A_c^-(V)=\{x\in\R^2\,:\, z_2-x_2 > c,\, \forall z\in V\}$. By Lemma \ref{suben3}, there exists $N_0>0$ such that:
\begin{itemize}
\item if $x\in [\ell_{1},\ell_{2}]$ then for all $n\geq N_0$ such that $f^n(x)\in [\ell_{1},\ell_{2}]$ we have $x_2-f^{n}(x)_2>M_1+c+K_0$, and 
\item if $y\in [\ell_0,\ell_{1}]$ then for all $n\geq N_0$ such that $f^n(y)\in [\ell_0,\ell_{1}]$ we have $f^n(y)_2-y_2> 0$.
\end{itemize}

The curve $\gamma$ has good intersection with the anchors of $[\ell_0,\ell_{1}]$ and with $\ell_{1}$, and then $\gamma$ contains an arc $\tl{\gamma}$ that realizes the good intersection (cf. Definition \ref{def.good2}). As $\tl{\gamma}(0)\in L_{\infty}^0$, we have $f^{N_0}(\tl{\gamma}(0))\in (\ell_0,\ell_{1})$, and by the choice of $N_0$,  
\begin{equation}    \label{k9}
f^{N_0}(\tl{\gamma}(0))_2>\tl{\gamma}(0)_2.
\end{equation}
We consider two cases. Either 
$$f^{N_0}(\tl{\gamma}(1))\in\ol{R}(\ell_{2}), \ \ \txt{or} \ \ f^{N_0}(\tl{\gamma}(1))\in L(\ell_{2}).$$ 
As $\tl{\gamma}(0)\in L_{\infty}^0$, we have that
\begin{equation} \label{ww1}
f^{N_0+n}(\tl{\gamma})\cap\ell_{2}\neq\empt \ \ \ \forall\,n\geq 0, \ \ \ \txt{if }f^{N_0}(\tl{\gamma}(1))\in\ol{R}(\ell_{2}).
\end{equation}

Now suppose that $f^{N_0}(\tl{\gamma}(1))\in L(\ell_{2})$. By the choice of $N_0$, and as $\tl{\gamma}(1)\in \ell_{1}$,
\begin{equation}   \label{k20}
\tl{\gamma}(1)_2-f^{N_0}(\tl{\gamma}(1))_2> M_1+c+K_0. 
\end{equation}   
Let $V_0\subset\R^2$ be the connected component of $X_0^c\cap[\ell_0,\ell_{1}]$ whose closure contains $\tl{\gamma}$. As $\txt{diam}_2(V_0)\leq M_1$, by (\ref{k20}) we have that
\begin{equation}    \label{k.1a}
z_2-f^{N_0}(\tl{\gamma}(1))_2>c+K_0 \ \ \ \ \forall\, z\in V_0.
\end{equation}
Then, by (\ref{k9}) and (\ref{k.1a}), and as $\tl{\gamma}(0)\in\ol{V}_0$ we have
\begin{equation}   \label{art12}
f^{N_0}(\tl{\gamma}(0))_2-f^{N_0}(\tl{\gamma}(1))_2 \geq c+K_0.
\end{equation}
By (\ref{art11}), (\ref{art12}), and as $\tl{\gamma}\subset \ol{V}_0$ we have that if $C_0$ is the connected component of $f^{N_0}(\tl{\gamma})\cap A_c^-(V_0)$ that contains $f^{N_0}(\tl{\gamma}(1))$, then $C_0\subset R(\ell_{1})$ and diam$_2(C_0)\geq K_0$. By the choice of the constant $K_0$, $f(C_0)\cap R(\ell_2)\neq\empt$, and then $f^{N_0+1}(\tl{\gamma})\cap R(\ell_{2})\neq\empt$. As $\tl{\gamma}(0)\in L(\ell_{1})$ we conclude that
\begin{equation} \label{ww2}
f^{N_0+1}(\tl{\gamma})\cap \ell_{2}\neq\emptyset, \ \ \ \txt{if }f^{N_0}(\tl{\gamma}(1))\in L(\ell_{2}).
\end{equation}
Combining (\ref{ww1}) and (\ref{ww2}), we have proved that
$$f^{N_0+1}(\tl{\gamma})\cap \ell_{2}\subset f^{N_0+1}(\gamma)\cap \ell_{2}\neq\emptyset.$$

Letting $N_2=N_0+1$ we have that for any curve $\gamma$ within the hypotheses of the lemma, $f^{N_2}(\gamma)\cap\ell_{2}\neq\empt$, and the lemma follows.
\end{prueba}

Now we are ready to prove the Main Lemma \ref{N6}.

\begin{prueba}[Proof of Main Lemma \ref{N6}]
Let $\gamma$ be a curve with good intersection with the anchors of $[\ell_0,\ell_{1}]$ and with $\ell_{1}$, and let $\tl{\gamma}\subset\gamma$ be an arc that realizes the good intersection (cf. Definition \ref{def.good2}).

Let $V$ be the connected component of $X_0^c\cap(\ell_0,\ell_{1})$ whose closure contains $\tl{\gamma}$. By Lemma \ref{limitadas} there exists a constant $M_1$ such that $V$ has diameter less than or equal to $M_1$.  By Corollary \ref{c} there is a constant $c$ such that
$$f^n(\tl{\gamma})\cap A_c^-(V) \subset R(\ell_{1}) \ \ \ \txt{ for all $n\geq 0$},$$
where $A_c^-(V) =\{x\in\R^2\, : \, y_2-x_2>c\txt{ for any }y\in V\}$.  

Let $p\in L_{\infty}^{1}\cap R_{\infty}^{1}$. Let $F$ be the anchor of $[\ell_{1},\ell_{2}]$ that contains $p$, and let $C$ be the connected component of $R_{\infty}^{1}\cap(\ell_{1},\ell_{2})$ that contains $p$. By Lemma \ref{propp2} there is $M_0>0$ such that $\txt{diam}_2(F)<M_0$ and $\txt{diam}_2(C)<M_0$, and then 
$$\txt{diam}_2(F\cup C)< 2M_0.$$ 
Let $s\in\Z$ be such that $T_2^s(F\cup C)\subset A_c^-(V)$ and $T_2^{s+1}(F\cup C)\cap A_c^-(V)^c\neq\empt$. Up to changing the point $p$ for a vertical integer translate of it, we may assume that $s=0$, and then
$$F\cup C\subset A_c^-(V) \ \ \ \ \txt{ and }\ \ \ \  T_2(F\cup C)\cap A_c^-(V)^c\neq\empt.$$
Therefore, for any $z\in C$,
\begin{equation}    \label{k11}
z+(0,2M_0+1)\in A_c^-(V)^c.
\end{equation}
By Lemma \ref{pasa} there is $N_2>0$ such that $f^{N_2}(\tl{\gamma})\cap\ell_{2}\neq\empt$. Let 
$$c_1=2M_0 + M_1 + c + N_2\|f-\txt{Id}\|_0 + 1.$$ 
By Lemma \ref{suben3}, there exists $N_0>0$ such that if $x$ and $f^{-N_0}(x)$ are contained in $(\ell_{1},\ell_{2})$ then $f^{-N_0}(x)_2-x_2> c_1$. In particular, 
\begin{equation}   \label{k10}
f^{-N_0}(z)_2-z_2>c_1 \ \ \ \txt{for any $z\in R_{\infty}^{1}\cap L(\ell_{2})$.}
\end{equation} 
As $\tl{\gamma}\subset \ol{V}$, if $y\in \ol{V}$ and $z\in f^{N_2}(\tl{\gamma})$, we have
\begin{equation}    \label{k12}
z_2-y_2 \leq M_1+N_2\left\| f-\txt{Id} \right\|_0.
\end{equation}
Then, by the definition of $c_1$ and $A_c^-(V)$, by (\ref{k11}), (\ref{k10}) and (\ref{k12}), we have that for any point $z$ in $C$, 
\begin{equation}    \label{k13}
f^{-N_0}(z)_2>y_2 \ \ \  \txt{for any $y\in f^{N_2}(\tl{\gamma})$}
\end{equation}
(see Fig. \ref{fig.leman3}). Remark that the constant $c_1$ is independent of $i$, and therefore the constant $N_0$ is also independent of $i$.

\begin{figure}[h] 
\begin{center} 
\includegraphics{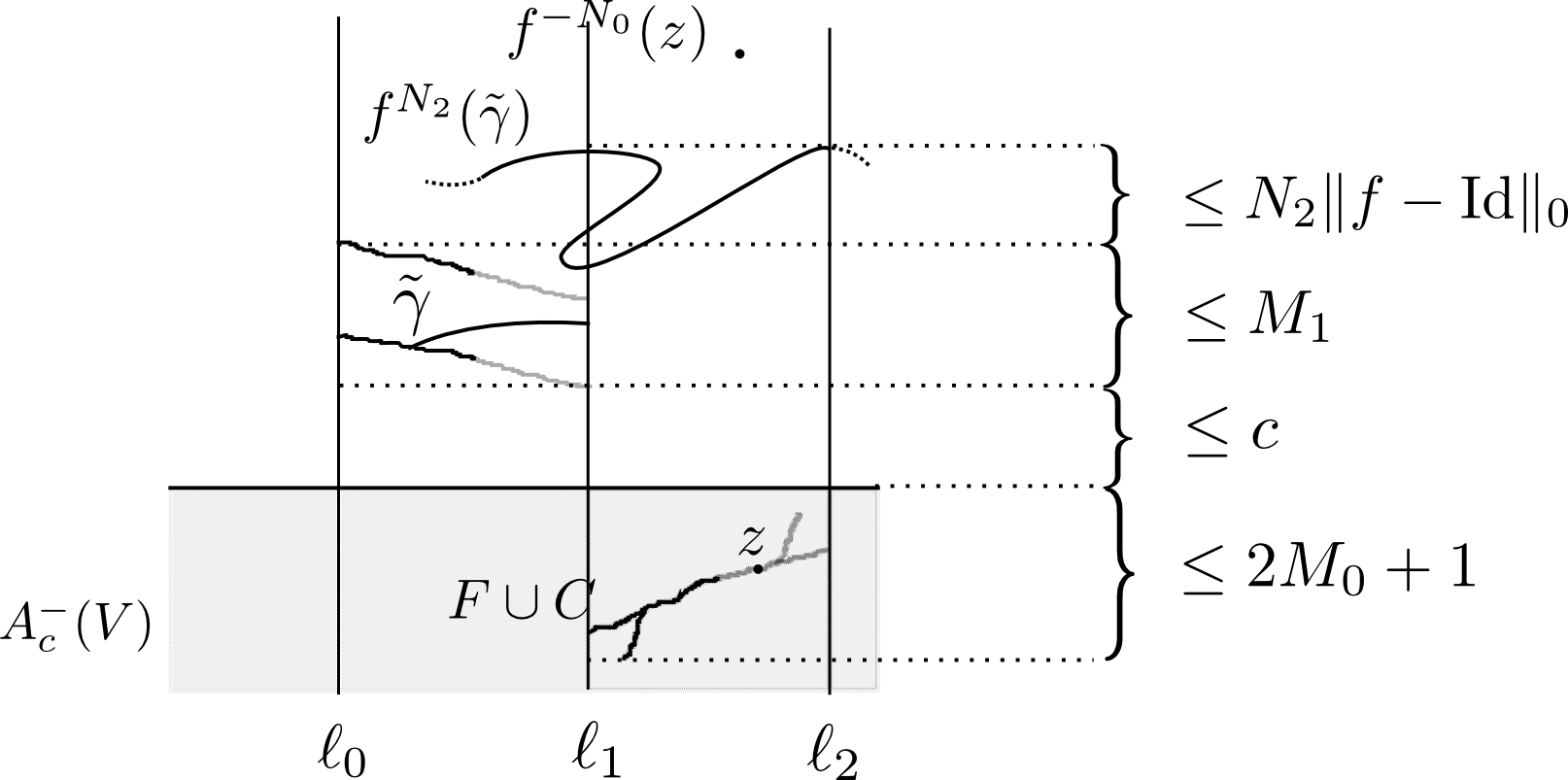}
\caption{If $z\in C$, then $f^{-N_0}(z)$ is above $f^{N_2}(\tl{\gamma})$.}
\label{fig.leman3}  
\end{center}  
\end{figure}

Now, let $\beta^1:(-\infty,0]\ra\R^2$ be a proper embedding such that
\begin{itemize}
\item $\beta^1(0)= f^{N_2}(\tl{\gamma}(0))\in L_{\infty}^0 \subset L(\ell_{1})$,
\item $\beta^1(t)\in L(f^{-N_0}(\ell_{1}))$ for all $t \leq 0$, and
\item $-\infty < \inf \, \{\beta^1(t)_2\, : \, t\leq 0\}$.
\end{itemize}

Let $\beta^2:[0,1]\ra\R^2$ be the arc contained in $f^{N_2}(\tl{\gamma})$ with endpoints $\beta^1(0)$ and $f^{N_2}(\tl{\gamma}(t_*))$, where $t_*=\min \{t \, : \, f^{N_2}(\tl{\gamma}(t))\in \ell_{2}\}$. Let $\beta^3:[0,\infty)\ra\R^2$ be the curve contained in $\ell_{2}$, starting in $\beta^2(1)$ and going upwards to infinity. Observe that the concatenation $\beta^1\beta^2\beta^3$ is a proper (not necessarily simple) map from $\R$ to $\R^2$ (see Fig. \ref{fig.encajeb}).

\begin{claim} \label{claim.mainlemma1}
The arc $f^{N_0}(\beta^2)\subset f^{N_0+N_2}(\tl{\gamma})$ has good intersection with the anchors of $[\ell_{1},\ell_{2}]$ and with $\ell_{2}$.
\end{claim}

\begin{prueba}
By (\ref{k13}) and by the definition of the arcs $\beta^j$, we have that $f^{-N_0}(C)\cap (\beta^1\cup\beta^2\cup\beta^3)=\empt$, and that for any $x\in f^{-N_0}(C)$, the straight vertical ray starting in $x$ and going upwards to infinity does not intersect the set $\beta^1\cup\beta^2\cup\beta^3$ either. Therefore the connected component $D$ of the complement of $\beta^1\cup\beta^2\cup\beta^3$ containing $f^{-N_0}(C)$ is unbounded from above (that is, $\sup\txt{pr}_2(D)=\infty$), and by the properties of the $\beta^j$ it is easy to see that it is bounded from below ($\inf\txt{pr}_2(D)>-\infty$). See Fig. \ref{fig.encajeb}.

\begin{figure}[h] 
\begin{center} 
\includegraphics{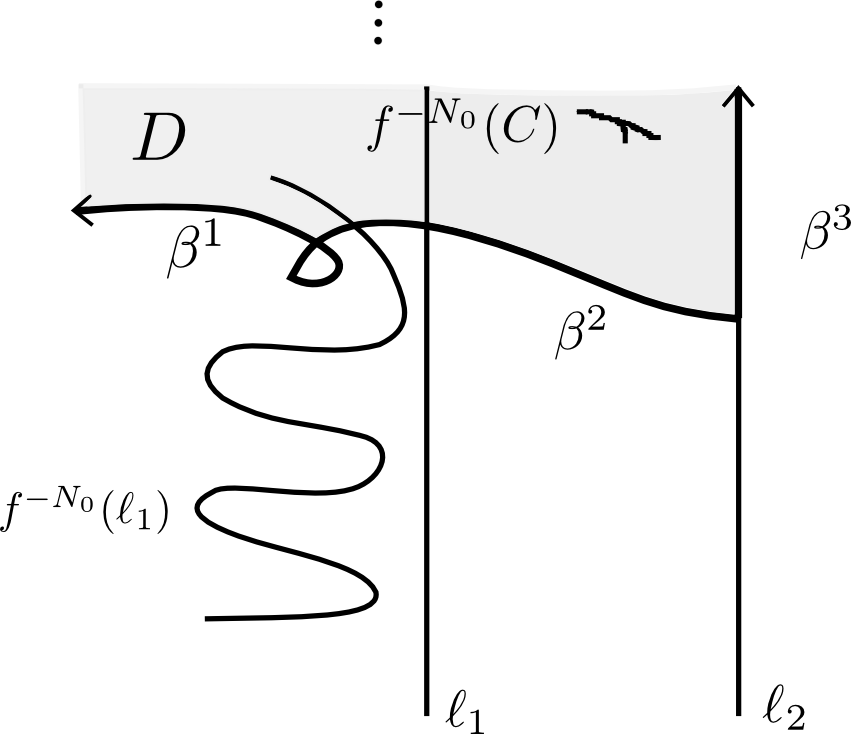}
\caption{The curves $\beta^j$ and the disk $D$.}
\label{fig.encajeb}  
\end{center}  
\end{figure}

As we have that $f^{-N_0}(C)\subset D$,
\begin{equation} \label {i6}
C\subset f^{N_0}(D).
\end{equation}
By the definition of the curves $\beta^1$ and $\beta^3$, 
$$(f^{N_0}(\beta^1)\cup f^{N_0}(\beta^3))\cap [\ell_{1},\ell_{2}]=\emptyset,$$   
and then, as $\pr D\subset \beta^1\cup\beta^2\cup\beta^3$, 
\begin{equation} \label{i20}
f^{N_0}(\pr D)\cap [\ell_{1},\ell_{2}] \subset f^{N_0}(\beta^2)\cap [\ell_{1},\ell_{2}].
\end{equation}
Since $f^{N_0}(\beta^2)\subset f^{N_0+N_2}(\tl{\gamma})$, by the definition of $A_c^-(V)\subset\R^2$ 
\begin{equation} \label{i5}
f^{N_0}(\beta^2)\cap A_c^-(V) \subset R(\ell_{1}).
\end{equation} 

Recall that $F$ and $C$ are the connected components of $L_{\infty}^{1}\cap[\ell_{1},\ell_{2}]$ and $R_{\infty}^{1}\cap(\ell_{1},\ell_{2})$, respectively, that contain the point $p$ defined above. By Lemma \ref{discos} the connected components of $L_{\infty}^{1}$ and $R_{\infty}^{1}$ that contain $p$ are unbounded to the left and to the right, resp. Therefore, the set $F\cup C$ separates the strip $(\ell_{1},\ell_{2})$, that is, $(\ell_{1},\ell_{2})\setminus (F\cup C)$ has more than one connected component. Moreover, since $F\cup C$ is bounded, the complement of $F\cup C$ in $(\ell_{1},\ell_{2})$ has exactly one connected component unbounded from below.

Let $U$ be such a component. Observe that $U$ is also bounded from above (since $F\cup C$ is bounded), and that $U\subset A_c^-(V)$. Let $w$ be the point with lowest second coordinate in $\ol{C}\cap \ell_{2}$. Then, $w\in\pr U$. By (\ref{i6}),  
$$w\in f^{N_0}(\ol{D}).$$ 
Since $f^{N_0}(D)$ is bounded from below (because $D$ is), if $w'$ is a point in $\ell_{2}$ with $w'_2 < \min\txt{pr}_2(f^{N_0}(\ol{D}))$, then 
$$w'\notin f^{N_0}(\ol{D}),$$
and also $w'\in\pr U$. Therefore $w$ and $w'$ belong to different connected components of $\ol{U}\setminus \pr f^{N_0}(D)$. As $\pr f^{N_0}(D)\cap \ol{U} \subset f^{N_0}(\beta^2)\cap\ol{U}$ (by (\ref{i20})), there is a connected component $\alpha$ of $f^{N_0}(\beta^2)\cap\ol{U}$ (which is a compact arc) such that $w$ and $w'$ belong to different connected components of $\ol{U}\setminus \alpha$. Note that as $\alpha\subset f^{N_0}(\beta^2)\subset f^{N_0+N_2}(\tl{\gamma})$, and as $\tl{\gamma}\subset \ol{L}(\ell_{1})$, by the definition of $R_{\infty}^{1}$ we have that
\begin{equation}    \label{i45}
\alpha \cap R_{\infty}^{1}=\empt.
\end{equation}

\begin{figure}[h] 
\begin{center} 
\includegraphics{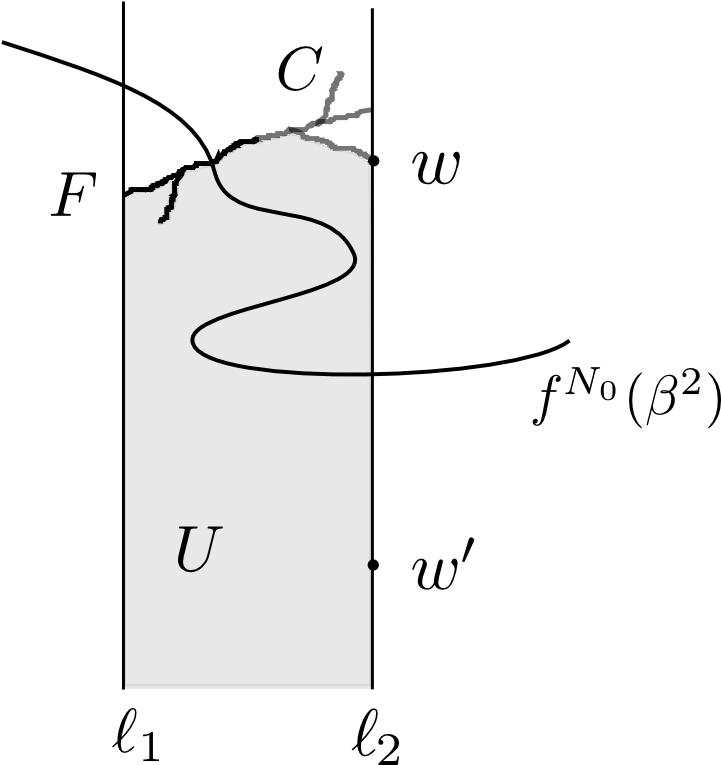}
\caption{The set U.}
\label{fig.U}  
\end{center}  
\end{figure}

\begin{claim}     \label{claim.100}
The arc $\alpha$ contains an arc $\alpha^1$ which has one endpoint in an anchor of $[\ell_{1},\ell_{2}]$ and the other in $\ell_{2}$.
\end{claim}

Before proving Claim \ref{claim.100} we use it to finish the proof of Claim \ref{claim.mainlemma1}. As the anchors of $[\ell_{1},\ell_{2}]$ are disjoint from $\ell_{2}$, there is a non-degenerate subarc $\tl{\alpha}$ of $\alpha^1$ which is minimal with respect to the property of having one endpoint in an anchor of $[\ell_{1},\ell_{2}]$ and the other in $\ell_{2}$. Note that $\tl{\alpha}\cap R_{\infty}^{1}=\empt$ by (\ref{i45}) and that $\tl{\alpha} \subset \ol{U}$. By this and by the minimality of $\tl{\alpha}$, the arc $\tl{\alpha}$ without its endpoints is contained in $(X_{1})^c\cap [\ell_{1},\ell_{2}]$.

As $\tl{\alpha}\subset\alpha^1\subset f^{N_0}(\beta^2)$, this means that $f^{N_0}(\beta^2)$ satisfies the definition of good intersection with the anchors of $[\ell_{1},\ell_{2}]$ and with $\ell_{2}$. This ends with the proof of Claim \ref{claim.mainlemma1}.
\end{prueba}

\begin{prueba}[Proof of Claim \ref{claim.100}]
First note that, as $\alpha$ is a connected component of $f^{N_0}(\beta^2)\cap \ol{U}$, and as the arc $f^{N_0}(\beta^2)$ has endpoints outside $\ol{U}$, we have that both endpoints of $\alpha$ must lie in $\pr U$.

Now, observe that 
$$\pr U\subset \ell_{1} \cup \ell_{2} \cup  C \cup F.$$ 
By (\ref{i5}), as $\alpha\subset f^{N_0}(\beta^2)$, and as $\alpha\subset \ol{U}\subset A_c^-(V)$, we have that $\alpha\cap\ell_{1}=\empt$. Also, by (\ref{i45}) and as $C\subset R_{\infty}^{1}$, $\alpha\cap C=\empt$. Therefore 
\begin{equation}    \label{i66}
\alpha\cap\pr U \subset F\cup\ell_{i+2}.
\end{equation}
We have then three possibilities:\\
\textbf{Case 1:} One endpoint of $\alpha$ lies in $F$ and the other in $\ell_{2}$.\\
\textbf{Case 2:} Both endpoints of $\alpha$ lie in $F$.\\
\textbf{Case 3:} Both endpoints of $\alpha$ lie in $\ell_{2}$.

If we are in Case 1, then the claim is proved, taking $\alpha^1=\alpha$. Suppose now we are in Case 2, and both endpoints of $\alpha$ lie in $F$. The arc $\alpha$ must intersect $\ell_{2}$; otherwise the whole subarc of $\ell_{2}$ from $w$ to $w'$ would belong to the same connected component in $\ol{U}\setminus\alpha$, which contradicts the definition of $\alpha$. Therefore $\alpha$ intersects $\ell_{2}$ and $\alpha$ contains an arc $\alpha^1$ with one endpoint in $F$ and the other in $\ell_{2}$, which proves the claim in Case 2. 

Finally, suppose we are in Case 3, and both endpoints of $\alpha$ lie in $\ell_{2}$. Then $\alpha$ must intersect $F$. To see this, note that if $\alpha\cap F=\empt$, then by (\ref{i66}) and by the definition of $w'$, $\alpha$ does not intersect the connected set $A:= \ol{C} \cup F\cup (\ell_{1} \cap \ol{U}) \cup B$, where $B$ is the set of points in $[\ell_{1},\ell_{2}]$ whose second coordinate is less than or equal to $w'_2$. As $A\subset \ol{U}\setminus\alpha$ and $w,w'\in A$, the points $w$ and $w'$ belong to the same connected component of $\ol{U}\setminus\alpha$, which contradicts the definition of $\alpha$. Therefore, the arc $\alpha$ intersects $F$ and it contains a subarc $\alpha^1$ with one endpoint in $F$ and the other in $\ell_{2}$, which proves the claim in Case 3. 

This finishes the proof of Claim \ref{claim.100}.
\end{prueba}

We now finish the proof of Main Lemma. By Claim \ref{claim.mainlemma1} we have that $f^{N_0+N_2}(\tl{\gamma})\subset f^{N_0+N_2}(\gamma)$ has good intersection with the anchors of $[\ell_{1},\ell_{2}]$ and with $\ell_{2}$. Setting $N_1=N_0+N_2$, the Main Lemma follows. 
\end{prueba}

\section{Proof of Addendum B.}               \label{sec.addendum}

We recall the statement of the addendum.

\begin{addb}
Let $\tl{f}$ and $f$ be as in Theorem A. Let $\{l_i\}_{i=0}^{r-1}$, and $\{l'_i\}_{i=0}^{r'-1}$ be two families of simple, closed, vertical and pairwise disjoint curves in $\T^2$, and suppose that both families satisfy the conslusion of Theorem A. Suppose also that the cardinalities $r$ and $r'$ are minimal with respect to this property.

Then, if we denote
$$\Theta_i = \bigcap_{n\in\Z} \tl{f}^n \left( [l_i,l_{i+1}] \right) \  \ \txt{for } i\in\Z/r\Z, \txt{ and}$$
$$\Theta'_i = \bigcap_{n\in\Z} \tl{f}^n \left( [l'_i,l'_{i+1}] \right)\ \  \txt{for } i\in\Z/r'\Z,$$
we have that $\{\Theta_i\}_{i=0}^{r'-1}= \{\Theta'_i\}_{i=0}^{r-1}$.
\end{addb}

We will work with the lift $f:\R^2\ra\R^2$, and we start with some notations. As in the proof of Theorem A (cf. Definition \ref{def.eles}) we fix consecutive lifts $\ell_i\subset\R^2$ of the curves $l_i$ for $0\leq i < r$, and consecutive lifts $\ell'_i\subset\R^2$ of the curves $l'_i$ for $0\leq i < r'$ (by consecutive we mean that $\ell_j\subset R(\ell_i)$ if $i < j$, and there is no other lift of $l_j$ contained in $(\ell_i,\ell_j)$). Also, for $0\leq i < r$ and $j\in\Z$ we define $\ell_{rj+i} = T_1^j(\ell_i)$, and for $0\leq i <r'$ and $j\in\Z$ we define $\ell'_{r'j+i}=T_1^j(\ell'_i)$. Finally, for $i\in\Z$ we let $\wt{\Theta}_i= \pi^{-1}(\Theta_i)\cap(\ell_i,\ell_{i+1})$, and $\wt{\Theta}'_i=\pi^{-1}(\Theta'_i)\cap(\ell'_i,\ell'_{i+1})$, where $\pi:\R^2\ra\T^2$ denotes the cannonical projection.  

The proof of the addendum will be by contradiction, and we begin with the following lemma. 

\begin{lema}    \label{lemma.addendum1}
Suppose that $\{\Theta_i\}_{i=0}^{r-1} \neq \{\Theta'_i\}_{i=0}^{r'-1}$. Then at least one of the two following holds:
\begin{itemize}
\item There is $i$ such that $\wt{\Theta}_i$ intersects both $\wt{\Theta}'_j$ and $\wt{\Theta}'_k$, for some $j\neq k$.
\item There is $i$ such that $\wt{\Theta}'_i$ intersects both $\wt{\Theta}_j$ and $\wt{\Theta}_k$, for some $j\neq k$.
\end{itemize}
\end{lema}

To prove the lemma, we start with the following weaker claim. 



\begin{claim}
If $\cup_{i\in\Z} \wt{\Theta}_i \neq \cup_{i\in\Z} \wt{\Theta}'_i$, then the conclusion of Lemma \ref{lemma.addendum1} holds.
\end{claim}

\begin{prueba} 
Without loss of generality we assume that $\cup_{i\in\Z} \wt{\Theta}_i \not\subset \cup_{i\in\Z} \wt{\Theta}'_i$. Let $i_0$ and $x$ be such that $x\in\wt{\Theta}_{i_0}$ and $x\in \R^2\setminus \cup_i \wt{\Theta}'_i$. 

Let $i_1$ be such that $x\in[\ell'_{i_1},\ell'_{i_1+1}]$, and without loss of generality assume $i_0=i_1=0$. By the choice of $x$, there are iterates of $x$ by $f$ which do not belong to $[\ell'_0,\ell'_1]$. Using the fact that the curves $\ell'_i$ are free for $f$, and that one of the sets $\Theta'_i$ is a vertical annular $\tl{f}$-invariant set, it is easy to see that there are $n_1$ and $i_2$ such that $f^n(x) \in (\ell'_{i_2}, \ell'_{i_2+1})$ for all $n\leq n_1$, and that there are $n_2$ and $i_3\neq i_2$ such that $f^n(x) \in (\ell'_{i_3},\ell'_{i_3+1})$ for all $n\geq n_2$. From this, it is also easy to see that $\wt{\Theta}_0$ intersects both $\wt{\Theta}'_{i_2}$ and $\wt{\Theta}'_{i_3}$. Therefore the conclusion of Lemma \ref{lemma.addendum1} holds, as desired.
\end{prueba}   

\begin{prueba}[Proof of Lemma \ref{lemma.addendum1}.]
Suppose that the conclusion of the lemma does not hold. Then, for all $i$, the set $\wt{\Theta}_i$ intersect at most one of the sets $\wt{\Theta}'_j$ and the set $\wt{\Theta}'_i$ intersects at most one of the sets $\wt{\Theta}_j$. Also, by last claim, $\cup_i \wt{\Theta}_i = \cup_i \wt{\Theta}'_i$. Therefore, for all $i$, $\wt{\Theta}_i\subset \wt{\Theta}'_j$ and $\wt{\Theta}'_i\subset \wt{\Theta}_k$ for some $j$ and $k$. This easily implies that $\{\Theta_i\}_{i=0}^{r-1} = \{\Theta'_i\}_{i=0}^{r'-1}$, and this proves the lemma.
\end{prueba}

We now proceed to the proof of the addendum. Suppose that $\{\Theta_i\}_{i=0}^{r-1} \neq \{\Theta'_i\}_{i=0}^{r'-1}$. Our objecitve from now on is to find a contradiction. 

As by hypothesis the cardinalities $r$ and $r'$ of the families $\{l_i\}$ and $\{l'_i\}$ are minimal in order that the conclusion of Theorem A holds, we must have that, for any $i$, if $\rho(\Theta_i,f)\subset \R^+$ then $\rho(\Theta_{i-1},f)\cup \rho(\Theta_{i+1},f) \subset \R^-$, and also for any $i$, if $\rho(\Theta_i,f)\subset \R^-$ then $\rho(\Theta_{i-1},f)\cup \rho(\Theta_{i+1},f) \subset \R^+$. The same statement also holds for the sets $\Theta'_i$.

By last lemma, we may assume without loss of generality that there is $i_0$ such that $\wt{\Theta}'_{i_0}$ intersects both $\wt{\Theta}_j$ and $\wt{\Theta}_k$, for some $j\neq k$. Also, without loss of generality, we assume that $i_0=j=0$, and that $\rho(\Theta'_0,f) \subset \R^+$. Therefore, $\rho(\Theta_0,f)\cup \rho(\Theta_k,f)\subset \R^+$. Observe that by the preceding paragraph, $k$ is even. 

\begin{figure}[h] 
\begin{center} 
\includegraphics{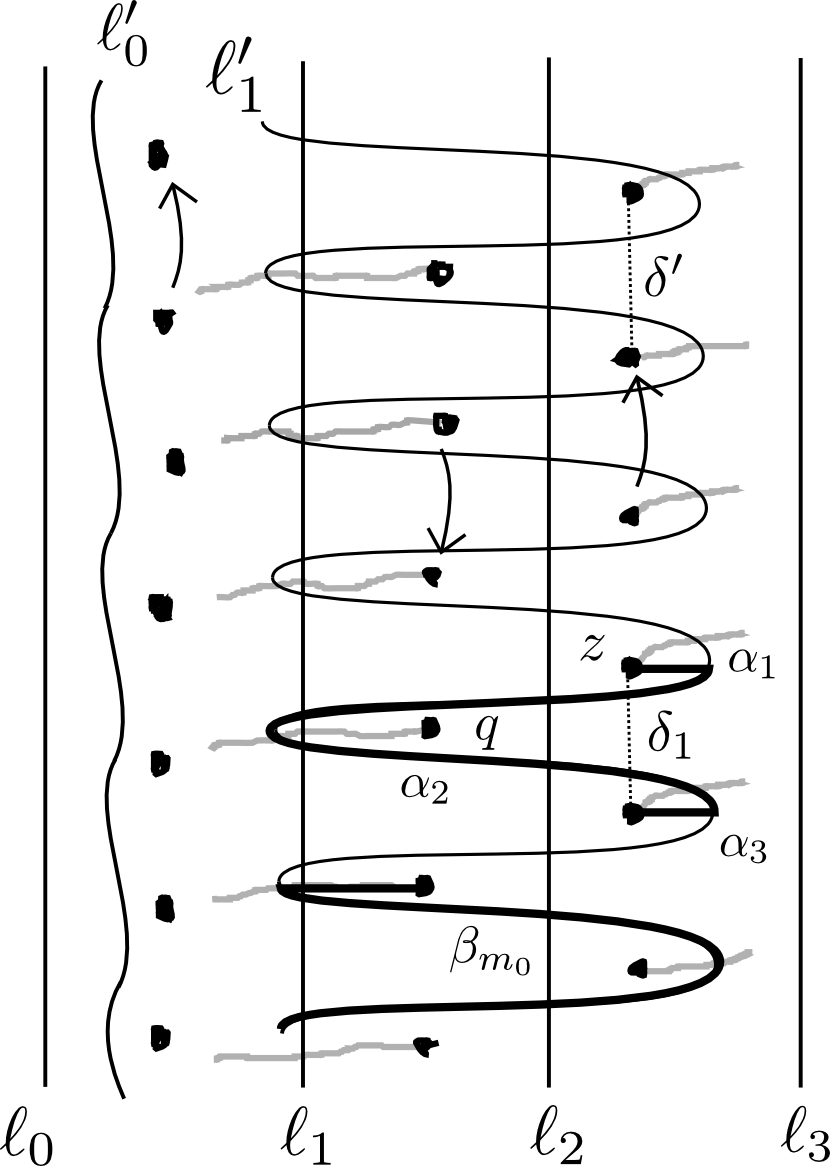}
\caption{Illustration for the assumption that $\wt{\Theta}'_0$ intersects both $\wt{\Theta}_0$ and $\wt{\Theta}_2$, and also that $\wt{\Theta}_1\cap R(\ell'_1)\neq\empt$. The sets $L_{\infty}^1$ and $R_{\infty}^2$ appear in gray, and the sets $\wt{\Theta}_0$, $\wt{\Theta}_1$, $\wt{\Theta}_2$ and $\wt{\Theta}'_0$ are represented by the thickened points.}
\label{fig.1a}  
\end{center}  
\end{figure}

We will work with the case $k=2$, the case $k\neq 2$ being analogous. We assume therefore that $\wt{\Theta}'_0$ intersects both $\wt{\Theta}_0$ and $\wt{\Theta}_2$. Observe that $\wt{\Theta}_1$ cannot be contained in $[\ell'_0,\ell'_1]$ (otherwise, by the invariance of $\wt{\Theta}_1$, we would have that $\wt{\Theta}_1\subset \wt{\Theta}'_0$, but $\rho(\Theta_1,f)\subset \R^-$ and $\rho(\Theta'_0,f)\subset \R^+$). Therefore $\wt{\Theta}_1\cap [\ell'_0,\ell'_1]^c\neq\empt$. We will assume that $\wt{\Theta}_1\cap R(\ell'_1)\neq\empt$ (see Fig. \ref{fig.1a}); in the subsequent arguments the case that $\wt{\Theta}_1\cap [\ell'_0,\ell'_1]^c\subset L(\ell'_0)$ will be symmetric.

Let $p\in \wt{\Theta}_1\cap R(\ell'_1)$. As the curve $\ell'_1$ is free for $f$, we have that, either $f^{n}(p)\in R(\ell'_1)$ for all $n \geq 0$, or $f^{-n}(p)\in R(\ell'_1)$ for all $n \geq 0$. Without loss of generality we will assume that we are in the first case, and then 
\begin{equation}  \label{zzz}
f^n(p)\in R(\ell'_1) \cap \wt{\Theta}_1  \ \ \ \ \forall \, n\geq 0.
\end{equation}
As the curves $\ell'_i$ are free for $f$, either $f(\ell'_1)\prec\ell'_1$ or $\ell'_1\prec f(\ell'_1)$. We will suppose that
\begin{equation}    \label{zxy}
f(\ell'_1) \prec \ell'_1,
\end{equation}
and the complementary case will follow from similar arguments (see Remark \ref{remark.add1}). We will further assume that for all $i$, the curve $\ell_i$ is a straight, vertical line.  

We will now define a closed curve $\Gamma$ such that Ind$(\Gamma,T_2^{n_0}(p))\neq 0$, for some $n_0\in\Z$. Let $z\in \wt{\Theta}_2\cap\wt{\Theta}'_0$. Let $\alpha_1:[0,1]\ra\R^2$ be the horizontal arc starting in $z$ and going to the right until the first moment of intersection with $\ell'_1$. Let $\alpha_2:[0,1]\ra\R^2$ be the arc contained in $\ell'_1$ starting in $\alpha_1(1)$ and ending in $\alpha_1(1) - (0,1)\in\ell'_1$. Define $\alpha_3:[0,1]\ra\R^2$ as the horizontal arc starting in $\alpha_2(1)$ and going to the left until $\alpha_1(0) - (0,1)$. Finally, let $\delta_1:[0,1]\ra\R^2$ be the straight arc from $\alpha_3(1)=\alpha_1(0)-(0,1)$ to $\alpha_1(0)$, and define $\Gamma=\alpha_1\alpha_2\alpha_3\delta_1$ (see Fig. \ref{fig.1a}). Recall that, for a point $x\in\R^2$, two arcs $\gamma_1$ and $\gamma_2$ in $\R^2\setminus\{x\}$ are said to be \textit{homotopic} Rel$(x)$ if there is a homotopy in $\R^2\setminus\{x\}$ from $\gamma_1$ to $\gamma_2$.

\begin{claim}            \label{claim2.a}
There is $n_0\in\Z$ such that Ind$(\Gamma,T_2^{n_0}(p))\neq 0$.
\end{claim}

\begin{prueba}
As $p\in L(\ell_2)$, as we assumed that $\ell_2$ is straight, and by the definition of $\Gamma$, we have that, for any $n\in\Z$, $\Gamma$ is homotopic Rel$(T_2^n(p))$ to the closed curve $\alpha_2\delta_2$, where $\delta_2:[0,1]\ra\R^2$ is the straight arc from $\alpha_2(0)$ to $\alpha_2(1)=\alpha_2(0)-(0,1)$. Therefore, to show that Ind$(\Gamma,T_2^{n_0}(p))\neq 0$ for some $n_0$, it suffices to show that Ind$(\alpha_2\delta_2,T_2^{n_0}(p))\neq 0$, or equivalently, that $\alpha_2$ is not homotopic with fixed endpoints Rel$(T_2^{n_0}(p))$ to $\delta_2$. 

To see this, suppose on the contrary that $\alpha_2$ is homotopic with fixed endpoints Rel$(T_2^n(p))$ to $\delta_2$ for all $n\in\Z$. Then, as $T_2^n(p)\in L(\ell_2)$ for all $n$, as $\delta_2\subset R(\ell_2)$, and as $\ell'_1=\cup_n T_2^n \alpha_2$, we get that $T_2^n(p)\in L(\ell'_1)$ for all $n$. However, by definition, $p$ belongs to $R(\ell'_1)$. This contradiction shows that there must be $n_0\in\Z$ such that $\alpha_2$ is not homotopic with fixed endpoints Rel$(T_2^{n_0}(p))$ to $\delta_2$, as we wanted. 
\end{prueba}

Continuing with the proof of the addendum, we now define auxiliary curves $\beta_m$. Let $n_0$ be as in last claim, let $q=T^{n_0}(p)$, and observe that by (\ref{zzz}), $f^n(q)\in R(\ell'_1) \cap \wt{\Theta}_1$ for all $n\geq 0$. 

For $m>0$, let $\beta_m=\beta_1\beta_2$, where $\beta_1:[0,1]\ra\R^2$ is the horizontal straight arc starting in $f^{m}(q)\in R(\ell'_1)\cap \wt{\Theta}_1$ and going to the left until the first point of intersection with $\ell'_1$, and where $\beta_2:[0,\infty)\ra\R^2$ is an injective curve contained in $\ell'_1$ such that $\beta_2(0)=\beta_1(1)$ and $\lim_{t\ra\infty}(\beta_2(t))_2 = -\infty$. Observe that as $\ell'_1$ is the lift of a vertical essential curve in $\T^2$, there is $k_0>0$ such that 
$$\sup \txt{pr}_2(\beta_m) - f^m(q)_2  < k_0 \ \ \ \txt{ for all } m>0.$$
On the other hand, as $q\in\wt{\Theta}_1$, and as $\rho(\Theta_1,f)\subset\R^-$, 
$$\lim_{n\ra\infty} f^n(q)_2 = -\infty.$$ 
Therefore, one easily verifies there is $m_1>0$ such that 
$$ \min \txt{pr}_2(\Gamma) > \max \txt{pr}_2(\beta_{m})   \ \ \ \txt{ for all } m\geq m_1.$$
Also, as $\alpha_1(0)=z$ and $\alpha_3(1)=z-(0,1)$ are contained in $\wt{\Theta}_2$, and as $\rho(\Theta_2,f)\subset\R^+$, there is $m_2>0$ such that 
$$\min \{ f^{m}(\{\alpha_1(0))_2, f^m(\alpha_3(1))_2 \} > \max \txt{pr}_2(\Gamma) \ \ \ \txt{ for all } m\geq m_2.$$ 
Let $m_0=\max\{m_1,m_2\}$.

\begin{claim}     \label{claim3.a}
Ind$(f^{m_0}(\Gamma), f^{m_0}(q))=0$.
\end{claim}

As Ind$(\Gamma,q)=\txt{Ind}(f^{m_0}(\Gamma), f^{m_0}(q))$, this claim contradicts the fact that Ind$(\Gamma,q)\neq 0$ given by Claim \ref{claim2.a}. This is the sought contradiction, and it will finish the proof of the addendum. We are therefore left with the proof of the claim.

\begin{prueba}[Proof of Claim \ref{claim3.a}.]
Recall that, as we assumed that the curves $\ell_i$ are straight, and as $\delta_1(0)=z$ and $\delta_1(1)=z-(0,1)$ are contained in $R(\ell_2)$, we have that $\delta_1\subset R(\ell_2)$. Then $\delta_1\cap L_{\infty}^1=\empt$, and by the invariance of $L_{\infty}^1$ we have 
\begin{equation}   \label{eq1a}
f^{m_0}(\delta_1)\cap L_{\infty}^1=\empt.
\end{equation}
Let $\delta'$ be the straight arc from $f^{m_0}(\delta_1(0))$ to $f^{m_0}(\delta_1(1))$. As $\ell_2$ is straight, and as $f^{m_0}(\delta_1(0)),f^{m_0}(\delta_1(1))\in \wt{\Theta}_2$, we have that 
\begin{equation}    \label{eq2a}
\delta'\subset R(\ell_2).
\end{equation}
As $\R^2\setminus L_{\infty}^1$ is simply connected (it is a nested union of simply connected sets), by (\ref{eq1a}), (\ref{eq2a}) and as $f^{m_0}(q)\in\wt{\Theta}_1\subset L_{\infty}^1$, we have that $f^{m_0}(\delta_1)$ is homotopic with fixed endpoints Rel$(f^{m_0}(q))$ to $\delta'$. This implies that
\begin{equation}    \label{u11}
\txt{Ind}(f^{m_0}(\Gamma),f^{m_0}(q))=\txt{Ind}(f^{m_0}(\alpha_1\alpha_2\alpha_3) \delta',f^{m_0}(q)).
\end{equation}
By the choice of $m_0$, 
$$\delta' \cap \beta_{m_0} = \empt.$$
Also, as $\alpha_1\alpha_2\alpha_3 \subset \ol{L}(\ell'_1)$, and as we assumed in (\ref{zxy}) that $f(\ell'_1)\prec \ell'_1$, we have $f^{m_0}(\alpha_1\alpha_2\alpha_3)\subset L(\ell'_1)$, and as $\beta_{m_0}\subset \ol{R}(\ell'_1)$, we get that 
$$f^{m_0}(\alpha_1\alpha_2\alpha_3)\cap\beta_{m_0}=\empt.$$ 
Therefore 
$$f^{m_0}(\alpha_1\alpha_2\alpha_3) \delta' \cap \beta_{m_0}=\empt.$$
As $\beta_{m_0}$ is unbounded and $f^{m_0}(q)\in\beta_{m_0}$, we have that $f^{m_0}(q)$ belongs to the unbounded component of $\R^2\minus  f^{m_0}(\alpha_1\alpha_2\alpha_3) \delta'$. This implies that 
$$\txt{Ind}(f^{m_0}(\alpha_1\alpha_2\alpha_3)\delta',f^{m_0}(q))=0,$$
and by (\ref{u11}), this yields Ind$(f^{m_0}(\Gamma),f^{m_0}(q))=0$. This proves our claim. 
\end{prueba}

\begin{remark}    \label{remark.add1}
We worked with the assumption (\ref{zxy}) that $f(\ell'_1)\prec \ell'_1$. The case $\ell'_1\prec f(\ell'_1)$ can be treated analogously. In a similar way as we constructed the arc $\alpha_1\alpha_2\alpha_3\subset \ol{L}(\ell'_1)$, one can construct an arc $\tl{\alpha}\subset \ol{R}(\ell'_1)$ with endpoints in $\wt{\Theta}_1$, and in a similar way as the curve $\beta_{m_0}\subset \ol{R}(\ell'_1)$ was constructed, one can construct a curve $\tl{\beta}\subset \ol{L}(\ell'_1)$ unbounded from above and with $\beta(0)\in\wt{\Theta}_2$. Then, all the arguments above work in a symmetric way to find a contradiction.
\end{remark}


\bibliographystyle{amsalpha}
\bibliography{bibtesis1}

\providecommand{\bysame}{\leavevmode\hbox to3em{\hrulefill}\thinspace}
\providecommand{\MR}{\relax\ifhmode\unskip\space\fi MR }
\providecommand{\MRhref}[2]{%
  \href{http://www.ams.org/mathscinet-getitem?mr=#1}{#2}
}
\providecommand{\href}[2]{#2}
\begin{thebibliography}{Han89}

\bibitem[Atk76]{atk}
G.~Atkinson, \emph{Recurrence of cocycles and random walks}, Journal of the
  London Mathematical Society \textbf{13} \textbf{no. 2} (1976), 486--488.

\bibitem[Bro12]{br}
L.~E.~J. Brouwer, \emph{Beweis des ebenen translationssatzes}, Math. Ann.
  \textbf{72} (1912), 37--54.

\bibitem[Cai51]{cairns}
S.~Cairns, \emph{An elementary proof of the {J}ordan-{S}choenflies theorem},
  Proc. of the Amer. Math. Soc. \textbf{91} \textbf{no. 2} (1951), 860--867.

\bibitem[Cal91]{lc3}
P.~Le Calvez, \emph{Propri\'et\'es dynamiques des diff\'eomorphisms de l'anneau
  et du tore}, 1991.

\bibitem[Cal04]{lc1}
\bysame, \emph{Une version feuillet\'ee du th\'eor\`eme de translation de
  {B}rouwer}, Comment. Math. Helv. \textbf{21} (2004), 229--259.

\bibitem[Cal05]{lc2}
\bysame, \emph{Une version feuillet\'ee \'equivariante du th\'eor\`eme de
  translation de {B}rouwer}, Publications Math\'ematiques de l'IH\'ES
  \textbf{102} \textbf{no. 1} (2005), 1--98.

\bibitem[Fay02]{fay1}
B.~Fayad, \emph{Weak mixing for reparametrized linear flows on the torus},
  Ergodic Theory and Dynamical Systems \textbf{22} (2002), 187--201.

\bibitem[Fra88]{f1}
J.~Franks, \emph{Recurrence and fixed points for surface homeomorphisms},
  Ergodic Theory \& Dynamical Systems \textbf{$8^*$} (1988), 99--107.

\bibitem[Fra89]{f2}
\bysame, \emph{Realizing rotation vectors for torus homeomorphisms},
  Transactions of the American Mathematical Society \textbf{311} \textbf{no. 1}
  (1989), 107--115.

\bibitem[Fra95]{f3}
\bysame, \emph{The rotation set and periodic points for torus homeomorphisms},
  Dynamical Systems \& Chaos (Aoki, Shiraiwa, and Takahashi, eds.), World
  Scientific, Singapore (1995), 41--48.

\bibitem[Hae62]{hf}
A.~Haefliger, \emph{Vari\'et\'es feuillet\'es}, Ann. Scuola Norm. Sup. Pisa
  \textbf{(3)} \textbf{16} (1962), 367--397.

\bibitem[Han89]{handel}
M.~Handel, \emph{Periodic point free homeomorphisms of $\mathbb{T}^2$}, Proc.
  of the Amer. Math. Soc. \textbf{107} \textbf{2} (1989), 511--515.

\bibitem[J\"09]{j1}
T.~J\"ager, \emph{Linearisation of conservative toral homeomorphisms}, Invent.
  Math. \textbf{176} \textbf{3} (2009), 601--616.

\bibitem[KK08]{kk}
A.~Kocksard and A.~Koropecki, \emph{Free curves and periodic points for torus
  homeomorphisms}, Ergodic Theory \& Dynamical Systems \textbf{28} (2008),
  1895--1915.

\bibitem[KK09]{kk2}
\bysame, \emph{A mixing-like property and inexistence of invariant foliations
  for minimal diffeomorphisms of the 2-torus}, Proc. of the Amer. Math. Soc.
  \textbf{137} \textbf{10} (2009), 3379--3386.

\bibitem[KT12]{korotal}
A.~Koropecki and Fabio~A. Tal, \emph{Area preserving irrotational
  diffeomorphisms of the torus with sublinear diffusion}, eprint
  arXiv:1206.2409 (2012).

\bibitem[Kur68]{kr}
K.~Kuratowski, \emph{Topology}, vol.~2, Academic Press, 1968.

\bibitem[LM91]{lm}
J.~Llibre and R.~S. Mackay, \emph{Rotation vectors and entropy for torus
  homeomorphisms isotopic to the identity}, Ergodic Theory \& Dynamical Sistems
  \textbf{(11)} (1991), 115--128.

\bibitem[MZ89]{mz}
M.~Misiurewicz and K.~Ziemian, \emph{Rotation set for maps of tori}, Journal of
  the London Mathematical Society \textbf{40} \textbf{no. 2} (1989), 490--506.

\bibitem[MZ91]{mz2}
\bysame, \emph{Rotation sets and ergodic measures for torus homeomorphisms},
  Fund. Math. \textbf{137} \textbf{1} (1991), 45--52.

\bibitem[Poi52]{po}
H.~Poincar\'e, \emph{Oeuvres compl\`etes}, {G}authier-{V}illars, Paris, 1952.

\end{thebibliography}
\footnotesize
UNIVERSIDADE DE S\~AO PAULO, INSTITUTO DE MATEM\'ATICA E ESTAT\'ISTICA, RUA DO MAT\~AO 1010, CIDADE UNIVERSIT\'ARIA, CEP 05508-090, BRASIL.\\
\ttfamily
email: pablod@ime.usp.br.

\end{document}